\documentclass[10pt,journal,compsoc]{IEEEtran}

\usepackage{multicol}

\usepackage{amssymb}
\usepackage{}
\usepackage{amsfonts}

\ifCLASSOPTIONcompsoc
  \usepackage[nocompress]{cite}
\else
  \usepackage{cite}
\fi
\usepackage{amsthm}
\usepackage{amsmath}

\usepackage{multirow}
\usepackage{color}

\usepackage{algorithm}
\usepackage{algorithmic}
\newtheorem{assumption}{Assumption}
\newtheorem{theorem}{Theorem}
\newtheorem{lemma}{Lemma}

\newtheorem{definition}{Definition}

\usepackage{hyperref}
\usepackage{cite}

\usepackage{graphicx}
\usepackage[tight,footnotesize]{subfigure}
\usepackage{slashbox}

\makeatletter

\newcommand{\Rmnum}[1]{\expandafter\@slowromancap\romannumeral #1@}
\makeatother

\ifCLASSINFOpdf
\else
\fi

\usepackage{hyperref}
\usepackage{cite}

\usepackage{diagbox}

\begin{document}
%
% paper title
% Titles are generally capitalized except for words such as a, an, and, as,
% at, but, by, for, in, nor, of, on, or, the, to and up, which are usually
% not capitalized unless they are the first or last word of the title.
% Linebreaks \\ can be used within to get better formatting as desired.
% Do not put math or special symbols in the title.
\title{Mini-batch Stochastic ADMMs for Nonconvex Nonsmooth Optimization}
%
%
% author names and IEEE memberships
% note positions of commas and nonbreaking spaces ( ~ ) LaTeX will not break
% a structure at a ~ so this keeps an author's name from being broken across
% two lines.
% use \thanks{} to gain access to the first footnote area
% a separate \thanks must be used for each paragraph as LaTeX2e's \thanks
% was not built to handle multiple paragraphs
%
%
%\IEEEcompsocitemizethanks is a special \thanks that produces the bulleted
% lists the Computer Society journals use for "first footnote" author
% affiliations. Use \IEEEcompsocthanksitem which works much like \item
% for each affiliation group. When not in compsoc mode,
% \IEEEcompsocitemizethanks becomes like \thanks and
% \IEEEcompsocthanksitem becomes a line break with idention. This
% facilitates dual compilation, although admittedly the differences in the
% desired content of \author between the different types of papers makes a
% one-size-fits-all approach a daunting prospect. For instance, compsoc
% journal papers have the author affiliations above the "Manuscript
% received ..."  text while in non-compsoc journals this is reversed. Sigh.

\author{Feihu Huang
        and Songcan Chen% <-this % stops a space
\IEEEcompsocitemizethanks{\IEEEcompsocthanksitem The authors are with the College of Computer Science and Technology,
Nanjing University of Aeronautics and Astronautics, Nanjing 210016, China. \protect
% note need leading \protect in front of \\ to get a newline within \thanks as
% \\ is fragile and will error, could use \hfil\break instead.
E-mail: \{huangfeihu, s.chen\}@nuaa.edu.cn}% <-this % stops a space
\thanks{}}

% note the % following the last \IEEEmembership and also \thanks -
% these prevent an unwanted space from occurring between the last author name
% and the end of the author line. i.e., if you had this:
%
% \author{....lastname \thanks{...} \thanks{...} }
%                     ^------------^------------^----Do not want these spaces!
%
% a space would be appended to the last name and could cause every name on that
% line to be shifted left slightly. This is one of those "LaTeX things". For
% instance, "\textbf{A} \textbf{B}" will typeset as "A B" not "AB". To get
% "AB" then you have to do: "\textbf{A}\textbf{B}"
% \thanks is no different in this regard, so shield the last } of each \thanks
% that ends a line with a % and do not let a space in before the next \thanks.
% Spaces after \IEEEmembership other than the last one are OK (and needed) as
% you are supposed to have spaces between the names. For what it is worth,
% this is a minor point as most people would not even notice if the said evil
% space somehow managed to creep in.

% The paper headers
\markboth{IEEE Transactions on Pattern Analysis and Machine Intelligence}%
{Shell \MakeLowercase{\textit{et al.}}: Bare Advanced Demo of IEEEtran.cls for IEEE Computer Society Journals}
% The only time the second header will appear is for the odd numbered pages
% after the title page when using the twoside option.
%
% *** Note that you probably will NOT want to include the author's ***
% *** name in the headers of peer review papers.                   ***
% You can use \ifCLASSOPTIONpeerreview for conditional compilation here if
% you desire.

% The publisher's ID mark at the bottom of the page is less important with
% Computer Society journal papers as those publications place the marks
% outside of the main text columns and, therefore, unlike regular IEEE
% journals, the available text space is not reduced by their presence.
% If you want to put a publisher's ID mark on the page you can do it like
% this:
%\IEEEpubid{0000--0000/00\$00.00~\copyright~2015 IEEE}
% or like this to get the Computer Society new two part style.
%\IEEEpubid{\makebox[\columnwidth]{\hfill 0000--0000/00/\$00.00~\copyright~2015 IEEE}%
%\hspace{\columnsep}\makebox[\columnwidth]{Published by the IEEE Computer Society\hfill}}
% Remember, if you use this you must call \IEEEpubidadjcol in the second
% column for its text to clear the IEEEpubid mark (Computer Society journal
% papers don't need this extra clearance.)

% use for special paper notices
%\IEEEspecialpapernotice{(Invited Paper)}

% for Computer Society papers, we must declare the abstract and index terms
% PRIOR to the title within the \IEEEtitleabstractindextext IEEEtran
% command as these need to go into the title area created by \maketitle.
% As a general rule, do not put math, special symbols or citations
% in the abstract or keywords.
\IEEEtitleabstractindextext{%
\begin{abstract}
With the large rising of complex data, the nonconvex models such as nonconvex loss function and nonconvex regularizer
are widely used in machine learning and pattern recognition.
In this paper, we propose a class of mini-batch stochastic ADMMs (alternating direction method of multipliers)
for solving large-scale nonconvex nonsmooth problems. We prove that, given an appropriate mini-batch size,
the mini-batch stochastic ADMM without variance reduction (VR) technique is convergent and reaches a convergence rate of $O(1/T)$
to obtain a stationary point of the nonconvex optimization,
where $T$ denotes the number of iterations.
Moreover, we extend the mini-batch stochastic gradient method to
both the nonconvex SVRG-ADMM and SAGA-ADMM proposed in our initial manuscript \cite{huang2016stochastic},
and prove these mini-batch stochastic ADMMs also reaches the convergence rate of $O(1/T)$
without condition on the mini-batch size.
In particular, we provide a specific parameter selection for step size $\eta$ of stochastic gradients and
penalty parameter $\rho$ of augmented Lagrangian function.
Finally, extensive experimental results on both simulated and real-world data demonstrate the effectiveness of the proposed algorithms.
\end{abstract}

% Note that keywords are not normally used for peerreview papers.
}

% make the title area
\maketitle

% To allow for easy dual compilation without having to reenter the
% abstract/keywords data, the \IEEEtitleabstractindextext text will
% not be used in maketitle, but will appear (i.e., to be "transported")
% here as \IEEEdisplaynontitleabstractindextext when compsoc mode
% is not selected <OR> if conference mode is selected - because compsoc
% conference papers position the abstract like regular (non-compsoc)
% papers do!
\IEEEdisplaynontitleabstractindextext
% \IEEEdisplaynontitleabstractindextext has no effect when using
% compsoc under a non-conference mode.

% For peer review papers, you can put extra information on the cover
% page as needed:
% \ifCLASSOPTIONpeerreview
% \begin{center} \bfseries EDICS Category: 3-BBND \end{center}
% \fi
%
% For peerreview papers, this IEEEtran command inserts a page break and
% creates the second title. It will be ignored for other modes.
\IEEEpeerreviewmaketitle

\ifCLASSOPTIONcompsoc
\IEEEraisesectionheading{\section{Introduction}\label{sec:introduction}}
\else
\section{Introduction}
\label{sec:introduction}
\fi
% Computer Society journal (but not conference!) papers do something unusual
% with the very first section heading (almost always called "Introduction").
% They place it ABOVE the main text! IEEEtran.cls does not automatically do
% this for you, but you can achieve this effect with the provided
% \IEEEraisesectionheading{} command. Note the need to keep any \label that
% is to refer to the section immediately after \section in the above as
% \IEEEraisesectionheading puts \section within a raised box.

% The very first letter is a 2 line initial drop letter followed
% by the rest of the first word in caps (small caps for compsoc).
%
% form to use if the first word consists of a single letter:
% \IEEEPARstart{A}{demo} file is ....
%
% form to use if you need the single drop letter followed by
% normal text (unknown if ever used by the IEEE):
% \IEEEPARstart{A}{}demo file is ....
%
% Some journals put the first two words in caps:
% \IEEEPARstart{T}{his demo} file is ....
%
% Here we have the typical use of a "T" for an initial drop letter
% and "HIS" in caps to complete the first word.

\IEEEPARstart {S}{t}ochastic optimization \cite{bottou2004stochastic}  is a class of powerful optimization tool for solving large-scale problems
in machine learning, pattern recognition and computer vision.
For example, stochastic gradient descent (SGD \cite{bottou2004stochastic}) is an efficient method
for solving the following optimization problem, which is a fundamental to machine learning,
\begin{align}
 \min_{x\in R^d} f(x) + g(x)
\end{align}
where $f(x)=\frac{1}{n}\sum_{i=1}^n f_i(x)$ denotes the loss function, and $g(x)$ denotes the regularization function.
The problem (1) includes many useful models such as support vector machine (SVM), logistic regression and neural network.
When sample size $n$ is large, even the first-order methods become computationally burdensome due to their per-iteration
complexity of $O(nd)$. While SGD only computes gradient of one sample instead of
all samples in each iteration, thus, it has only per-iteration complexity of $O(d)$.
Despite its scalability, due to the existence of variance in stochastic process, the
stochastic gradient is much noisier than the batch gradient.
Thus, the step size has to be decreased gradually as stochastic
learning proceeds, leading to slower convergence than the batch method.
Recently, a number of accelerated algorithms have successfully been proposed to reduce this variance.
For example, stochastic average gradient (SAG \cite{roux2012stochastic}) obtains a fast convergence rate
by incorporating the old gradients estimated in the previous iterations.
Stochastic dual coordinate ascent (SDCA \cite{shalev2013stochastic}) performs the stochastic coordinate ascent on the dual problems
to obtain also a fast convergence rate.
Moreover, an accelerated randomized proximal coordinate gradient (APCG \cite{lin2015accelerated}) method
accelerates the SDCA method by using Nesterov's accelerated method \cite{nesterov2004introductory}.
However, these fast methods require much space to store old gradients or dual variables.
Thus, stochastic variance reduced gradient (SVRG \cite{johnson2013accelerating,xiao2014proximal}) methods are proposed,
and enjoy a fast convergence rate with no extra space to store the intermediate gradients or dual variables.
Moreover,\cite{defazio2014saga} proposes the SAGA method, which extends the SAG method and enjoys better theoretical convergence rate
than both SAG and SVRG. Recently, \cite{nitanda2014stochastic} presents an accelerated SVRG
by using the Nesterov's acceleration technique \cite{nesterov2004introductory}.
Moreover, \cite{allen2016katyusha} proposes a novel momentum accelerated SVRG method (Katyusha) via using the strongly
convex parameter, which reaches a faster convergence rate.
In addition, \cite{zhang2017sparse} specially proposes a class of stochastic composite optimization methods for sparse learning, when
$g(\cdot)$ is a sparsity-inducing regularizer such as $\ell_1$-norm and nuclear norm.

\begin{table*}
  \centering
  \caption{Summary of the existing stochastic ADMMs for the nonconvex optimization.
  \checkmark denotes that the proposed methods can optimize the corresponding nonconvex problems.}\label{tab:1}

\begin{tabular}{|c|c|c|}
  \hline
  % after \\: \hline or \cline{col1-col2} \cline{col3-col4} ...
  \diagbox{Methods}{Convergence rate}{Problems} & $\min_x \frac{1}{n}\sum_{i=1}^n f_i(x)+g(x) $ & $\min_{x} \frac{1}{n}\sum_{i=1}^n f_i(x)+g(Ax)$\\ \hline
  Nonconvex incremental ADMM \cite{hong2014distributed} & \checkmark, Unknown &  \\ \hline
  NESTT \cite{hajinezhad2016nestt} & \checkmark, $O(1/T)$ &  \\ \hline
  Nonconvex mini-batch stochastic ADMM (ours) & \checkmark, $O(1/T)$ & \checkmark, $O(1/T)$ \\ \hline
  Nonconvex SVRG-ADMM (ours and \cite{zheng2016stochastic}) &\checkmark, $O(1/T)$ & \checkmark, $O(1/T)$ \\  \hline
  Nonconvex SAGA-ADMM (ours)& \checkmark, $O(1/T)$ & \checkmark, $O(1/T)$ \\
  \hline
\end{tabular}

\end{table*}

Though the above methods can effectively solve many problems in machine learning,
they are still difficultly to be competent for some complicated problems with
the \emph{nonseparable} and \emph{nonsmooth }regularization function as follows
\begin{align}
 \min_{x\in R^d} f(x) + g(Ax)
\end{align}
where $A\in R^{d\times p}$ is a given matrix, $f(x)=\frac{1}{n}\sum_{i=1}^n f_i(x)$ denotes the loss function,
and $g(x)$ denotes the regularization function. With regard to $g(\cdot)$, we are interested in a sparsity-inducing regularization functions,
e.g. $\ell_1$-norm and nuclear norm. The problem (2) includes the graph-guided fuzed Lasso \cite{kim2009multivariate},
the overlapping group Lasso\cite{yuan2013efficient}, and generalized Lasso \cite{tibshirani2011solution}.
It is well known that the alternating direction method
of multipliers (ADMM \cite{gabay1976dual,boyd2011distributed,lu2018unified}) is an efficient optimization method for the problem (2).
Specifically, we can use auxiliary variable $y=Ax$ to make the problem (2) be suitable for the general ADMM form.
When sample size $n$ is large, due to the need of computing the empirical risk loss function on all training samples
at each iteration, the offline or batch ADMM is unsuitable for large-scale learning problems.
Thus, the online and stochastic versions of ADMM \cite{wang2012online,suzuki2013dual,ouyang2013stochastic} have been successfully developed
for the large-scale problems.
Due to the existence of variance in the stochastic process, these stochastic ADMMs also suffer from the slow convergence rate.
Recently, some accelerated stochastic ADMMs are effectively proposed to reduce this variance.
For example, SAG-ADMM \cite{zhong2014fast} is proposed by additionally using the previous estimated gradients.
An accelerated stochastic ADMM \cite{azadi2014towards} is proposed
by using Nesterov's accelerated method \cite{nesterov2004introductory}.
SDCA-ADMM \cite{suzuki2014stochastic} obtains linearly convergence rate for the strong problem
by solving its dual problem.
SCAS-ADMM \cite{zhao2015scalable} and SVRG-ADMM \cite{zheng2016fast} are developed,
and reach the fast convergence rate with no extra space for the previous gradients or dual variables.
Moreover, \cite{liu2017accelerated} proposes an accelerated SVRG-ADMM
by using the momentum accelerated technique.
More recently, \cite{fang2017faster} proposes a fast stochastic ADMM,
which achieves a \emph{non-ergodic} convergence rate of $O(1/T)$ for the convex problem.
In addition, an adaptive stochastic ADMM \cite{zhao2015adaptive} is proposed by using the \emph{adaptive gradients}.
Due to that the penalty parameter in ADMM can affect convergence \cite{nishihara2015general},
another adaptive stochastic ADMM \cite{xu2017admm} is proposed by using the \emph{adaptive penalty parameters}.

So far, the above study on stochastic optimization methods relies heavily on strongly convex or convex problems.
However, there exist many useful nonconvex models in machine learning such as nonconvex empirical risk minimization models\cite{aravkin2016smart}
and deep learning \cite{lecun2015deep}.
Thus, the study of nonconvex optimization methods is much needed.
Recently, some works focus on studying the stochastic gradient methods for the large-scale nonconvex optimizations.
For example, \cite{ghadimi2016accelerated,ghadimi2016mini} have established the iteration complexity of $O(1/\epsilon^2)$ for the SGD
to obtain an $\epsilon$-stationary solution of the nonconvex problems.
\cite{allen2016variance,reddi2016stochastic,reddi2016fast1} have proved that the variance reduced stochastic gradient methods such as
the nonconvex SVRG and SAGA
reach the iteration complexity of $O(1/\epsilon)$.
At the same time, \cite{reddi2016fast2} has proved that the variance reduced stochastic gradient methods
also reach the iteration complexity of $O(1/\epsilon)$ for
the \emph{nonconvex nonsmooth} composite problems.
More recently, \cite{allen2017natasha} propose a faster nonconvex stochastic optimization method (Natasha)
via using the strongly non-convex parameter.
\cite{paquette2017catalyst} proposes a faster gradient-based nonconvex optimization
by using \emph{catalyst} approach in \cite{lin2015universal}.

Similarly, the above nonconvex methods are difficult to be competent
to some complicated nonconvex problems, such as nonconvex graph-guided regularization risk loss minimizations \cite{huang2016stochastic} and
tensor decomposition \cite{jiang2016structured}.
Recently, some works \cite{wang2015convergence,yang2015alternating,wang2015global,hong2016convergence,jiang2016structured}
have begun to study the ADMM method for the nonconvex optimization,
but they only focus on studying the deterministic ADMMs for the nonconvex optimization.
Due to the need of computing the empirical loss function on
all the training examples at each iteration, these nonconvex ADMMs are not yet well competent to the large-scale learning problems.
Recently, \cite{hong2014distributed} has proposed a distributed, asynchronous and incremental algorithm based on the ADMM method
for the large-scale nonconvex problems, but this method is difficult for the nonconvex problem (2)
with the \emph{nonseparable} and \emph{nonsmooth} regularizers such as graph-guided fused lasso and overlapping group lasso.
A nonconvex primal dual splitting (NESTT  \cite{hajinezhad2016nestt}) method is proposed for the distributed and stochastic optimization,
but it is also difficult for the nonconvex problem (2).
More recently, our initial manuscript \cite{huang2016stochastic} proposes the stochastic ADMMs with variance reduction (e.g., nonconvex SVRG-ADMM and
nonconvex SAGA-ADMM) for optimizing these nonconvex problems with some complicated structure regularizers
such as graph-guided fuzed Lasso, overlapping group Lasso, sparse plus low-rank penalties.
In addition, our initial manuscript \cite{huang2016stochastic} and \emph{Zheng and Kwok}'s paper \cite{zheng2016stochastic} \emph{simultaneously} propose
the nonconvex SVRG-ADMM method\footnote{ The first version of our manuscript\cite{huang2016stochastic}(https://arxiv.org/abs/1610.02758v1) proposes both \textbf{non-convex }SVRG-ADMM and SAGA-ADMM, which is online available in {\color{red}{\textbf{Oct. 10, 2016} }}.  \textbf{The first version} of \cite{zheng2016stochastic} ({\color{red}{https://arxiv.org/abs/1604.07070v1}}) only proposes the convex SVRG-ADMM, which is online available in Apr. 24, 2016 and named as '\emph{Fast-and-Light Stochastic ADMM}'. While, \textbf{the second version} of \cite{zheng2016stochastic} ({\color{red}{https://arxiv.org/abs/1604.07070v2}}) adds the \textbf{non-convex} SVRG-ADMM, which is online available in {\color{red}{\textbf{Oct. 12, 2016}}} and
\textbf{renamed} as '\emph{Stochastic Variance-Reduced ADMM}'. }.
At present, to our knowledge, there still exist two important problems needing to be addressed:
\begin{itemize}
  \item[1)] \emph{Whether the general stochastic ADMM without VR technique is convergent for the nonconvex optimization?}
  \item[2)]\emph{What is convergence rate of the general stochastic ADMM for the nonconvex optimization, if convergent?}
\end{itemize}

In the paper, we provide the positive answers to them by developing
a class of mini-batch stochastic ADMMs for the nonconvex optimization.
Specifically, we study the mini-batch stochastic ADMMs for optimizing the \emph{nonconvex nonsmooth} problem below:
\begin{align}
 \min_{x,y}  & \ \frac{1}{n}\sum_{i=1}^n f_i(x) + g(y)   \label{eq:1} \\
 \mbox{s.t.} & \ Ax + By =c, \nonumber
\end{align}
where $x\in R^d,y \in R^p$, $f(x)=\frac{1}{n}\sum_{i=1}^n f_i(x)$, each $f_i(x)$ is a \emph{nonconvex} and smooth loss function,
$g(y)$ is \emph{nonsmooth} and possibly nonconvex,
and $A\in R^{q\times d}$, $B\in R^{q\times p}$ and $c \in R^q$ denote the given matrices and vector, respectively.
The problem \eqref{eq:1} is inspired by the structural risk minimization in machine learning \cite{vapnik2013nature}.
In summary, our main contributions are four-fold as follows:
\begin{itemize}
\item[1)] We propose the mini-batch stochastic ADMM for the nonconvex nonsmooth optimization.
          Moreover, we prove that, given an appropriate mini-batch size, the mini-batch stochastic ADMM
          reaches a fast convergence rate of $O(1/T)$ to obtain a stationary point.
\item[2)] We extend the mini-batch stochastic gradient method to both the nonconvex SVRG-ADMM and SAGA-ADMM,
          proposed in our initial manuscript \cite{huang2016stochastic}.
          Moreover, we prove that these stochastic ADMMs also reach a convergence rate of $O(1/T)$ without condition on the mini-batch size.
\item[3)] We provide a specific parameter selection for step size $\eta$ of stochastic gradients and
          penalty parameter $\rho$ of the augmented Lagrangian function.
\item[4)] Some numerical experiments demonstrate the effectiveness of the proposed algorithms.
\end{itemize}
In addition, Table \ref{tab:1} shows the convergence rate summary of the stochastic/incremental ADMMs for optimizing the nonconvex problems.

\subsection{Notations}
$\|\cdot\|$ denotes the Euclidean norm of a vector or the spectral norm of a matrix.
$I_p$ denotes an $p$-dimensional identity matrix.
$H\succ0$ denotes a positive definite matrix $H$, and $\|x\|^2_H = x^THx$.
Let $A^+$ denote the generalized inverse of matrix $A$.
$\phi_{\min}^A$ denotes the smallest eigenvalues of matrix $A^T A$.
$\phi_{\max}^H$ and $\phi_{\min}^H$ denotes the largest and smallest eigenvalues of positive matrix $H$, respectively.
The other notations used in this paper is summarized as follows:
$\tilde{L}=L+1$, $\phi^H = (\phi^H_{\min})^2+20(\phi^H_{\max})^2$,  $\zeta = \frac{5(L^2\eta^2 + (\phi_{\max}^H)^2)}{\phi_{\min}^A\eta^2}$,
and $\zeta_1 = \frac{ 5(\phi_{\max}^H)^2}{\phi_{\min}^A\eta^2}$.

% needed in second column of first page if using \IEEEpubid
%\IEEEpubidadjcol

\section{Nonconvex Mini-batch Stochastic ADMM without VR}
In this section, we propose a mini-batch stochastic ADMM to optimize the nonconvex problem \eqref{eq:1}.
Moreover, we study convergence of the mini-batch stochastic ADMM.
In particular, we prove that, given an appropriate mini-batch size,
it reaches the convergence rate of $O(1/T)$.

First, we review the deterministic ADMM for solving the problem \eqref{eq:1}.
The augmented Lagrangian function of \eqref{eq:1} is defined as follows:
\begin{align}
 \mathcal {L}_{\rho}(x,y,\lambda) = & f(x) + g(y) - \langle\lambda, Ax+By-c\rangle \nonumber \\
 & + \frac{\rho}{2} \|Ax+By-c\|^2,
\end{align}
where $\lambda$ is the Lagrange multiplier, and $\rho$ is the penalty parameter.
At $t$-th iteration, the ADMM executes the update:
\begin{align}
& y_{t+1} = \arg\min_y \mathcal {L}_{\rho}(x_t,y,\lambda_t),  \label{eq:3}\\
& x_{t+1} = \arg\min_x \mathcal {L}_{\rho}(x,y_{t+1},\lambda_t), \\
& \lambda_{t+1} = \lambda_t - \rho(Ax_{t+1}+By_{t+1}-c).  \label{eq:5}
\end{align}

Next, we give a mild assumption, as in the general stochastic optimization \cite{ghadimi2016accelerated,ghadimi2016mini}
and the initial convex stochastic ADMM \cite{ouyang2013stochastic}.
\begin{assumption}
For smooth function $f(x)$, there exists a stochastic first-order oracle that returns a noisy estimation to the gradient of
$f(x)$, and the noisy estimation $G(x,\xi)$ satisfies
\begin{align}
 & \mathbb{E} [G(x,\xi)] = f(x), \\
 & \mathbb{E} \big[ \|G(x,\xi) - \nabla f(x)\|^2 \big] \leq \sigma^2,
\end{align}
where the expectation is taken with respect to the random variable $\xi$.
\end{assumption}

Let $M$ be the size of mini-batch $\mathcal{I}$,
and $\xi_{\mathcal{I}}=\{\xi_1,\xi_2,\cdots,\xi_M\}$ denotes a set of i.i.d.
random variables, and the stochastic gradient is given by
\begin{align}
G(x,\xi_{\mathcal{I}}) = \frac{1}{M} \sum_{i\in \mathcal{I}} G(x,\xi_i). \nonumber
\end{align}
Clearly, we have
\begin{align}
 & \mathbb{E}[G(x,\xi_{\mathcal{I}})] = \nabla f(x),  \label{eq:10}  \\
 & \mathbb{E}\big[ \|G(x,\xi_{\mathcal{I}})-\nabla f(x)\|^2\big] \leq \sigma^2/M. \label{eq:11}
\end{align}

\begin{algorithm}[htb]
   \caption{ Mini-batch Stochastic ADMM (STOC-ADMM) for Nonconvex Nonsmooth Optimization }
   \label{alg:1}
\begin{algorithmic}[1]
   \STATE {\bfseries Input:} Number of iteration $T$, Mini-batch size $0<M<n$ and $\rho>0$;
   \STATE {\bfseries Initialize:} $x_0$, $y_0$ and $\lambda_0$;
   \FOR {$t=0,1,\cdots,T-1$}
   \STATE{}  Uniformly randomly pick a mini-batch $\mathcal{I}_t$ from $\{1,2,\cdots,n\}$;
   \STATE{}  $y_{t+1}=\arg\min_y \mathcal {L}_{\rho}(x_t,y,\lambda_t)$;
   \STATE{}  $x_{t+1}=\arg\min_x \tilde{\mathcal {L}}_{\rho}\big(x;y_{t+1},\lambda_t,x_t,G(x_t,\xi_{\mathcal{I}_t})\big)$;
   \STATE{}  $\lambda_{t+1} = \lambda_t-\rho(Ax_{t+1}+By_{t+1}-c)$;
   \ENDFOR
   \STATE {\bfseries Output:} Iterate $x$ and $y$ chosen uniformly random from $\{x_{t},y_{t}\}_{t=1}^{T}$.
\end{algorithmic}
\end{algorithm}

In the stochastic ADMM algorithm, we can update $y$ and $\lambda$ by \eqref{eq:3} and \eqref{eq:5}, respectively,
as in the deterministic ADMM.
However, to update the variable $x$,
we will define an \emph{approximated} function of the form:
\begin{align} \label{eq:6}
\tilde{\mathcal {L}}_{\rho}\big(x;y_{t+1},\lambda_t,& x_t,G(x_t, \xi_{\mathcal{I}_t})\big) = f(x_t)+ G(x_t,\xi_{\mathcal{I}_t})^T(x-x_t)  \nonumber \\
  & + \frac{1}{2\eta}\|x-x_t\|_H^2 -\langle\lambda_t, Ax+By_{t+1}-c\rangle \nonumber \\
  & + \frac{\rho}{2} \|Ax+By_{t+1}-c\|^2,
\end{align}
where $\mathbb{E}[G(x_t,\xi_{\mathcal{I}_t})] = \nabla f(x_t)$, $\eta>0$ and $H \succ 0$.
By minimizing \eqref{eq:6} on the variable $x$, we have
\begin{align}
 x_{t+1} \!=\! (\frac{H}{\eta} \!+\! \rho A^TA)^{-1} \big[ \frac{H}{\eta}x_t \!-\! G(x_t,\xi_{\mathcal{I}_t}) \!-\! \rho A^T( By_{t+1} \!-\! c \!-\! \frac{\lambda_t}{\rho})\big]. \nonumber
\end{align}
When $A^TA$ is large, computing $(\frac{H}{\eta} + \rho A^TA)^{-1}$ is expensive, and
storage of this matrix may still be
problematic.
To avoid them, we can use the inexact Uzawa method \cite{zhang2011unified} to
linearize the last term in \eqref{eq:6}. In other words, we set $H=r I - \rho\eta A^TA$ with
\begin{align}
 r \geq r_{\min} \equiv \eta\rho \|A^TA\| +1    \nonumber
\end{align}
to ensure $H\succeq I$.
Then we have
\begin{align}
 x_{t+1} = x_t -\frac{\eta}{r}\big[ G(x_t,\xi_{\mathcal{I}_t}) + \rho A^T(A_t+By_{t+1}-c-\frac{\lambda_t}{\rho}) \big].
\end{align}
Finally, we give
the algorithmic framework of the mini-batch stochastic ADMM (STOC-ADMM) in Algorithm \ref{alg:1}.

\subsection{Convergence Analysis of Nonconvex Mini-batch STOC-ADMM }
In the subsection, we study the convergence and
iteration complexity of the nonconvex mini-batch STOC-ADMM.
First, we give some mild assumptions as follows:
\begin{assumption} \label{ass:2}
For smooth function $f(x)$, its gradient is
Lipschitz continuous with the constant $L>0$, such that
\begin{align}\label{eq:12}
\|\nabla f(x_1)-\nabla f(x_2)\|  \leq L \|x_1 - x_2\|, \ \forall x_1,x_2 \in R^d,
\end{align}
and this is equivalent to
\begin{align}\label{eq:13}
f(x_1) \leq f(x_2) + \nabla f(x_2)^T(x_1-x_2) + \frac{L}{2}\|x_1-x_2\|^2.
\end{align}
\end{assumption}
\begin{assumption} \label{ass:3}
Gradient of loss function $f(x)$ is bounded, i.e., there exists a constant $\delta >0$ such that for all $x$,
it follows $\|\nabla f(x)\|^2 \leq \delta^2$.
\end{assumption}
\begin{assumption} \label{ass:4}
$f(x) $ and $g(y)$ are all lower bounded, and denoting $f^*=\inf_x f(x)$ and $g^*=\inf_y g(y)$.
\end{assumption}
\begin{assumption} \label{ass:5}
$A$ is a full row or column rank.
\end{assumption}

Assumption \ref{ass:2} has been widely used in the convergence analysis
of nonconvex algorithms \cite{allen2016variance,reddi2016stochastic}.
Assumptions \ref{ass:3},\ref{ass:4} have been used in study of ADMM for nonconvex optimzations \cite{jiang2016structured}.
Assumption \ref{ass:5} has been used in the convergence analysis of ADMM \cite{jiang2016structured,deng2016global}. 
Assumption \ref{ass:5} guarantees the matrix $A^TA$ or $AA^T$ is non-singular.
Without loss of generality, we will use the full column rank matrix $A$ below.
Next, we define the $\epsilon$-stationary point of the nonconvex problem \eqref{eq:1} below:
\begin{definition} \label{def:1}
For $\epsilon>0$, the point $(x^*,y^*,\lambda^*)$ is said to be an $\epsilon$-stationary point of the nonconvex problem \eqref{eq:1} if it holds that
 \begin{equation} \label{eq:def}
  \left\{
  \begin{aligned}
  &\mathbb{E}\|Ax^*+By^*-c\|^2\leq \epsilon,   \\
  &\mathbb{E}\|\nabla f(x^*)-A^T\lambda^*\|^2\leq \epsilon,   \\
  &\mathbb{E}[\mbox{dist}\big(B^T\lambda^*,\partial g(y^*)\big)^2]\leq \epsilon,
  \end{aligned}
  \right.\end{equation}
where $dist\big(y_0,\partial g(y)\big):=\inf \{\|y_0-z\|: \ z\in \partial g(y) \}$, and
$\partial g(y)$ denotes the subgradient of $g(y)$.
If $\epsilon=0$, the point $(x^*,y^*,\lambda^*)$ is said to be a stationary point of \eqref{eq:1}.
\end{definition}
Note that the above inequalities \eqref{eq:def}
are equivalent to $\mathbb{E}\big[ \mbox{dist}\big(0,\partial L(x^*,y^*,\lambda^*)\big)^2 \big] \leq \epsilon$, where
\begin{align}
   \partial L(x,y,\lambda) = \left [ \begin{matrix}
     \partial L(x,y,\lambda)/\partial x \\
     \partial L(x,y,\lambda)/\partial y \\
     \partial L(x,y,\lambda)/\partial \lambda
 \end{matrix}
 \right ],\nonumber
\end{align}
where $L(x,y,\lambda)=f(x) + g(y) - \langle\lambda, Ax+By-c\rangle$ is the Lagrangian function of \eqref{eq:1}.
In the following, based the above assumptions and definition, we study the convergence and iteration complexity of the
mini-batch stochastic ADMM.

\begin{lemma} \label{lem:2}
 Suppose the sequence $\{x_t,y_t,\lambda_t\}_{t=1}^T$ is generated by Algorithm \ref{alg:1}. The following inequality holds
 \begin{align}
   \mathbb{E}\|\lambda_{t+1} \!-\! \lambda_{t}\|^2 \!\leq  \zeta \|x_{t}\!-\!x_{t-1}\|^2 \!+\! \zeta_1 \mathbb{E}\|x_{t+1}\!-\!x_t\|^2
   \!+\! \frac{10\sigma^2}{M\phi_{\min}^A}, \nonumber
 \end{align}
 where $\zeta = \frac{5(L^2\eta^2 + (\phi_{\max}^H)^2)}{\phi_{\min}^A\eta^2}$ and $\zeta_1 = \frac{ 5(\phi_{\max}^H)^2}{\phi_{\min}^A\eta^2} $.
\end{lemma}
A detailed proof of Lemma \ref{lem:2} is provided in \hyperref[app:lem-2]{Appendix A.1}.
Lemma \ref{lem:2} gives the upper bound of $\mathbb{E}\|\lambda_{t+1}-\lambda_{t}\|^2$.
Given a sequence $\{x_t,y_t,\lambda_t\}_{t=1}^T$ generated from Algorithm \ref{alg:1},
then we define a useful sequence $\big\{\Psi_t\big\}_{t=1}^T$ as follows:
\begin{align} \label{eq:16}
\Psi_t =  \mathbb{E}\big[ \mathcal{L}_{\rho}(x_t,y_t,\lambda_t) + \frac{5(L^2\eta^2 + (\phi_{\max}^H)^2)}{\rho \phi_{\min}^A\eta^2}\|x_t-x_{t-1}\|^2 \big].
\end{align}
For notational simplicity, let $\tilde{L}=L+1$, $\phi^H = (\phi^H_{\min})^2+20(\phi^H_{\max})^2$ and
$\varphi=(\tilde{L}+10L^2/(\rho \phi^A_{\min}))-\phi^A_{\min}\rho$.
\begin{lemma} \label{lem:3}
 Suppose that the sequence $\{x_t,y_t,\lambda_t\}_{t=1}^T$ is generated by Algorithm \ref{alg:1}.
 Let $\rho_*=\frac{\tilde{L}+\sqrt{40L^2+\tilde{L}^2}}{2\phi_{\min}^A}$, $\triangle = (\phi^H_{\min})^2 + \frac{20(\phi^H_{\max})^2}{\rho \phi^A_{\min}}\big(\phi^A_{\min}\rho -(\tilde{L}+\frac{10L^2}{\rho\phi^A_{\min}})\big)$, and
 $$\rho_0=\frac{10\phi^H_{\max}\big(\tilde{L}\phi^H_{\max}+\sqrt{\tilde{L}^2(\phi^H_{\max})^2+2L^2 \phi^H}\big)}
 {\phi^A_{\min}\phi^H}$$
 and suppose the parameters $\rho$ and $\eta$, respectively, satisfy
 \begin{equation} \label{eq:le3}
  \left\{
  \begin{aligned}
  & \eta \in \big( \frac{\phi^H_{\min}-\sqrt{\triangle}}{\varphi}, \frac{\phi^H_{\min}+\sqrt{\triangle}}{\varphi} \big), \quad \rho \in \big(\rho_0, \rho_* \big); \\
  & \eta\in \big(\frac{10(\phi_{\max}^H)^2}{\rho \phi_{\min}^A\phi_{\min}^H}, \frac{r-1}{\rho\|A^TA\|}\big], \quad \rho = \rho_*;\\
  & \eta \in \big(\frac{\phi^H_{\min}-\sqrt{\triangle}}{\varphi},\frac{r-1}{\rho\|A^TA\|}\big] , \quad \rho\in (\rho_*,+\infty).
  \end{aligned}
  \right.\end{equation}
 Then we have $\gamma=  \frac{\phi^H_{\min}}{\eta} + \frac{\phi^A_{\min}\rho}{2} - \frac{\tilde{L}}{2}- \frac{5(L^2\eta^2 + 2(\phi^H_{\max})^2)}{\rho\phi^A_{\min}\eta^2}>0$, and it holds that
 \begin{align}\label{eq:17}
 \frac{1}{T} \sum_{t=0}^{T-1} \mathbb{E}\|x_t-x_{t+1}\|^2 \leq \frac{\Psi_0-\Psi^*}{\gamma T} + \frac{(\phi_{\min}^A\rho+20)\sigma^2}{2\gamma\phi_{\min}^A\rho M}.
 \end{align}
 where $\Psi^*$ is a lower bound of sequence $\big\{\Psi_t\big\}_{t=1}^T$. 
\end{lemma}
A detailed proof of Lemma \ref{lem:3} is provided in \hyperref[app:lem-3]{Appendix A.2}.
Lemma \ref{lem:3} gives a property of the sequence $\big\{\Psi_t\big\}_{t=1}^T$.
Moreover, \eqref{eq:le3} provides a specific parameter selection on the step size $\eta$ and
the penalty parameter $\rho$, in which selection of the step size $\eta$ depends on
the parameter $\rho$. Next, we define a useful variable $\theta_{t}$ defined by:
\begin{align}
\theta_{t} = \big[\|x_{t+1}-x_t\|^2+\|x_t-x_{t-1}\|^2 \big]. 
\end{align}

\begin{theorem} \label{th:1}
 Suppose the sequence $\{x_t,y_t,\lambda_t\}_{t=1}^T$ is generated by Algorithm \ref{alg:1}. Define
 $ \kappa_1=3(L^2+\frac{(\phi^H_{\max})^2}{\eta^2})$, $\kappa_2 = \frac{\zeta}{\rho^2}$,
 $\kappa_3=\rho^2\|B\|^2\|A\|^2$, and $\kappa_4=\frac{\phi_{\min}^A\rho+20}{2\phi_{\min}^A\rho }$.
 Let
 \begin{align} \label{eq:th5}
  & M \geq \frac{2\sigma^2}{\epsilon} \max \{\kappa_1\kappa_4+3,\kappa_2\kappa_4+\frac{10}{\phi_{\min}^A\rho^2},\kappa_3\kappa_4\}, \nonumber \\
  & T =\frac{\max\{\kappa_1,\kappa_2,\kappa_3\}}{\epsilon\gamma}(\Psi_1 - \Psi^*), \nonumber
 \end{align}
 where $\Psi^*$ is a lower bound of the sequence $\{\Psi_t\}_{t=1}^T$.
 Let $t^* = \mathop{\arg\min}_{2\leq t \leq T+1 }\theta_{t}$,
 then $(x_{t^*},y_{t^*})$ is an $\epsilon$-stationary point of the problem \eqref{eq:1}.
\end{theorem}

A detailed proof of Theorem \ref{th:1} is provided in \hyperref[app:th-1]{Appendix A.4}.
Theorem \ref{th:1} shows that, given an mini-batch size $M=O(1/\epsilon)$,
the mini-batch stochastic ADMM has the convergence rate of $O(\frac{1}{T})$ to obtain  an $\epsilon$-stationary point of the
nonconvex problem \eqref{eq:1}.
Moreover, the \emph{IFO}(Incremental First-order Oracle \cite{reddi2016stochastic}) complexity of
the mini-batch stochastic ADMM is $O(M/\epsilon)=O(1/\epsilon^2)$ for obtaining an $\epsilon$-stationary point.
While, the \emph{IFO} complexity of the
deterministic proximal ADMM \cite{jiang2016structured} is $O(n/\epsilon)$ for obtaining an $\epsilon$-stationary point.
When $n>\frac{1}{\epsilon}$, the mini-batch stochastic ADMM needs less \emph{IFO} complexity than
the deterministic ADMM.

In the convergence analysis,
given an appropriate mini-batch size $M$ satisfies the condition \eqref{eq:th5},
the step size $\eta$ only need satisfies the condition \eqref{eq:le3} instead of
$\eta=O(\frac{1}{\sqrt{t}})$ used in the convex stochastic ADMM \cite{ouyang2013stochastic}.
\section{ Nonconvex Minin-batch SVRG-ADMM }
In the subsection, we propose a mini-batch nonconvex stochastic variance reduced gradient ADMM (SVRG-ADMM)
to solve the problem \eqref{eq:1}, which uses a multi-stage strategy
to progressively reduce the variance of stochastic gradients.

Algorithm \ref{alg:2} gives an algorithmic framework of mini-batch SVRG-ADMM for nonconvex optimizations.
In Algorithm \ref{alg:2}, the stochastic gradient $\hat{\nabla} f(x_{t}^{s+1}) =
\frac{1}{M}\sum_{i_t\in \mathcal{I}_t} \big(\nabla f_{i_t}(x_{t}^{s+1})-\nabla f_{i_t}(\tilde{x}^s)\big)+\nabla f(\tilde{x}^s)$ is unbiased, i.e.,
$\mathbb{E}[\hat{\nabla} f(x_{t}^{s+1})]=\nabla f(x_{t}^{s+1})$.
In the following, we give an upper bound of variance of
the stochastic gradient $\hat{\nabla} f(x_{t}^{s+1})$.

\begin{lemma} \label{lem:6}
 In Algorithm \ref{alg:2}, set $\Delta^{s+1}_t=\hat{\nabla}f(x^{s+1}_t)-\nabla f(x^{s+1}_t)$,
 then it holds
 \begin{align}
  \mathbb{E}\|\Delta^{s+1}_t\|^2 \leq \frac{L^2}{M}\|x^{s+1}_t-\tilde{x}^{s}\|^2,
 \end{align}
 where $\mathbb{E}\|\Delta^{s+1}_t\|^2$ denotes variance of the stochastic gradient $\hat{\nabla}f(x^{s+1}_t)$.
\end{lemma}
A detailed proof of Lemma \ref{lem:6} is provided in \hyperref[app:lem-6]{Appendix B.1}.
Lemma \ref{lem:6} shows that the variance of the stochastic gradient $\hat{\nabla}f(x^{s+1}_t)$
has an upper bound $O(\|x_t^{s+1}-\tilde{x}^s\|^2)$.
Due to $\tilde{x}^s = x^s_m $, as number of iterations increases,
both $x_t^{s+1}$ and $\tilde{x}^s$ approach the same stationary point,
thus the variance of stochastic gradient vanishes.
In fact, the variance of stochastic gradient $\hat{\nabla}f(x^{s+1}_t)$ is progressively reduced.

\begin{algorithm}[htb]
   \caption{ Mini-batch SVRG-ADMM for Nonconvex Nonsmooth Optimization }
   \label{alg:2}
\begin{algorithmic}[1]
   \STATE {\bfseries Input:} Mini-batch size $M$, epoch length $m$, $T$, $S=[T/m]$, $\rho>0$;
   \STATE {\bfseries Initialize:} $\tilde{x}^0=x_m^0$, $y_m^0$ and $\lambda_m^0$;
   \FOR {$s=0,1,\cdots,S-1$}
   \STATE{} $x_0^{s+1}=x_{m}^s$, $y_0^{s+1}=y_{m}^s$ and $\lambda_0^{s+1}=\lambda_{m}^s$;
   \STATE{} $\nabla f(\tilde{x}^s)=\frac{1}{n}\sum_{i=1}^n\nabla f_i(\tilde{x}^s)$;
   \FOR {$t=0,1,\cdots,m-1$}
   \STATE{} Uniformly randomly pick a mini-batch $\mathcal{I}_t$ from $\{1,2,\cdots,n\}$;
   \STATE{} $y^{s+1}_{t+1}=\arg\min_y \mathcal {L}_{\rho}(x^{s+1}_t,y,\lambda_t^{s+1})$;
   \STATE{} $\hat{\nabla} f(x_{t}^{s+1}) = \frac{1}{M}\sum_{i_t\in \mathcal{I}_t} \big(\nabla f_{i_t}(x_{t}^{s+1})-\nabla f_{i_t}(\tilde{x}^s)\big)+\nabla f(\tilde{x}^s)$;
   \STATE{} $x^{s+1}_{t+1}=\arg\min_x \tilde{\mathcal {L}}_{\rho}\big(x;y_{t+1}^{s+1},\lambda_{t}^{s+1},x_{t}^{s+1},\hat{\nabla} f(x_{t}^{s+1})\big)$;
   \STATE{} $\lambda_{t+1}^{s+1} = \lambda_{t}^{s+1}-\rho(Ax_{t+1}^{s+1}+By_{t+1}^{s+1}-c)$;
   \ENDFOR
   \STATE{} $\tilde{x}^{s+1}= x_m^{s+1}$;
   \ENDFOR
   \STATE {\bfseries Output:} Iterate $x$ and $y$ chosen uniformly random from $\{(x_{t}^s,y_{t}^s)_{t=1}^{m}\}_{s=1}^S$.
\end{algorithmic}
\end{algorithm}

\subsection{Convergence Analysis of Nonconvex Mini-batch SVRG-ADMM }
In the subsection, we study the convergence and iteration complexity of the mini-batch nonconvex SVRG-ADMM.
First, we give an upper bound of $\mathbb{E} \|\lambda^{s+1}_{t+1}-\lambda^{s+1}_{t}\|^2$.
\begin{lemma} \label{lem:7}
 Suppose the sequence $\{(x^{s}_t,y^{s}_t,\lambda^{s}_t)_{t=1}^m\}_{s=1}^S$ is generated by Algorithm \ref{alg:2}. The following inequality holds
 \begin{align}
 \mathbb{E}\|\lambda^{s+1}_{t+1} \! -\! \lambda^{s+1}_{t}\|^2 & \!\leq \frac{5L^2}{\phi^A_{\min}M} \mathbb{E} \|x^{s+1}_{t}\!-\!\tilde{x}^{s}\|^2
   \!+\! \frac{5L^2}{\phi^A_{\min}M}\|x^{s+1}_{t-1}\!-\!\tilde{x}^{s}\|^2 \nonumber \\
  & + \zeta \|x^{s+1}_{t}-x^{s+1}_{t-1}\|^2 + \zeta_1 \mathbb{E}\|x^{s+1}_{t+1}-x^{s+1}_t\|^2. \nonumber
 \end{align}
\end{lemma}
A detailed proof of Lemma \ref{lem:7} is provided in \hyperref[app:lem-7]{Appendix B.2}.
Given the sequence $\{(x^s_t,y^s_t,\lambda^s_t)_{t=1}^m\}_{s=1}^S$ generated from Algorithm \ref{alg:2}, then we define
a useful sequence $\big\{(\Phi^{s}_{t})_{t=1}^m\big\}_{s=1}^S$ as follows:
\begin{align} \label{eq:88}
 \Phi^s_t
  = & \mathbb{E}\big[ \mathcal {L}_{\rho}(x^{s}_{t},y^{s}_{t},\lambda^{s}_{t})+ h^{s}_{t}(\|x^{s}_{t}-\tilde{x}^{s-1}\|^2+ \|x^{s}_{t-1}-\tilde{x}^{s-1}\|^2) \nonumber \\
   &+ \frac{\zeta}{ \rho}\|x^{s}_{t}-x^{s}_{t-1}\|^2\big],
\end{align}
where $\{(h^s_t)_{t=1}^m \}_{s=1}^S$ is a positive sequence.

\begin{lemma} \label{lem:8}
 Suppose the sequence $\{(x^{s}_t,y^{s}_t,\lambda^{s}_t)_{t=1}^m\}_{s=1}^S$ is generated from Algorithm \ref{alg:2},
 and suppose the positive sequence $\{(h_t^s)_{t=1}^m\}_{s=1}^S$ satisfies, for $s =1,2,\cdots,S$
 \begin{equation}
  h^{s}_t= \left\{
  \begin{aligned}
  & (2+\beta)h^{s}_{t+1}+ \frac{(10+\phi^A_{\min}\rho)L^2}{2\rho\phi^A_{\min}M}, \ 1 \leq t \leq m-1, \\
  & \frac{10L^2}{\phi^A_{\min}\rho M}, \quad t=m,
  \end{aligned}
  \right.\end{equation}
 where $\beta>0$. Let $\hat{h}=\min_t\{(1+\frac{1}{\beta})h_{t+1}^s,h_1^{s+1}\}$, $\triangle_1 = (\phi^H_{\min})^2 + \frac{20(\phi^H_{\max})^2}{\rho \phi^A_{\min}}\big(\phi^A_{\min}\rho -(\tilde{L}+2\hat{h}+\frac{10L^2}{\rho\phi^A_{\min}})\big)$,
 $\rho_*=\frac{\tilde{L}+2\hat{h}+\sqrt{40L^2+(\tilde{L}+2\hat{h})^2}}{2\phi_{\min}^A}$,
 and
 $$\rho_0=\frac{10\phi^H_{\max}\bigg((\tilde{L}+2\hat{h})\phi^H_{\max}+\sqrt{(\tilde{L}+2\hat{h})^2(\phi^H_{\max})^2+2L^2\phi^H }\bigg)}
 {\phi^A_{\min} \phi^H }$$
 and suppose the parameters $\rho$ and $\eta$, respectively, satisfy
 \begin{equation} \label{eq:le8}
  \left\{
  \begin{aligned}
  & \eta \in \big( \frac{\phi^H_{\min}-\sqrt{\triangle_1}}{\varphi_1}, \frac{\phi^H_{\min}+\sqrt{\triangle_1}}{\varphi_1} \big), \quad \rho \in \big(\rho_0, \rho_* \big); \\
  & \eta\in \big(\frac{10(\phi_{\max}^H)^2}{\rho \phi_{\min}^A\phi_{\min}^H}, \frac{r-1}{\rho\|A^TA\|}\big], \quad \rho = \rho_*;\\
  & \eta \in \big(\frac{\phi^H_{\min}-\sqrt{\triangle_1}}{\varphi_1},\frac{r-1}{\rho\|A^TA\|}\big] , \quad \rho\in (\rho_*,+\infty).
  \end{aligned}
  \right.\end{equation}
 where $\varphi_1 = (\tilde{L}+2\hat{h}+10L^2/(\rho \phi^A_{\min}))-\phi^A_{\min}\rho$.
 Then it holds that the sequence $\{(\Gamma^{s}_t)_{t=1}^m\}_{s=1}^S$ is positive, defined by
 \begin{equation}
 \Gamma^{s}_t \!=\! \left\{
 \begin{aligned}
   & \frac{\phi^H_{\min}}{\eta} \!+\! \frac{\phi^A_{\min}\rho}{2}\!-\! \frac{\tilde{L}}{2} \!-\!\frac{\zeta\!+\!\zeta_1}{\rho}\!-\!(1\!+\!\frac{1}{\beta})h_{t+1},
   \ 1 \!\leq\! t \!\leq\! m-1  \\
   & \frac{\phi^H_{\min}}{\eta} \!+\! \frac{\phi^A_{\min}\rho}{2}\!-\!\frac{\tilde{L}}{2}-\frac{\zeta+\zeta_1}{\rho}\!-\! h_1^{s+1}, \quad t=m
 \end{aligned}
 \right.\end{equation}
  and the sequence $\big\{(\Phi^{s}_{t})_{t=1}^m\big\}_{s=1}^S$  monotonically decreases.
\end{lemma}

A detailed proof of Lemma \ref{lem:8} is provided in \hyperref[app:lem-8]{Appendix B.3}.
Lemma \ref{lem:8} shows that the sequence $\big\{(\Phi^{s}_{t})_{t=1}^m\big\}_{s=1}^S$
monotonically decreases.
Moreover, \eqref{eq:le8} provides a specific parameter selection on the step size $\eta$ and
the penalty parameter $\rho$ in Algorithm \ref{alg:2}.

\begin{lemma} \label{lem:9}
Suppose the sequence $\{(x^{s}_t,y^{s}_t,\lambda^{s}_t)_{t=1}^m\}_{s=1}^S$ is generated by Algorithm \ref{alg:2}.
Under the same conditions as in Lemma \ref{lem:8}, the sequence
$\big\{(\Phi^{s}_{t})_{t=1}^m\big\}_{s=1}^S$ has a lower bound.
\end{lemma}

Lemma \ref{lem:9} shows that the sequence
$\big\{(\Phi^{s}_{t})_{t=1}^m\big\}_{s=1}^S$ has a lower bound.
The proof of Lemma \ref{lem:9} is the same as the proof of Lemma 7 in \cite{huang2016stochastic}.
Next, we define a useful variable $\hat{\theta}^s_t$ as follows:
 \begin{align} \label{eq:48}
  \hat{\theta}^{s}_t\!= & \mathbb{E} \big[ \|x^{s}_{t} \!-\! \tilde{x}^{s-1}\|^2 \!+ \|x^{s}_{t-1}\!-\!\tilde{x}^{s-1}\|^2 \!+ \|x^{s}_{t+1}\!-\!x^{s}_t\|^2 \nonumber \\
  & \!+ \|x^{s}_{t}-x^{s}_{t-1}\|^2 \big].
 \end{align}
In the following, we will analyze the convergence properties
of the nonconvex SVRG-ADMM based on the above lemmas.

\begin{theorem} \label{th:2}
 Suppose the sequence $\{(x^{s}_t,y^{s}_t,\lambda^{s}_t)_{t=1}^m\}_{s=1}^S$ is generated by Algorithm \ref{alg:2}.
 Denote
 $\kappa_1=3\big(L^2+\frac{(\phi^H_{\max})^2}{\eta^2}\big)$, $\kappa_2=\frac{\zeta}{\rho^2}$, $\kappa_3=\rho^2\|B\|^2\|A\|^2$,
 and $\gamma=\min_{(t,s)} \Gamma^s_t$ and $\omega = \min_{(s,t)}\{(2+\beta)h^{s}_{t+1}+\frac{L^2}{2M},\frac{5L^2}{\phi^A_{\min}\rho M} \}$.
 Let
 \begin{align}
  mS = T = \frac{\max\{\kappa_1,\kappa_2,\kappa_3\}}{\tau \epsilon}(\Phi^{1}_{1}- \Phi^*),
 \end{align}
 where $\tau=\min(\gamma,\omega)$, and $\Phi^*$ is a lower bound of the sequence $\big\{(\Phi^{s}_{t})_{t=1}^m\big\}_{s=1}^S$ .
 Let
 \begin{align}
  (t^*,s^*) = \mathop{\arg\min}_{1 \leq t\leq m,\ 1 \leq s\leq S}\hat{\theta}^{s}_t, \nonumber
 \end{align}
 then $(x_{t^*}^{s^*}, y_{t^*}^{s^*})$ is an $\epsilon$-stationary point of the problem \eqref{eq:1}.
\end{theorem}

A detailed proof of Theorem \ref{th:2} is provided in \hyperref[app:th-2]{Appendix B.4}.
Theorem \ref{th:2} shows that the mini-batch SVRG-ADMM for nonconvex optimizations has
a convergence rate of $O(\frac{1}{T})$. Moreover, the \emph{IFO} complexity of the mini-batch SVRG is
$O\big( (\frac{n}{m}+M)/\epsilon \big)$. 
When $\frac{n}{m}+ M <n$, the mini-batch SVRG-ADMM needs less \emph{IFO} complexity than
the deterministic ADMM.

Since the mini-batch SVRG-ADMM uses VR technique, its convergence
does not depend on the mini-batch size $M$. In other words, when $M=1$, the mini-batch nonconvex SVRG-ADMM reduces to
the initial nonconvex SVRG-ADMM in \cite{huang2016stochastic}, which also has a convergence rate of $O(\frac{1}{T})$.
However, by Lemma \ref{lem:6},
the variance of stochastic gradient in the mini-batch SVRG-ADMM decreases faster than that
in the initial nonconvex SVRG-ADMM.

\section{ Nonconvex Mini-batch SAGA-ADMM }
In the subsection, we propose a mini-batch nonconvex stochastic average gradient ADMM (SAGA-ADMM)
by additionally using the old gradients estimated in the previous iteration,
which is inspired by the SAGA method \cite{defazio2014saga}.

The algorithmic framework of the SAGA-ADMM is given in Algorithm \ref{alg:3}.
In Algorithm \ref{alg:3}, the stochastic gradient $\hat{\nabla} f(x_{t})=\frac{1}{M}\sum_{i_t\in \mathcal{I}_t}
\big(\nabla f_{i_t}(x_{t})-\nabla f_{i_t}(z^t_{i_t}) \big)+\psi_t$ is unbiased
(i.e., $\mathbb{E}[\hat{\nabla} f(x_{t})]=\nabla f(x_{t})$), where $\psi_t=\frac{1}{n}\sum_{i=1}^n\nabla f_i(z^t_i)$.
In the following, we give an upper bound of the variance of the stochastic gradient $\hat{\nabla} f(x_{t})$.
\begin{lemma} \label{lem:11}
 For Algorithm \ref{alg:3}, Let $\Delta_t=\hat{\nabla}f(x_t)-\nabla f(x_t)$, then it holds
 \begin{align}
  \mathbb{E}\|\Delta_t\|^2 \leq \frac{L^2}{Mn}\sum_{i=1}^n \|x_t-z_i^t\|^2,
 \end{align}
 where $\mathbb{E}\|\Delta_t\|^2$
 denotes variance of the stochastic gradient $\hat{\nabla}f(x_t)$.
\end{lemma}
A detailed proof of Theorem \ref{lem:11} is provided in \hyperref[app:lem-11]{Appendix C.1}.
Lemma \ref{lem:11} shows that the variance of the stochastic gradient $\hat{\nabla}f(x_t)$
has an upper bound $O(\frac{1}{n}\sum_{i=1}^n \|x_t-z_i^t\|^2)$.
As the number of iteration increases, both $x_t$ and the stored points $\{z^t\}_{i=1}^n$ approach the same stationary point,
so the variance of stochastic gradient progressively reduces.
In fact, the variance of stochastic gradient $\hat{\nabla}f(x_t)$
is progressively reduced via additionally using the old gradients in the previous iterations.
\begin{algorithm}[htb]
   \caption{ Mini-batch SAGA-ADMM for Nonconvex Nonsmooth Optimization }
   \label{alg:3}
\begin{algorithmic}[1]
   \STATE {\bfseries Input:} $x_0 \in R^d$, $y_0\in R^q$, $z_i^0=x_0$ for $i\in \{1,2,\cdots,n\}$, number of iterations $T$;
   \STATE {\bfseries Initialize:} $\psi_0=\frac{1}{n}\sum_{i=1}^n\nabla f_i(z^0_i)$;
   \FOR {$t=0,1,\cdots,T-1$}
   \STATE{} Uniformly randomly pick a mini-batch $\mathcal{I}_t$ from $\{1,2,\cdots,n\}$;
   \STATE{} $y_{t+1}=\arg\min_y \mathcal {L}_{\rho}(x_t,y,\lambda_t)$;
   \STATE{} $\hat{\nabla} f(x_{t}) = \frac{1}{M}\sum_{i_t\in \mathcal{I}_t} \big(\nabla f_{i_t}(x_{t})-\nabla f_{i_t}(z^t_{i_t}) \big)+\psi_t$
            with $\psi_t=\frac{1}{n}\sum_{i=1}^n\nabla f_i(z^t_i)$;
   \STATE{} $x_{t+1}=\arg\min_x \tilde{\mathcal {L}}_{\rho}\big(x;y_{t+1},\lambda_{t},x_{t},\hat{\nabla} f(x_{t})\big)$;
   \STATE{} $\lambda_{t+1} = \lambda_{t}-\rho(Ax_{t+1}+By_{t+1}-c)$;
   \STATE{} $z^{t+1}_{i_t}= x_{t+1}$ and $z_i^{t+1}=z^t_i$ for $i\neq i_t$, for all $i_t\in \mathcal{I}_t$;
   \STATE{} $\psi_{t+1}=\psi_t-\frac{1}{n}\sum_{i_t\in \mathcal{I}_t} \big(\nabla f_{i_t}(z^t_{i_t})-\nabla f_{i_t}(z^{t+1}_{i_t})\big)$;
   \ENDFOR
   \STATE {\bfseries Output:} Iterate $x$ and $y$ chosen uniformly random from $\{x_{t},y_{t}\}_{t=1}^{T}$.
\end{algorithmic}
\end{algorithm}

\subsection{Convergence Analysis of Nonconvex Mini-batch SAGA-ADMM}

In the subsection, we study the convergence and iteration complexity of the nonconvex mini-batch SAGA-ADMM.
First, we give same useful lemmas as follows:
\begin{lemma} \label{lem:12}
 Suppose the sequence $\{x_t,y_t,\lambda_t\}_{t=1}^T$ is generated by Algorithm \ref{alg:3}. The following inequality holds
 \begin{align}
  \mathbb{E}\|\lambda_{t+1}\!-\lambda_{t}\|^2 \leq & \frac{5L^2}{\phi^A_{\min}M n} \sum_{i=1}^n \big( \mathbb{E} \|x_{t}\!-z^{t}_i\|^2
   + \|x_{t-1}\!-z^{t-1}_i\|^2 \big) \nonumber \\
  & + \zeta\|x_{t}-x_{t-1}\|^2 + \zeta_1\mathbb{E}\|x_{t+1}-x_t\|^2. \nonumber
 \end{align}
\end{lemma}
Lemma \ref{lem:12} gives an upper bound of $\mathbb{E}\|\lambda_{t+1}-\lambda_{t}\|^2$.
Its proof is the same as that of Lemma \ref{lem:7}.
Given the sequence $\{x_t,y_t,\lambda_t\}_{t=1}^T$ generated by Algorithm \ref{alg:3},
then we define a useful sequence $\{\Theta_t\}_{t=1}^T$ below:
\begin{align} \label{eq:89}
\Theta_{t} =  & \mathbb{E} \big[ \mathcal {L}_{\rho}(x_{t},y_{t},\lambda_{t})
  + \frac{\alpha_{t}}{n}\sum_{i=1}^n(\|x_{t}-z^{t}_i\|^2 +\|x_{t-1}-z^{t-1}_i\|^2) \nonumber \\
 & + \frac{\zeta}{\rho}\|x_{t}-x_{t-1}\|^2 \big],
\end{align}
where $\{\alpha_t\}_{t=1}^T$ is a decreasing positive sequence.

\begin{lemma} \label{lem:13}
 Suppose that the sequence $\{x_t,y_t,\lambda_t\}_{t=1}^T$ is generated by Algorithm \ref{alg:3},
 and the positive sequence $\{\alpha_{t}\}_{t=1}^T$ satisfy
 \begin{align}
  \alpha_t=\frac{10L^2+\phi^A_{\min}\rho L^2}{2\rho\phi^A_{\min}M} + (\frac{2n-M}{n}+\frac{n-M}{n}\beta)\alpha_{t+1},
 \end{align}
 where $\beta>0$. Let $\hat{\alpha}=\min_t\{\frac{n-M}{n}(1+\frac{1}{\beta})\alpha_{t+1}\}$,
 $\triangle_2 = (\phi^H_{\min})^2 + \frac{20(\phi^H_{\max})^2}{\rho \phi^A_{\min}}\big(\phi^A_{\min}\rho -(\tilde{L}+2\hat{\alpha}+\frac{10L^2}{\rho\phi^A_{\min}})\big)$,
 $\rho_*=\frac{\tilde{L}+2\hat{\alpha}+\sqrt{40L^2+(\tilde{L}+2\hat{\alpha})^2}}{2\phi_{\min}^A}$, and
 $$\rho_0\!=\!\frac{10\phi^H_{\max}\bigg((\tilde{L}+2\hat{\alpha})\phi^H_{\max}\!+\!\sqrt{(\tilde{L}+2\hat{\alpha})^2(\phi^H_{\max})^2\!+\!2L^2\phi^H}\bigg)}
 {\phi^A_{\min}\phi^H}$$
 and suppose the parameters $\rho$ and $\eta$, respectively, satisfy
 \begin{equation} \label{eq:le13}
  \left\{
  \begin{aligned}
  & \eta \in \big( \frac{\phi^H_{\min}-\sqrt{\triangle_2}}{\varphi_2}, \frac{\phi^H_{\min}+\sqrt{\triangle_2}}{\varphi_2} \big), \quad \rho \in \big(\rho_0, \rho_* \big); \\
  & \eta\in \big(\frac{10(\phi_{\max}^H)^2}{\rho \phi_{\min}^A\phi_{\min}^H}, \frac{r-1}{\rho\|A^TA\|}\big], \quad \rho = \rho_*;\\
  & \eta \in \big(\frac{\phi^H_{\min}-\sqrt{\triangle_2}}{\varphi_2},\frac{r-1}{\rho\|A^TA\|}\big] , \quad \rho\in (\rho_*,+\infty),
  \end{aligned}
  \right.\end{equation}
 where $\varphi_2 = (\tilde{L}+2\hat{\alpha}+10L^2/(\rho \phi^A_{\min}))-\phi^A_{\min}\rho$.
 Then it holds the sequence $\{\Gamma_t\}_{t=1}^T$ is positive, defined by
 \begin{align}
\Gamma_{t} = \frac{\phi^H_{\min}}{\eta} + \frac{\phi^A_{\min}\rho}{2}-\frac{\tilde{L}}{2} - \frac{\zeta+\zeta_1}{\rho}
   -\frac{n-M}{n}(1+\frac{1}{\beta})\alpha_{t+1},
 \end{align}
 and the sequence $\{\Theta_{t}\}_{t=1}^T $ monotonically decreases.
\end{lemma}
A detailed proof of Lemma \ref{lem:13} is provided in \hyperref[app:lem-13]{Appendix C.2}.
Lemma \ref{lem:13} shows that the sequence $\{\Theta_{t}\}_{t=1}^T $ monotonically decreases.
Moreover, \eqref{eq:le13} provides a specific parameter selection on the step size $\eta$ and
the penalty parameter $\rho$ in Algorithm \ref{alg:3}.

\begin{lemma} \label{lem:14}
Suppose the sequence $\{x_t,y_t,\lambda_t\}_{t=1}^T$ is generated by Algorithm \ref{alg:3}. Under the same conditions as in Lemma \ref{lem:13}, the sequence
$\{\Theta_{t}\}_{t=1}^T$ has a lower bound.
\end{lemma}
Lemma \ref{lem:14} shows that the sequence
$\{\Theta_{t}\}_{t=1}^T$ has a lower bound.
Its proof is the same as the proof of Lemma 7 in \cite{huang2016stochastic}.
In the following, we will study the convergence and iteration complexity of the SAGA-ADMM based on the above lemmas.
We define a useful variable $\tilde{\theta}_{t}$ defined by:
 \begin{align}
  \tilde{\theta}_{t} = & \big[\|x_{t+1}-x_t\|^2+\|x_t-x_{t-1}\|^2 + \frac{1}{n}\sum_{i=1}^n ( \|x_t-z^t_i\|^2  \nonumber \\
  & + \|x_{t-1}-z^{t-1}_i\|^2 )\big]. \label{eq:72}
 \end{align}

\begin{theorem} \label{th:3}
 Suppose the sequence $\{x_t,y_t,\lambda_t\}_{t=1}^T$ is generated by Algorithm \ref{alg:3}. Denote
 $ \kappa_1=3\big(L^2+\frac{(\phi^H_{\max})^2}{\eta^2}\big)$, $\kappa_2 = \frac{\zeta}{\rho^2}$,
 $\kappa_3=\rho^2\|B\|^2\|A\|^2$, and $\gamma=\min_t\Gamma_t$ and $\omega=\min_t \big\{ \frac{L^2}{2M} + (\frac{2n-M}{n}+\frac{n-M}{n}\beta)\alpha_{t+1}\big\}$.
 Let
 \begin{align}
 T =\frac{\max\{\kappa_1,\kappa_2,\kappa_3\}}{\tau\epsilon}(\Theta_{1}- \Theta^*),
 \end{align}
 where $\tau = \min\big\{ \gamma, \omega \big\}>0$, and $\Theta^*$ is a lower bound of
 the sequence $\big\{\Theta_{t}\big\}_{t=1}^T$.
 Let $t^* = \mathop{\arg\min}_{2\leq t \leq T+1 }\tilde{\theta}_{t}$,
 then $(x_{t^*},y_{t^*})$ is an $\epsilon$-stationary point of the problem \eqref{eq:1}.
\end{theorem}

A detailed proof of Theorem \ref{th:3} is provided in \hyperref[app:th-3]{Appendix C.3}.
Theorem \ref{th:3} shows that the mini-batch SAGA-ADMM for nonconvex optimizations has a convergence rate of
$O(\frac{1}{T})$. Moreover, the \emph{IFO} complexity of the mini-batch SAGA-ADMM is $O(M/\epsilon)$ for obtaining an $\epsilon$-stationary point.
Clearly, due to $1\leq M <n$, the mini-batch SAGA-ADMM needs less \emph{IFO} complexity than
the deterministic ADMM.

Since the mini-batch SAGA-ADMM also uses VR technique, its convergence
does not depend on the mini-batch size $M$.
In other words, when $M=1$, the mini-batch nonconvex SAGA-ADMM reduces to
the initial nonconvex SAGA-ADMM in \cite{huang2016stochastic}, which also has the convergence rate of $O(\frac{1}{T})$.
However, by Lemma \ref{lem:11},
the variance of stochastic gradient in the mini-batch nonconvex SAGA-ADMM
decreases faster than that in the initial nonconvex SAGA-ADMM.

\begin{table}
  \centering
  \caption{Comparing the best \emph{IFO} and \emph{EI} complexity of different algorithms. The complexity is measured in terms
  of the number of \emph{oracle} calls required to achieve an $\epsilon$-stationary point (see Definition \ref{def:2}).}\label{tab:2}
  \resizebox{1.0\columnwidth}{!}{
\begin{tabular}{|c|c|c|c|}
  \hline
  % after \\: \hline or \cline{col1-col2} \cline{col3-col4} ...
  Algorithms & \emph{IFO} & \emph{EI} & Fixed Step Size ?\\ \hline
  Deterministic ADMM  & $O(n/\epsilon)$ & $O(1/\epsilon)$ & \checkmark \\ \hline
  Mini-batch STOC-ADMM  & $O(1/\epsilon^2)$  & $O(1/\epsilon)$ & \checkmark\\ \hline
  Mini-batch SVRG-ADMM  &$O\big(n +M/\epsilon\big)$& $O(1/\epsilon)$ & \checkmark \\  \hline
  Mini-batch SAGA-ADMM  & $O(n+M/\epsilon)$  & $O(1/\epsilon)$ &  \checkmark \\
  \hline
\end{tabular}
}
\end{table}

Finally, in Table \ref{tab:2},
we give the \emph{IFO}(Incremental First-order Oracle \cite{reddi2016stochastic}) and \emph{EI} (Effective Iteration)
of both the mini-batch stochastic ADMMs and the deterministic (or batch) ADMM.
Specifically, the definition of \emph{EI} is given in Definition \ref{def:2}.
From Table \ref{tab:2}, we can find that though both the mini-batch stochastic and deterministic ADMMs have the same \emph{EI} complexity,
the mini-batch stochastic ADMMs has lower \emph{IFO} complexity than the deterministic ADMM when $M<n$.
In the above theoretical analysis, the mini-batch size $M$ of the nonconvex STOC-ADMM may be very large when $\epsilon$ is small.
However, the following extensive experimental results show that
STOC-ADMM still has good performances given a moderate $M$, and is comparable with
both SVRG-ADMM and SAGA-ADMM.

\begin{definition} \label{def:2}
 For ADMM and its variants, an \emph{EI} describes the fact that all the primal and dual variables in the algorithm
 are updated once.
\end{definition}

\section{ Experiments }
In this section, we perform some numerical experiments on both simulated and real-world data to examine
performances of the proposed algorithms for the nonconvex nonsmooth optimization\footnote{We will put our code online once this paper is accepted.}.
In the experiments, we compare nonconvex mini-batch stochastic ADMM (\emph{STOC-ADMM}) with
nonconvex mini-batch \emph{SVRG-ADMM}, nonconvex mini-batch \emph{SAGA-ADMM}
and deterministic ADMM (\emph{DETE-ADMM}).
In the experiments, we use the inexact Uzawa method to both mini-batch stochastic ADMMs
and deterministic (or batch) ADMM.
In the following, all algorithms are implemented in MATLAB,
and all experiments are performed on a PC with an Intel E5-2630 CPU and 32GB memory.

\subsection{ Simulated Data }
In the subsection, we compare the performances in some synthetic data.
Here we focus on the binary classification task problem with the graph-guided fused lasso and
the overlapping group lasso regularization functions, respectively.
Given a set of training samples $(a_i,b_i)_{i=1}^n$,
where $a_i\in R^d$, $b_i \in \{-1,+1\}$,
then we solve the following nonconvex nonsmooth optimization problem:
\begin{align} \label{eq:80}
 \min_{x\in R^d} \frac{1}{n}\sum_{i=1}^n f_i(x) + \nu\|Ax\|_1,
\end{align}
where $f_i(x)=\frac{1}{1+\exp(b_i a_i^Tx)}$ is the \emph{sigmoid loss} function \cite{allen2016variance},
which is \emph{nonconvex} and smooth,
and $\nu$ denotes a nonegative regularization parameter.
When using graph-guided fused lasso \cite{kim2009multivariate} in \eqref{eq:80},
we let $A$ to decode the sparsity pattern of graph,
which is obtained by sparse precision matrix estimation \cite{friedman2008sparse,hsieh2014quic}.
When using overlapping group lasso \cite{kim2009multivariate} in \eqref{eq:80},
we let $Ax$ be concatenation of $k-$repetitions of $x$ (i.e., $Ax=[x;\ldots;x]$ ) as in \cite{suzuki2013dual},
where $k$ denotes the number of overlapping group of unknown parameter.

\subsubsection{ Graph-guided Fused Lasso }
Here we compare the performances in some simulated data, where a graph-guided fused lasso
regularization is imposed.
First, we generate a sparse precision matrix $\Lambda\in R^{d\times d}$ with elements
\begin{equation}
\Lambda_{ij} \mathop{\sim}^{i.i.d.} \begin{cases}  0,  & \mbox{prob.}\ 0.95 \\
 Unif([-0.75,-0.25]\bigcup [0.25,0.75]), & \mbox{otherwise}. \end{cases}  \nonumber
\end{equation}
Then the input feature vectors $\{a_i\}_{i=1}^n$ are i.i.d. generated from
multivariate normal distribution $N(0,\Lambda^{-1})$.
The true vector parameter $x^*\in R^d$ is generated from the standard normal distribution.
The output label is generated as $b_i = sign(a_i^Tx^*+\epsilon_i)$, where $\epsilon_i$
is chosen uniformly at random from $[0,1]$.

\begin{figure*}[htbp]
\centering
\subfigure[$n=20,000$]{\includegraphics[width=0.32\textwidth]{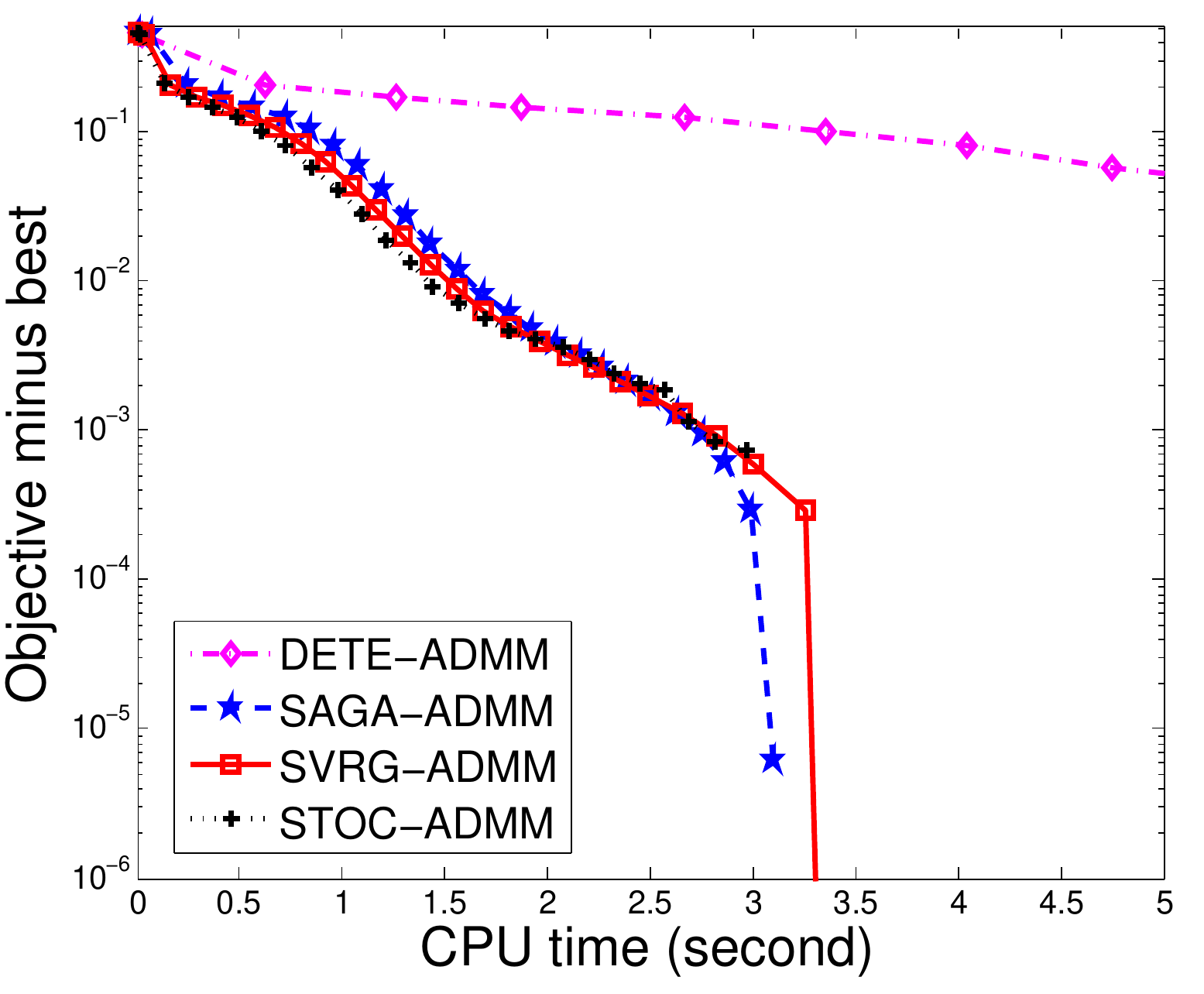}}
\subfigure[$n=40,000$]{\includegraphics[width=0.32\textwidth]{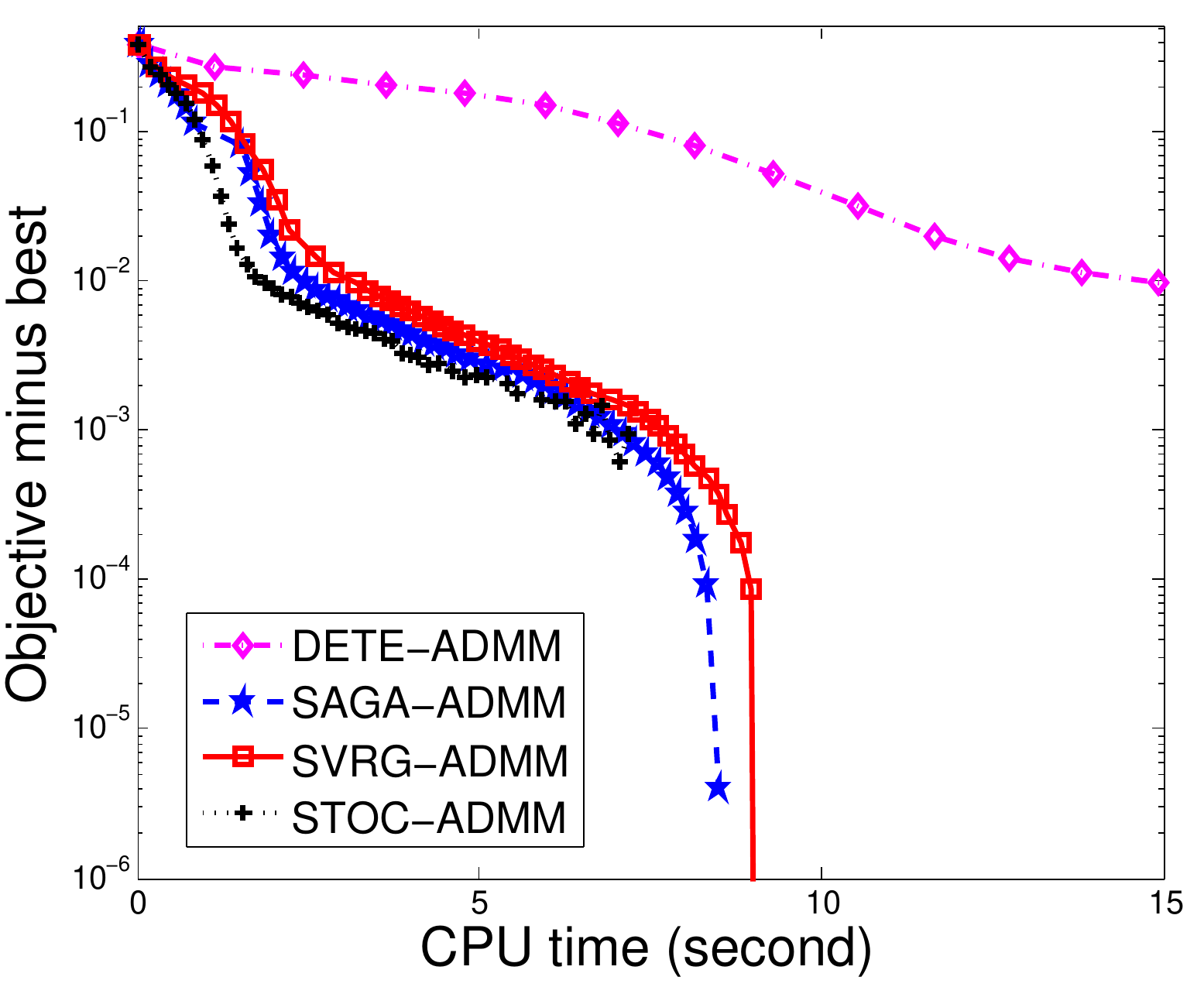}}
\subfigure[$n=60,000$]{\includegraphics[width=0.32\textwidth]{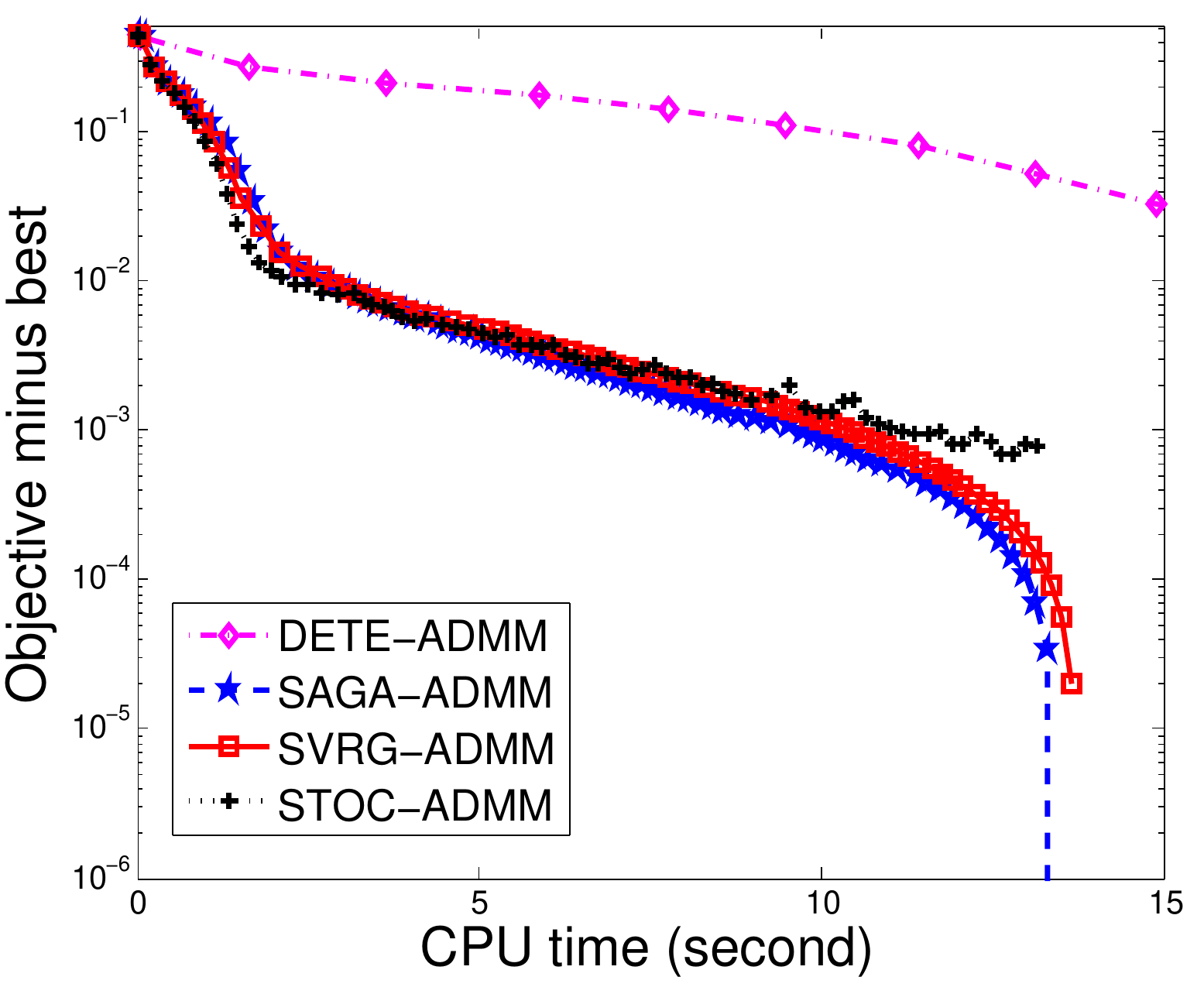}}
\caption{Objective value \emph{versus} CPU time on the simulated \emph{nonconvex} model with \textbf{graph-guided fused Lasso}.}
\label{fig:1}
\vspace{-1em}
\end{figure*}

\begin{figure*}[htbp]
\centering
\subfigure[$n=20,000$]{\includegraphics[width=0.32\textwidth]{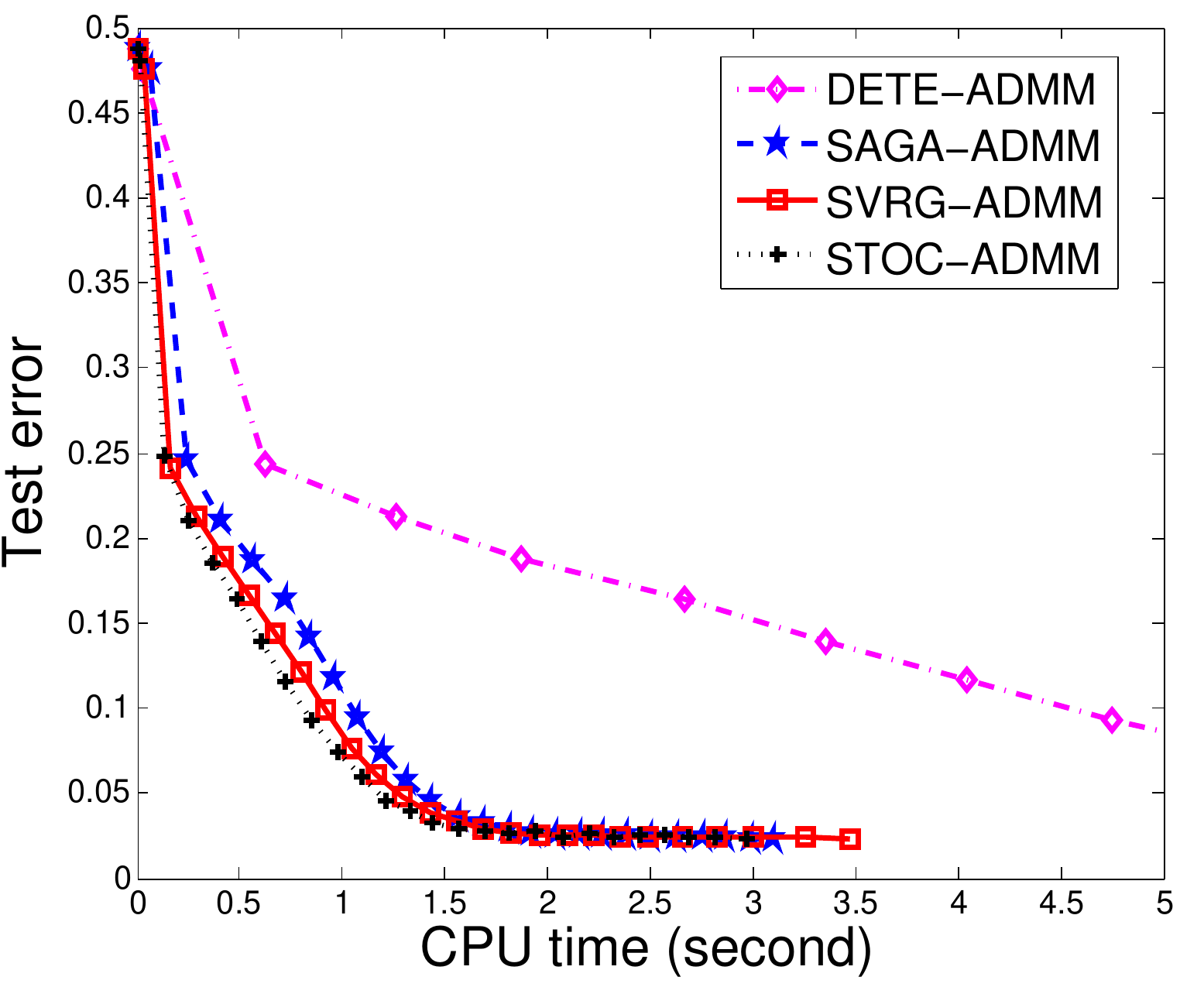}}
\subfigure[$n=40,000$]{\includegraphics[width=0.32\textwidth]{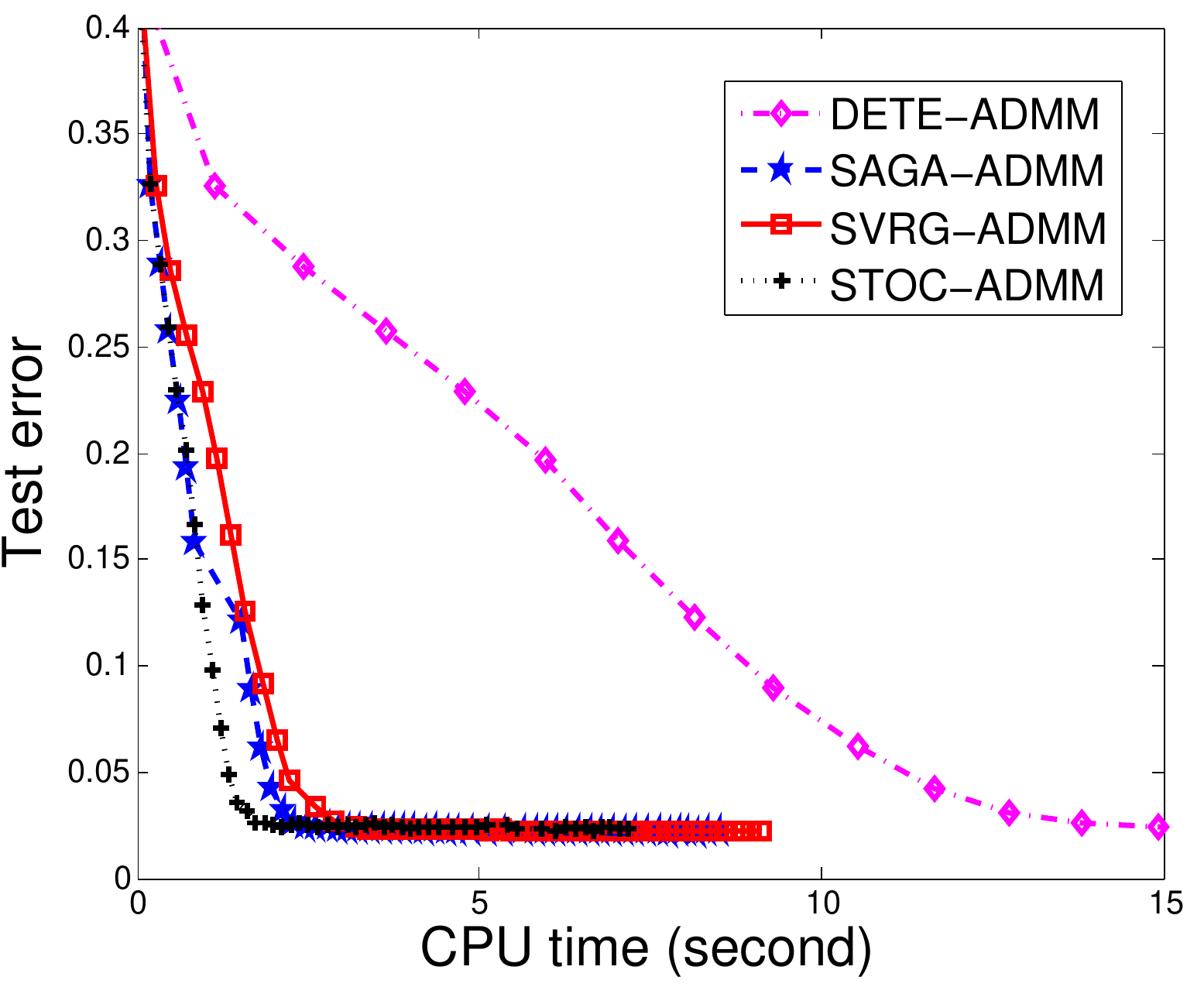}}
\subfigure[$n=60,000$]{\includegraphics[width=0.32\textwidth]{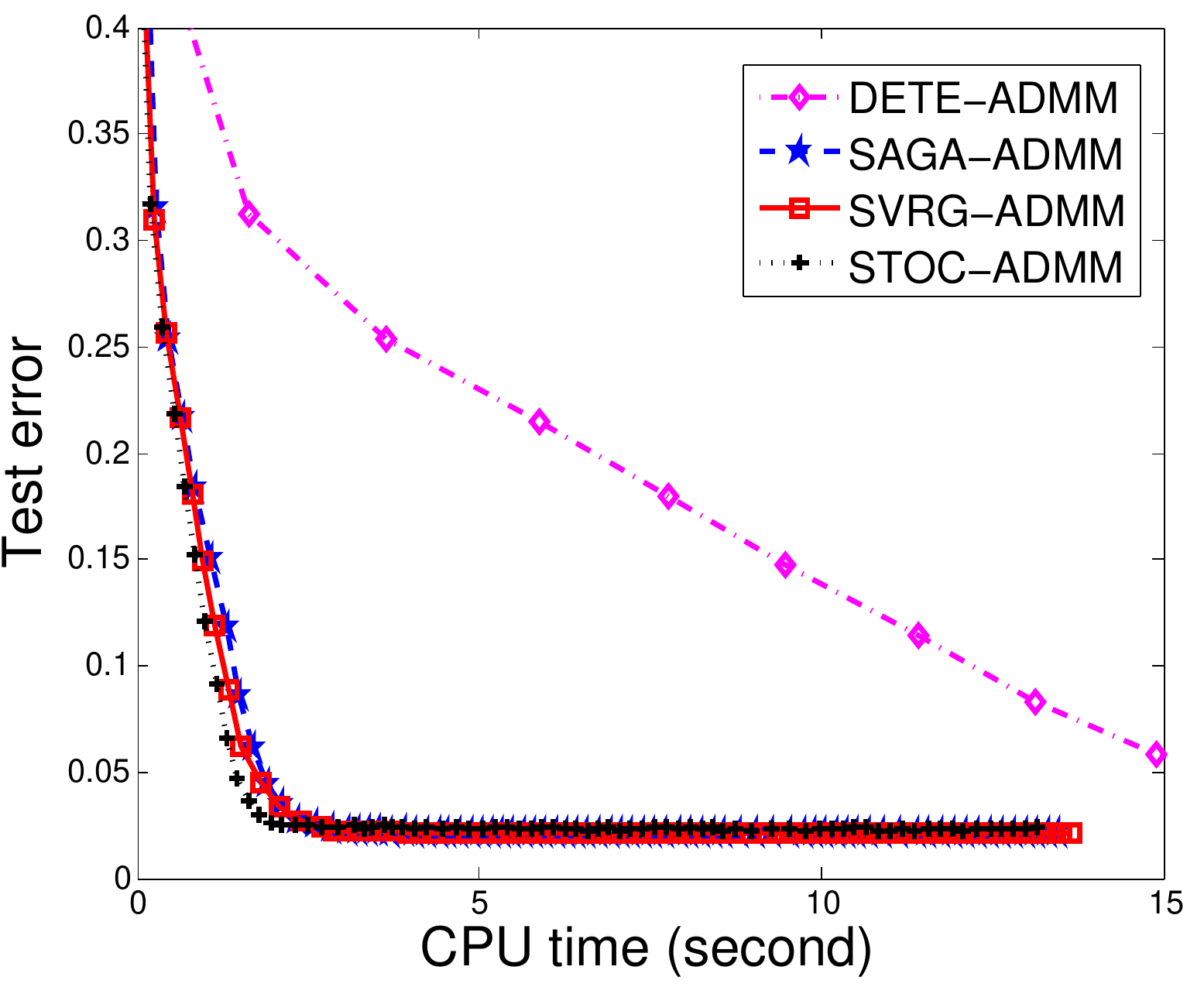}}
\caption{Test error \emph{versus} CPU time on the simulated \emph{nonconvex} model with \textbf{graph-guided fused Lasso}.}
\label{fig:2}
\vspace{-1em}
\end{figure*}

In the experiment, we set $d=200$, and then generate $n=\{20000, 40000, 60000\}$ samples $(a_i,b_i)_{i=1}^n$, respectively.
For each dataset, we choose half of the samples as training data, while use the rest as testing data.
In the problem \eqref{eq:80}, we use a graph-guided fused lasso, and fix the regularization parameter $\nu=10^{-5}$.
In the algorithms, we use the same initial solution $x_0$ from the standard normal distribution
and choose the step size $\eta = 1$.
In addition, we choose the mini-batch size $M=100$ in the stochastic algorithms, and
$m=[n/M]$ in the mini-batch SVRG-ADMM.
Finally, all experimental results are averaged over 10 repetitions.

Figs.\ref{fig:1} and \ref{fig:2} show that both the objective values and test loss of these
stochastic ADMMs faster decrease than those of the deterministic ADMM,
as CPU time consumed increases.
In particular, though the nonconvex \emph{STOC-ADMM} uses a fixed step size $\eta$,
it shows good performance in the nonconvex optimization with graph-guided fused lasso regularization,
and is comparable with both the nonconvex \emph{SVRG-ADMM} and \emph{SAGA-ADMM}.

\subsubsection{ Overlapping Group Lasso }
Here we compare the performances in some simulated data,
where an overlapping group lasso regularization is imposed.
First, we generate $n$ input feature vector $\{a_i\}_{i=1}^n$ with the dimension $d=400$,
where each feature is i.i.d. generated from the standard normal distribution.
Next, we generate a sparse matrix $X\in R^{20\times 20}$ , where only
the first column is non-zero (generated i.i.d. from standard normal distribution)
and other columns are zero, and the true parameter vector $x^*$ is vectorization of the matrix $X$.
The output label $b_i$ is generated as $b_i = sign(a_i^Tx^*+\epsilon_i)$, where $\epsilon_i$
is the standard normal distribution.
Then, we generate $n=\{20000, 40000, 60000\}$ samples $(a_i,b_i)_{i=1}^n$, respectively.
In the experiment, we choose the mini-batch size $M=200$ in these stochastic algorithms,
and the other settings are the similar as the above graph-guided fused lasso task.
In the problem \eqref{eq:80}, we use a overlapping group lasso
penalty function $g(x) = \nu(\sum_{i=1}^{20}\|X_{i,.}\| + \sum_{j=1}^{20}\|X_{.,j}\|)$.
Then we let $A = [I;I]$ as in \cite{suzuki2013dual}, and fix the parameter $\nu=10^{-5}$.

\begin{figure*}[htbp]
\centering
\subfigure[$n=20,000$]{\includegraphics[width=0.32\textwidth]{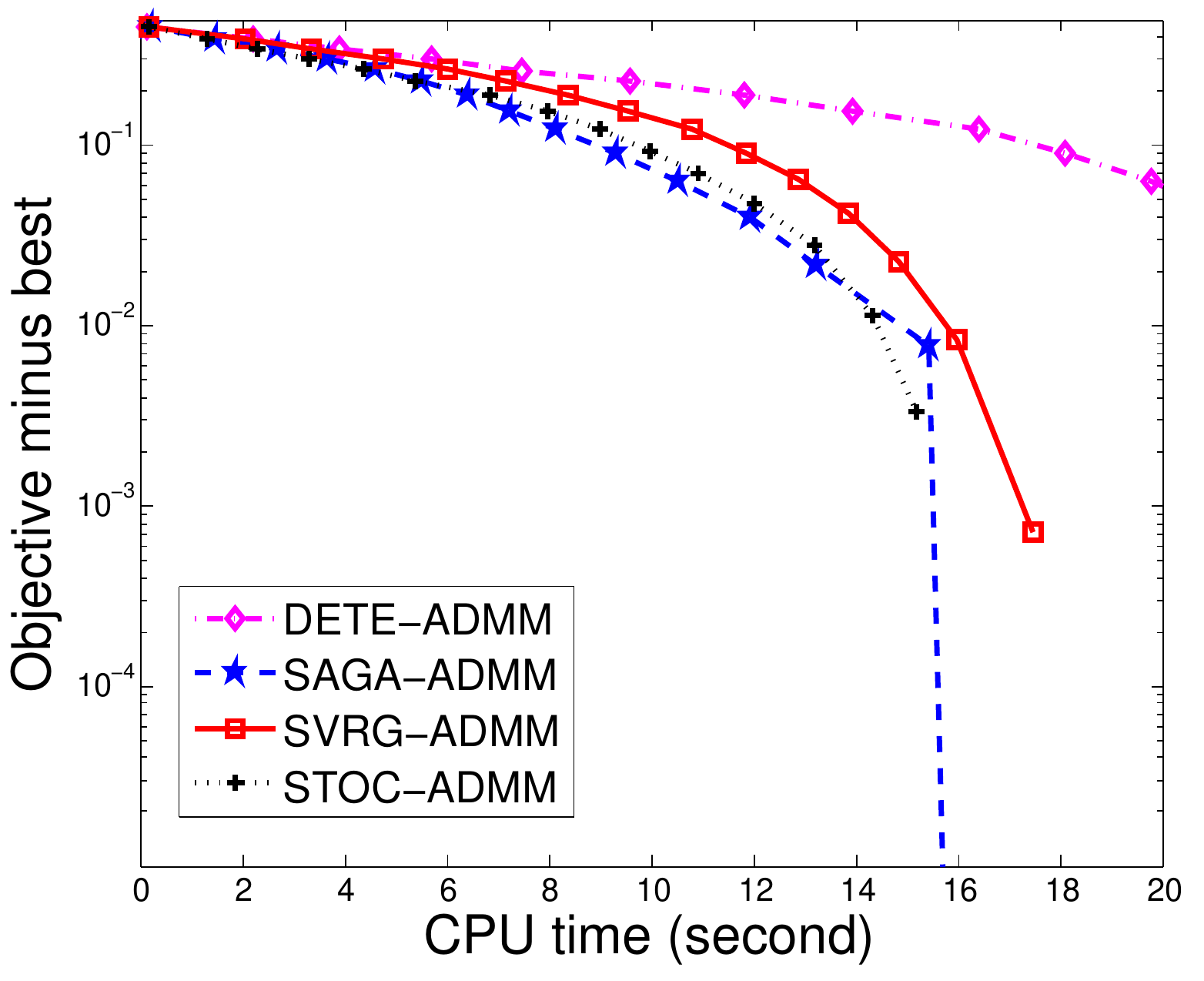}}
\subfigure[$n=40,000$]{\includegraphics[width=0.32\textwidth]{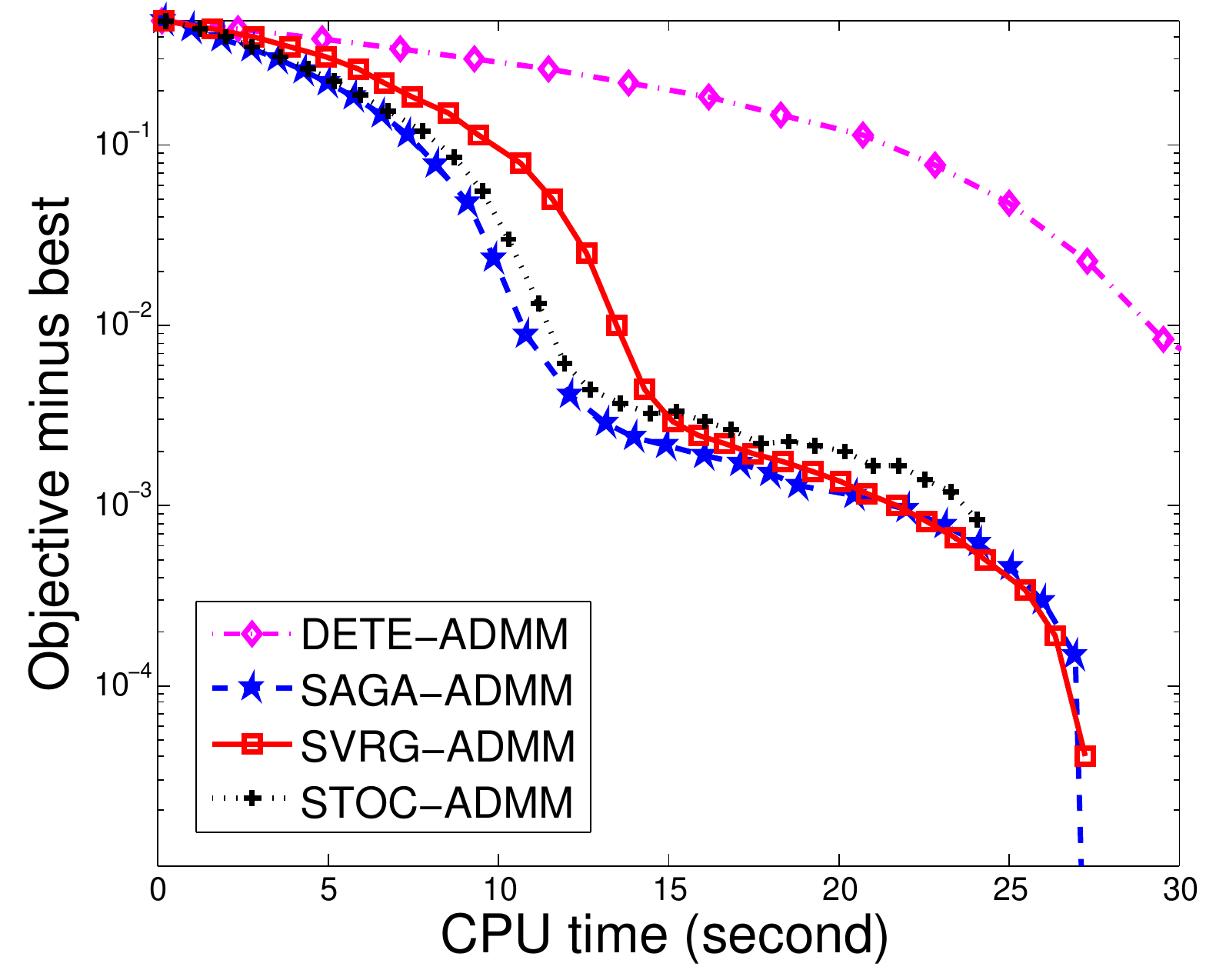}}
\subfigure[$n=60,000$]{\includegraphics[width=0.32\textwidth]{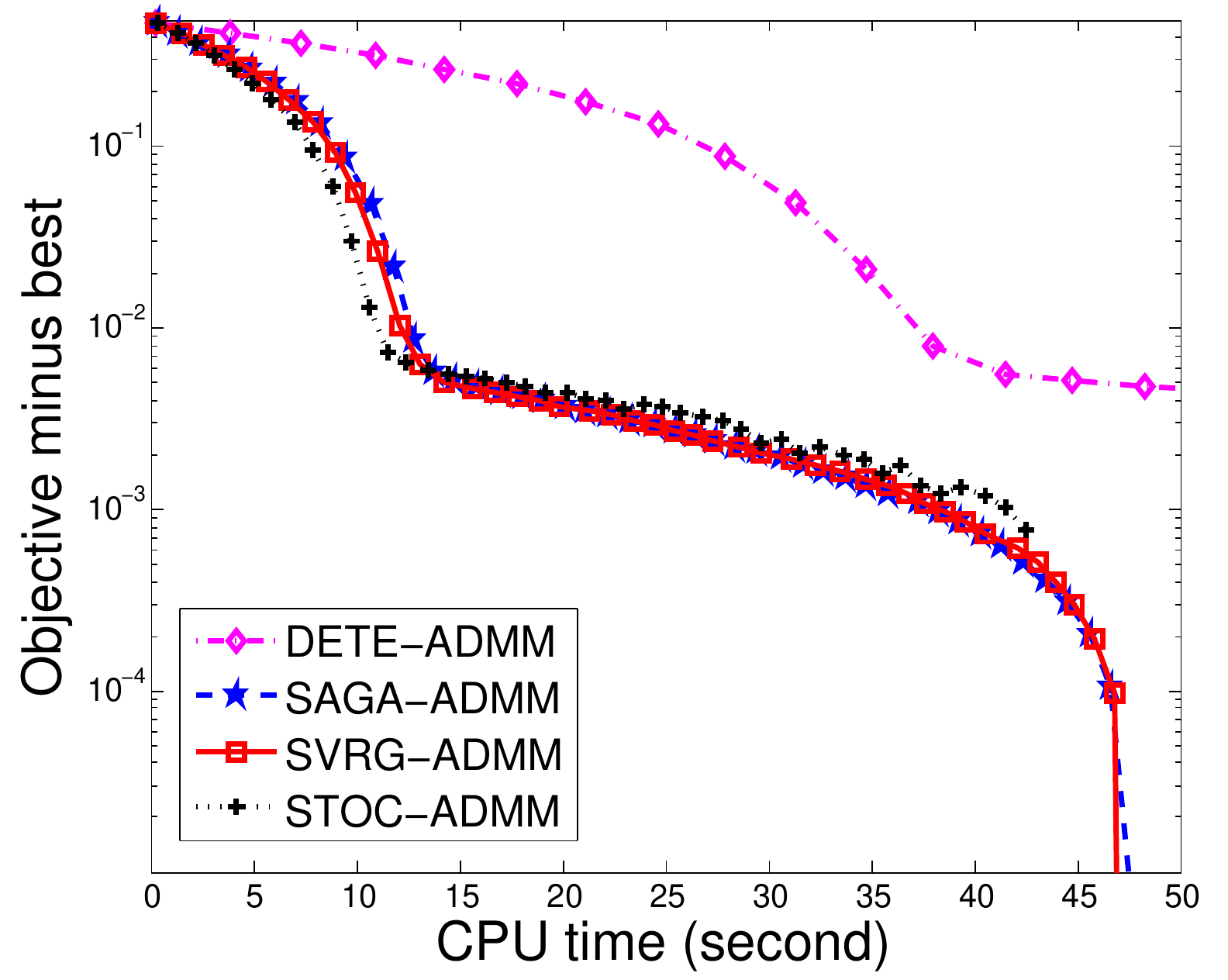}}
\caption{Objective value \emph{versus} CPU time on the simulated \emph{nonconvex} model with \textbf{overlapping group Lasso}.}
\label{fig:3}
\vspace{-1em}
\end{figure*}

\begin{figure*}[htbp]
\centering
\subfigure[$n=20,000$]{\includegraphics[width=0.32\textwidth]{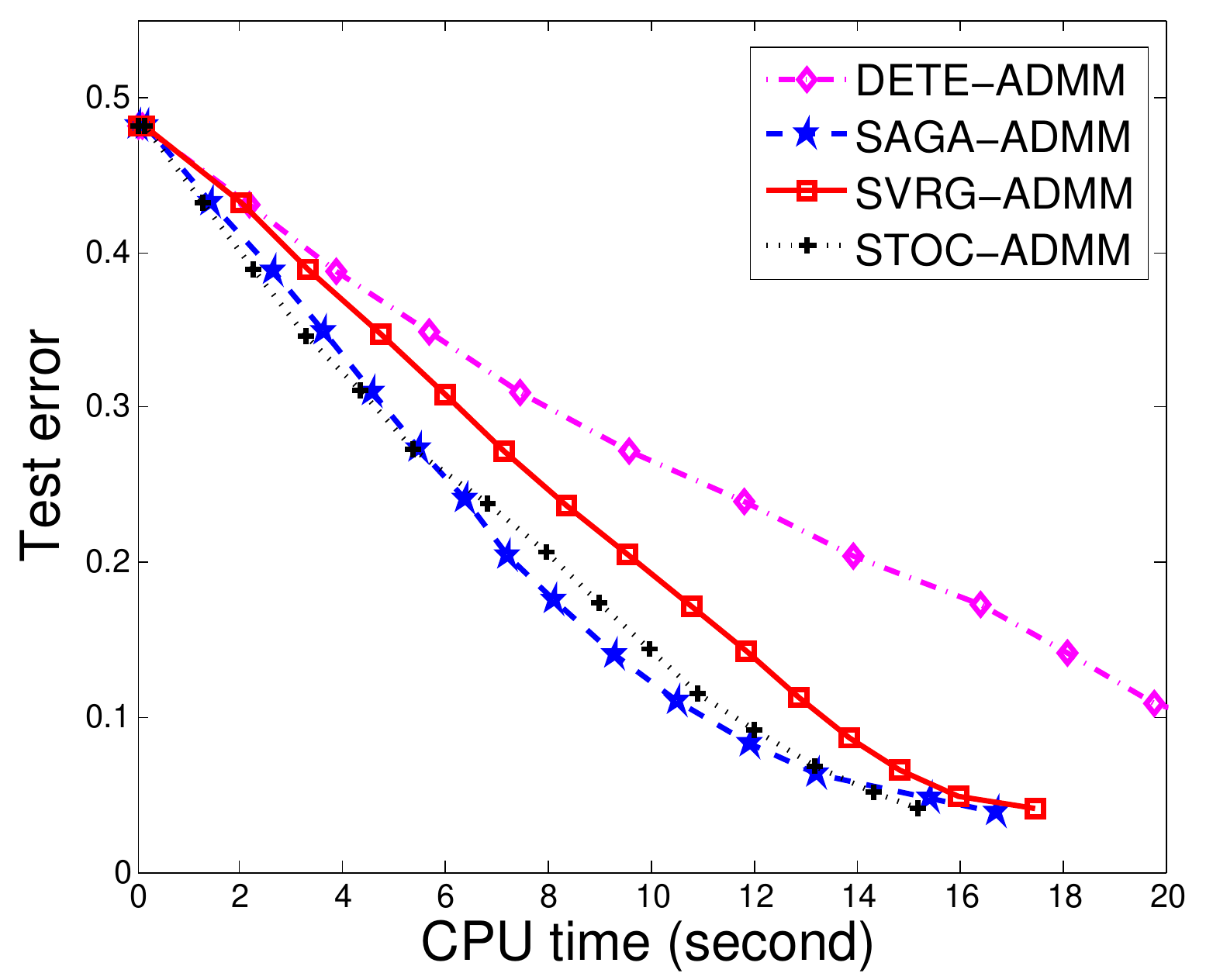}}
\subfigure[$n=40,000$]{\includegraphics[width=0.32\textwidth]{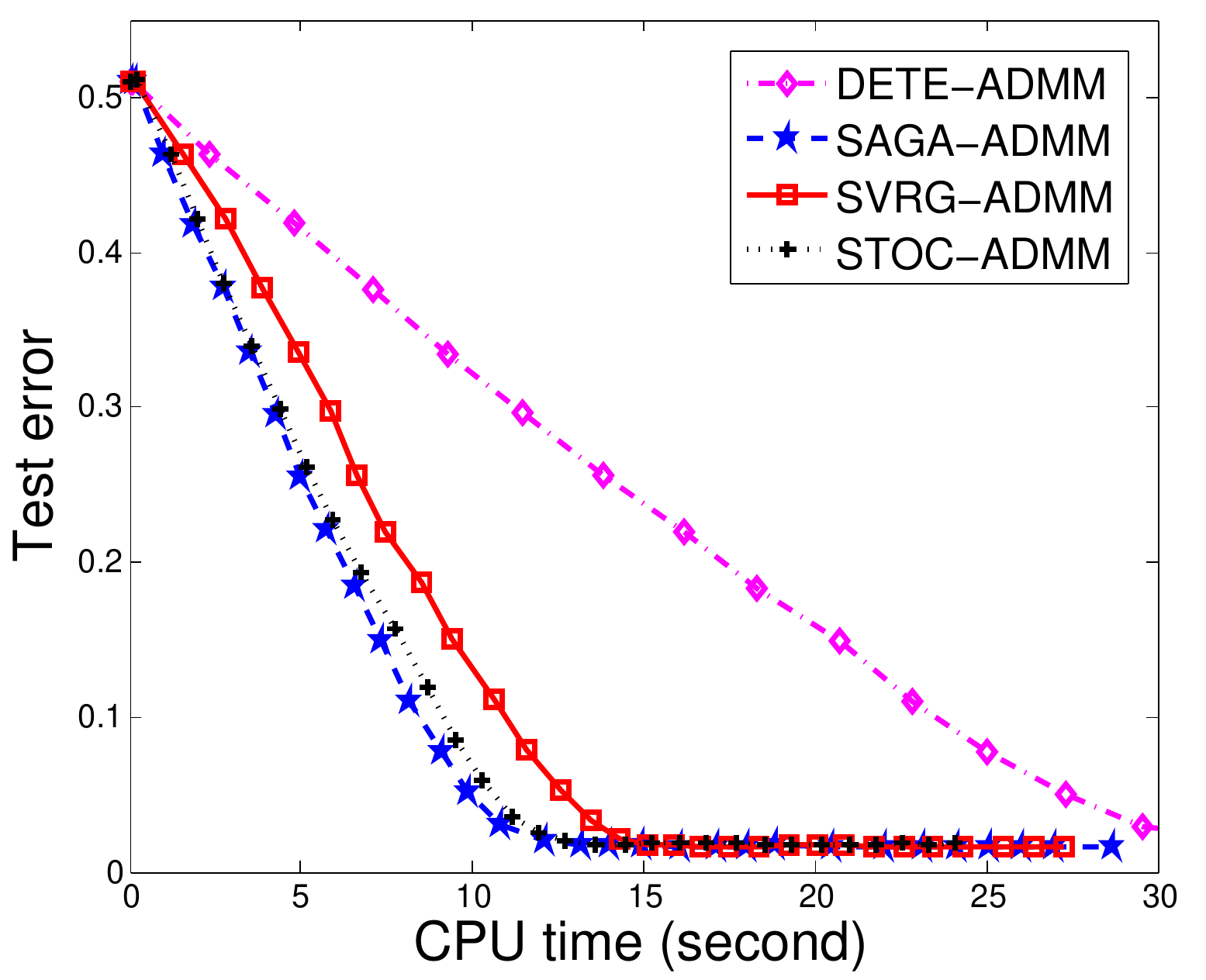}}
\subfigure[$n=60,000$]{\includegraphics[width=0.32\textwidth]{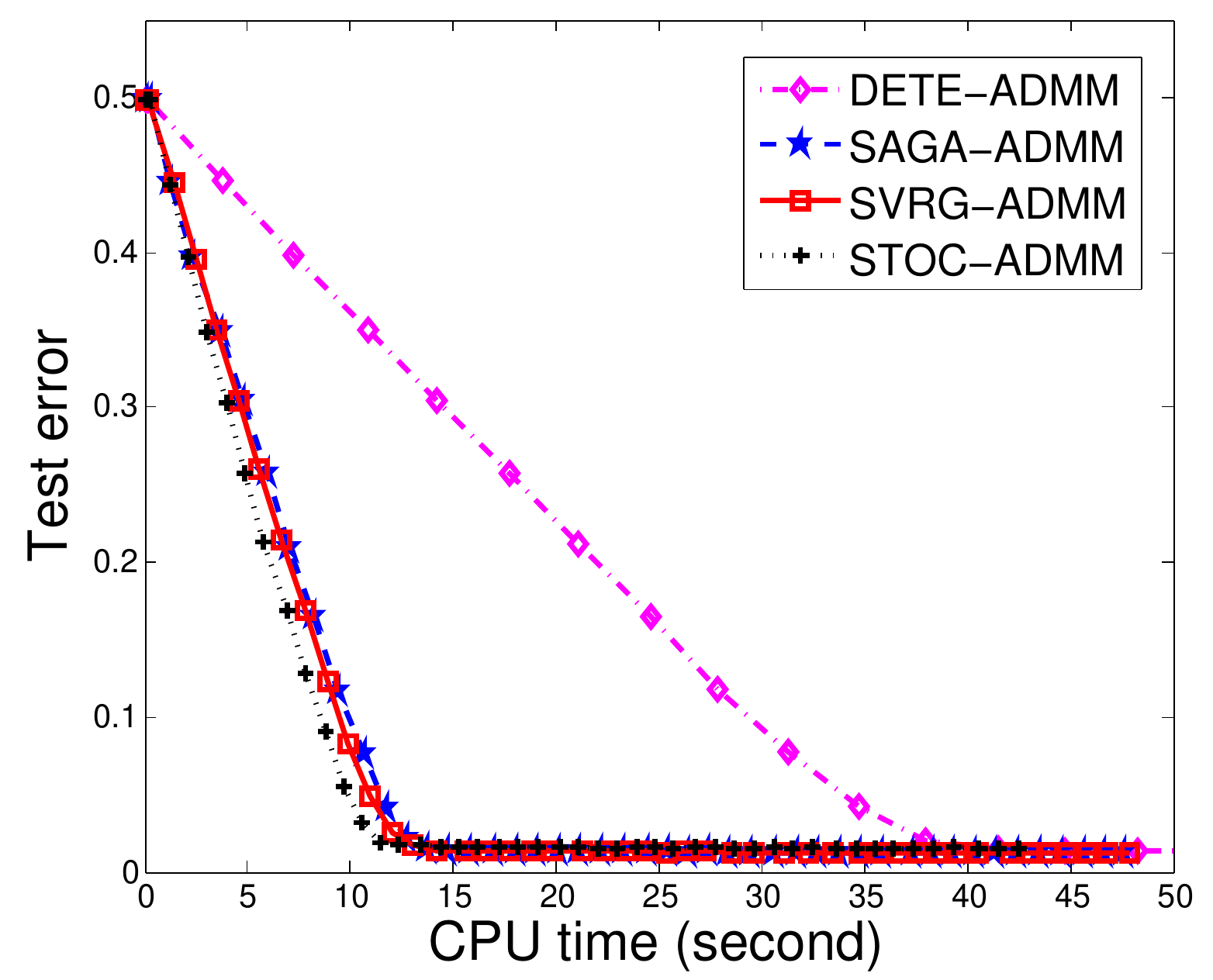}}
\caption{Test error \emph{versus} CPU time on the simulated \emph{nonconvex} model with \textbf{overlapping group Lasso}.}
\label{fig:4}
\vspace{-1em}
\end{figure*}

Figs. \ref{fig:3} and \ref{fig:4} show that both the objective values and test loss of these
stochastic ADMMs faster decrease than those of the deterministic ADMM,
as CPU time consumed increases.
In particular, though the nonconvex \emph{STOC-ADMM} uses a fixed step size $\eta$,
it shows good performance in the nonconvex optimization with overlapping group lasso regularization,
and is comparable with both the nonconvex \emph{SVRG-ADMM} and \emph{SAGA-ADMM}.

\subsection{ Real Data }
In the subsection, we compare the performances in some real data.
Specifically, we perform the binary classification task and multitask learning, respectively.

\subsubsection{ Graph-guided Fused Lasso }
Here we perform the binary classification task with the graph-guided fused lasso
penalty function as in \eqref{eq:80}.
In the experiment, we use some publicly available datasets\footnote{\emph{20news} is from the website (https://cs.nyu.edu/~roweis/data.html);
\emph{a9a}, \emph{w8a}, \emph{ijcnn1} and \emph{covtype.binary}
are from the LIBSVM website (www.csie.ntu.edu.tw/~cjlin/libsvmtools/datasets/).},
which are summarized in Table \ref{tab:3}.
In the algorithms, we use the same initial solution $x_0$ from the standard normal distribution
and choose a fixed step size $\eta = 1$.
In the problem \eqref{eq:80}, we fix the parameter $\nu=10^{-5}$.
In addition, we choose the mini-batch size $M=100$ in these stochastic algorithms.
The following experimental results are averaged over 10 repetitions.

\begin{figure*}[htbp]
\centering
\subfigure[\emph{20news}]{\includegraphics[width=0.19\textwidth]{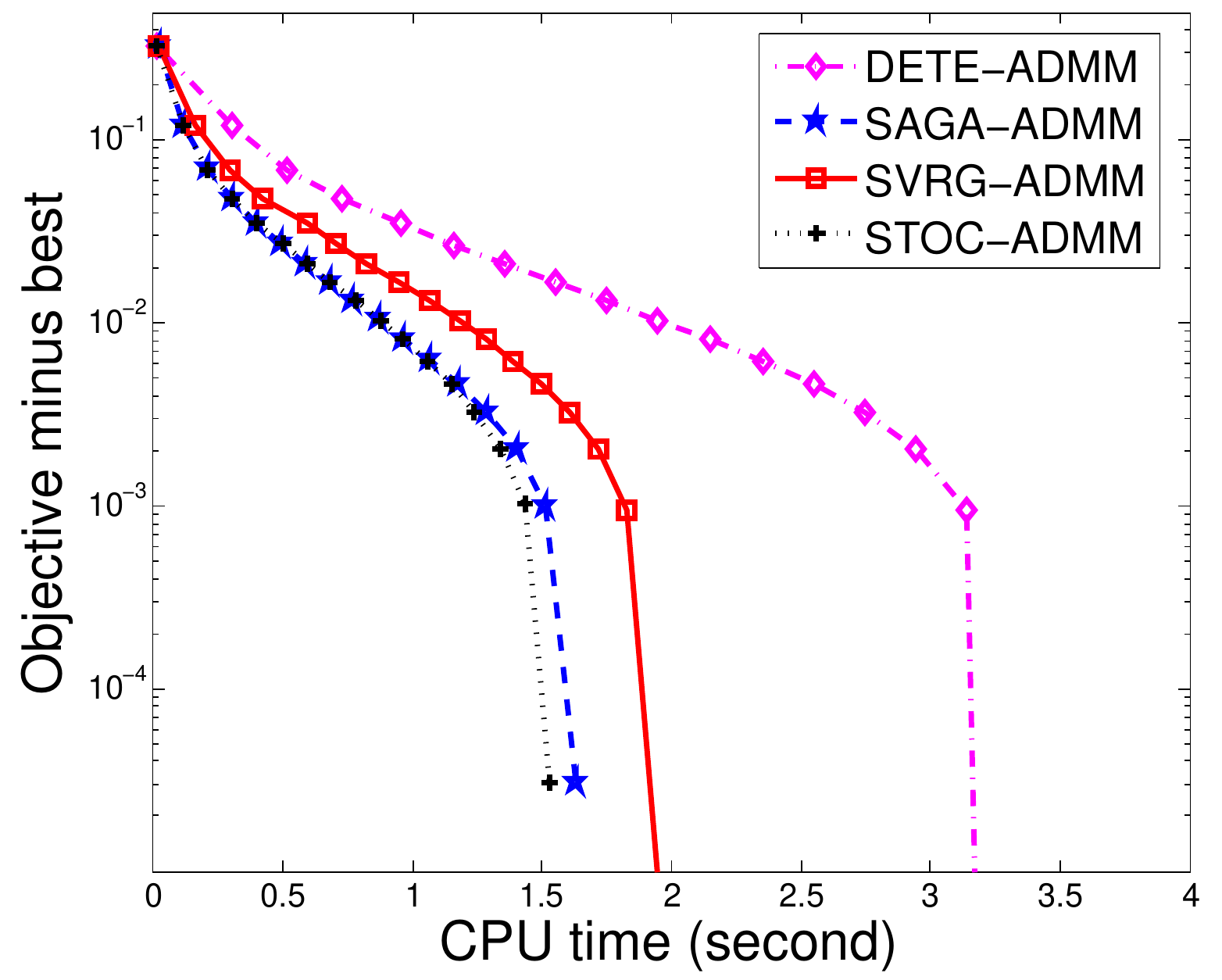}}
\subfigure[\emph{a9a}]{\includegraphics[width=0.19\textwidth]{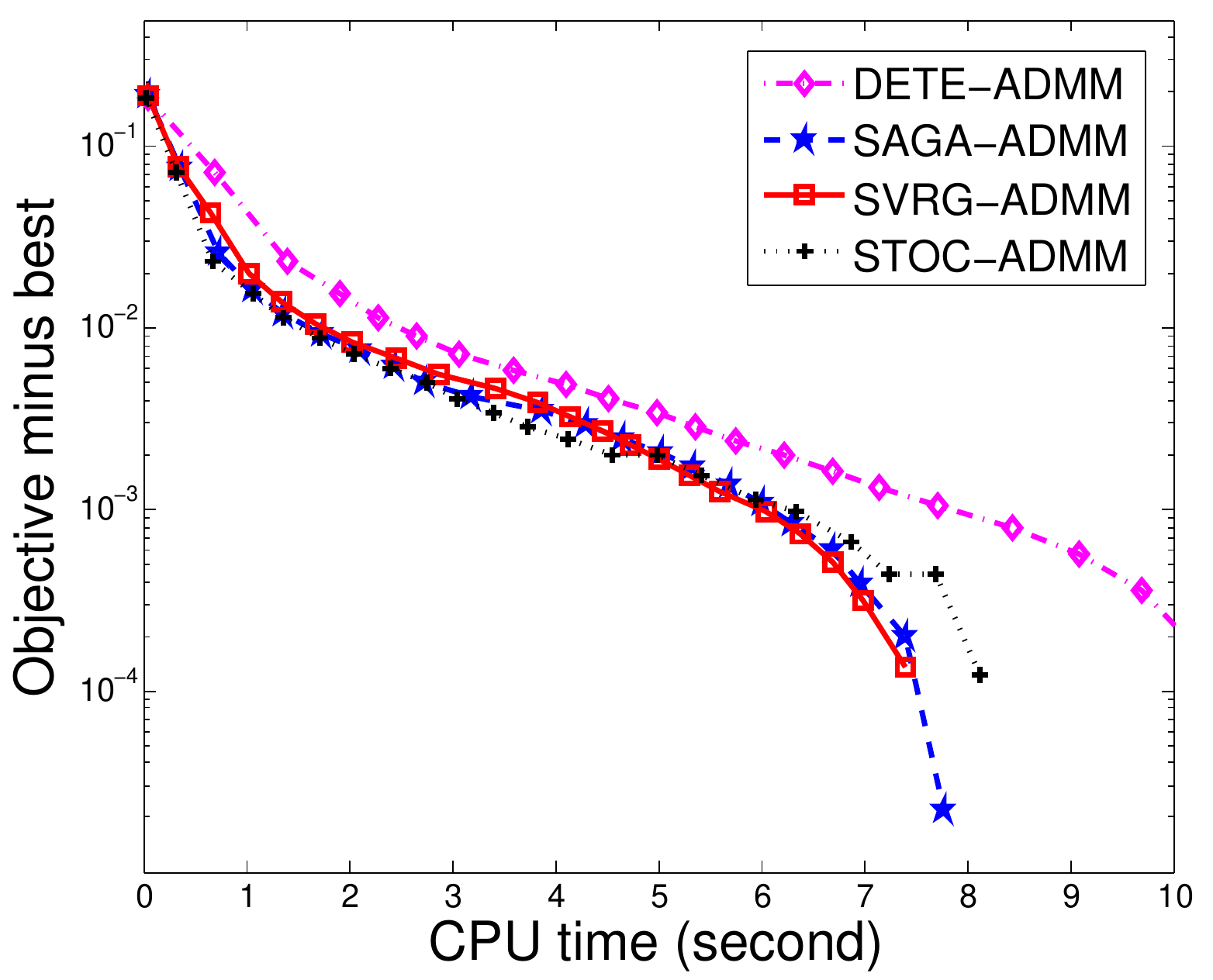}}
\subfigure[\emph{w8a}]{\includegraphics[width=0.19\textwidth]{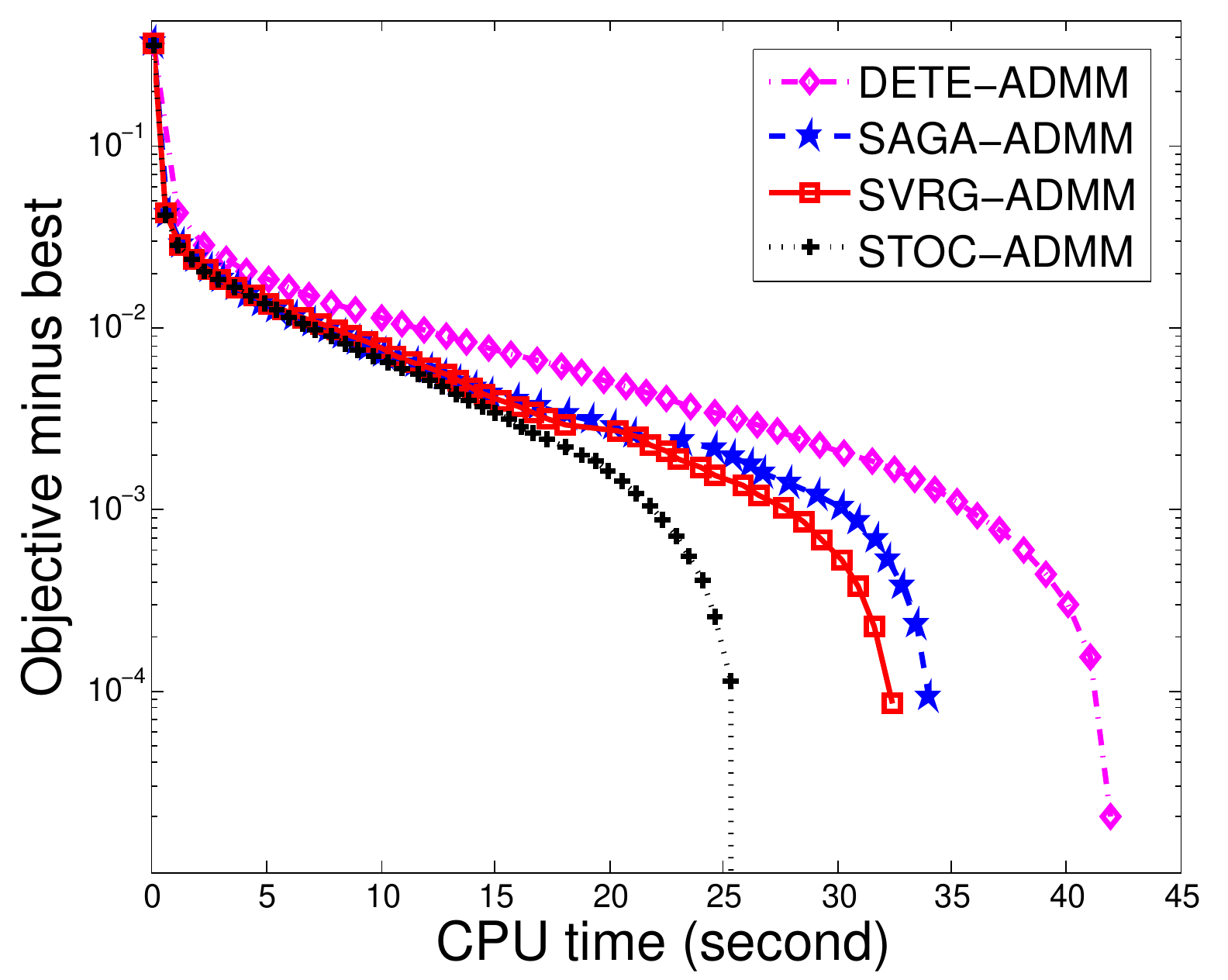}}
\subfigure[\emph{ijcnn1}]{\includegraphics[width=0.19\textwidth]{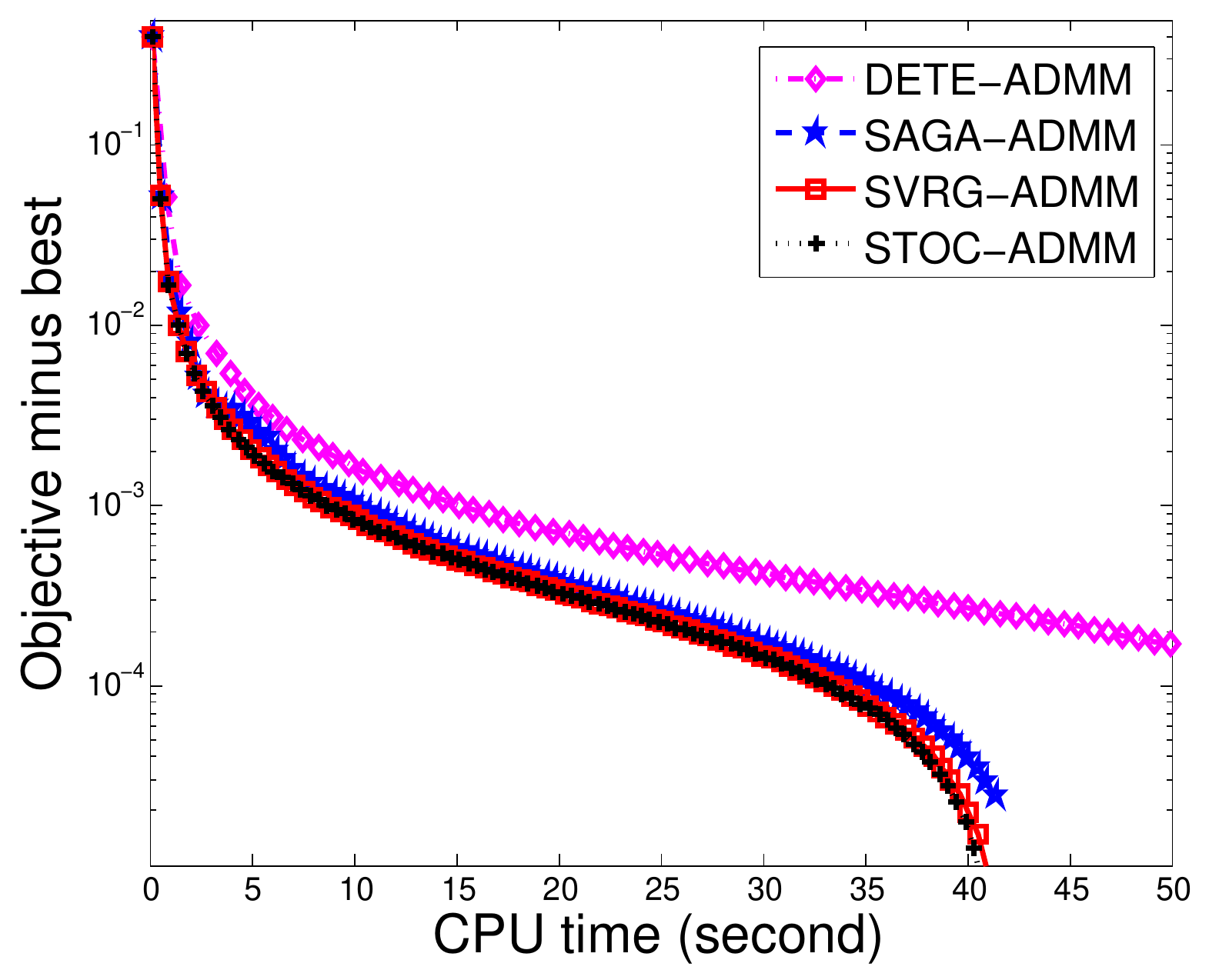}}
\subfigure[\emph{covtype.binary}]{\includegraphics[width=0.19\textwidth]{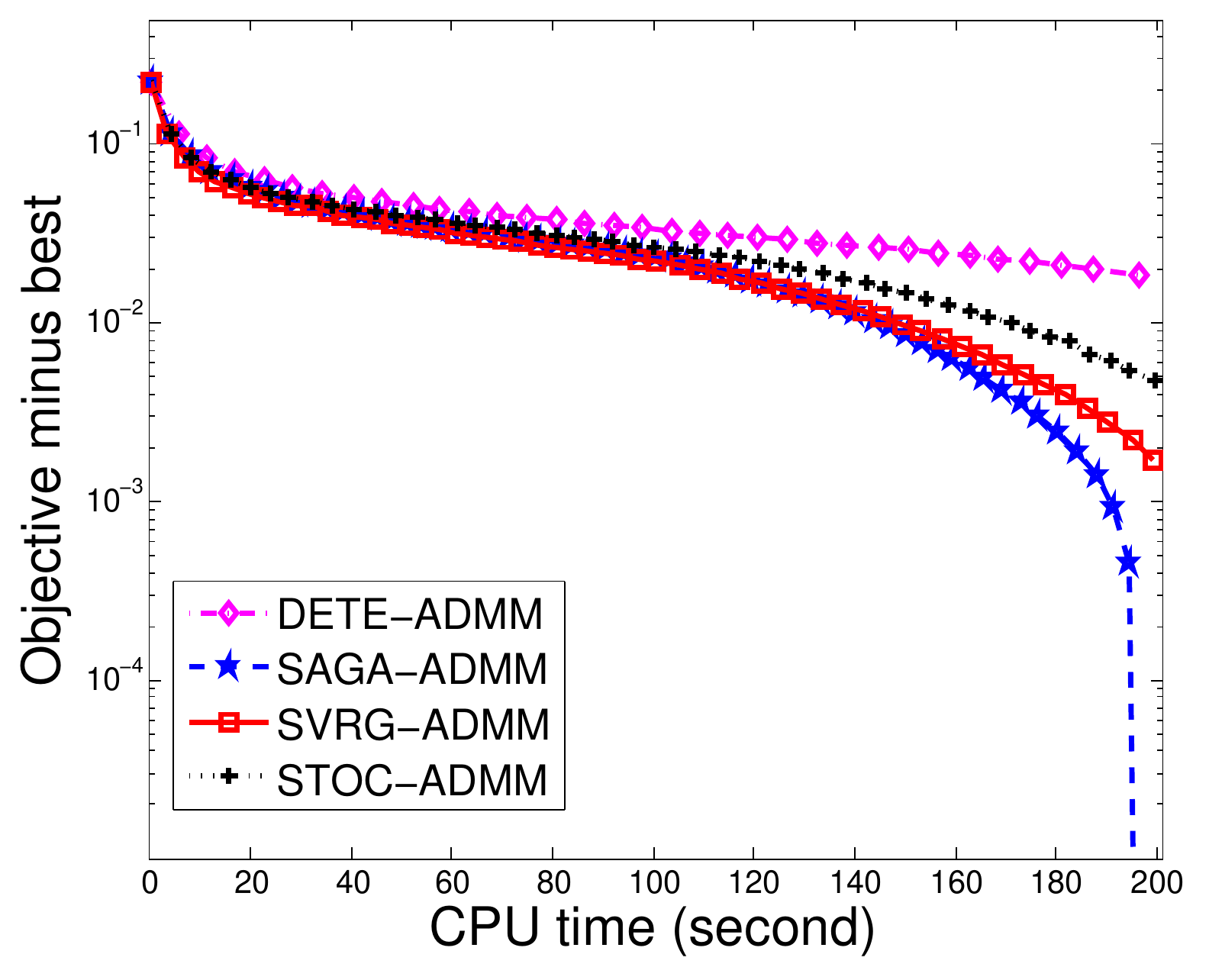}}
\caption{Objective value \emph{versus} CPU time of the \emph{nonconvex} \textbf{graph-guided binary classification model} on some real datasets.}
\label{fig:5}
\vspace{-1em}
\end{figure*}

\begin{figure*}[htbp]
\centering
\subfigure[\emph{20news}]{\includegraphics[width=0.19\textwidth]{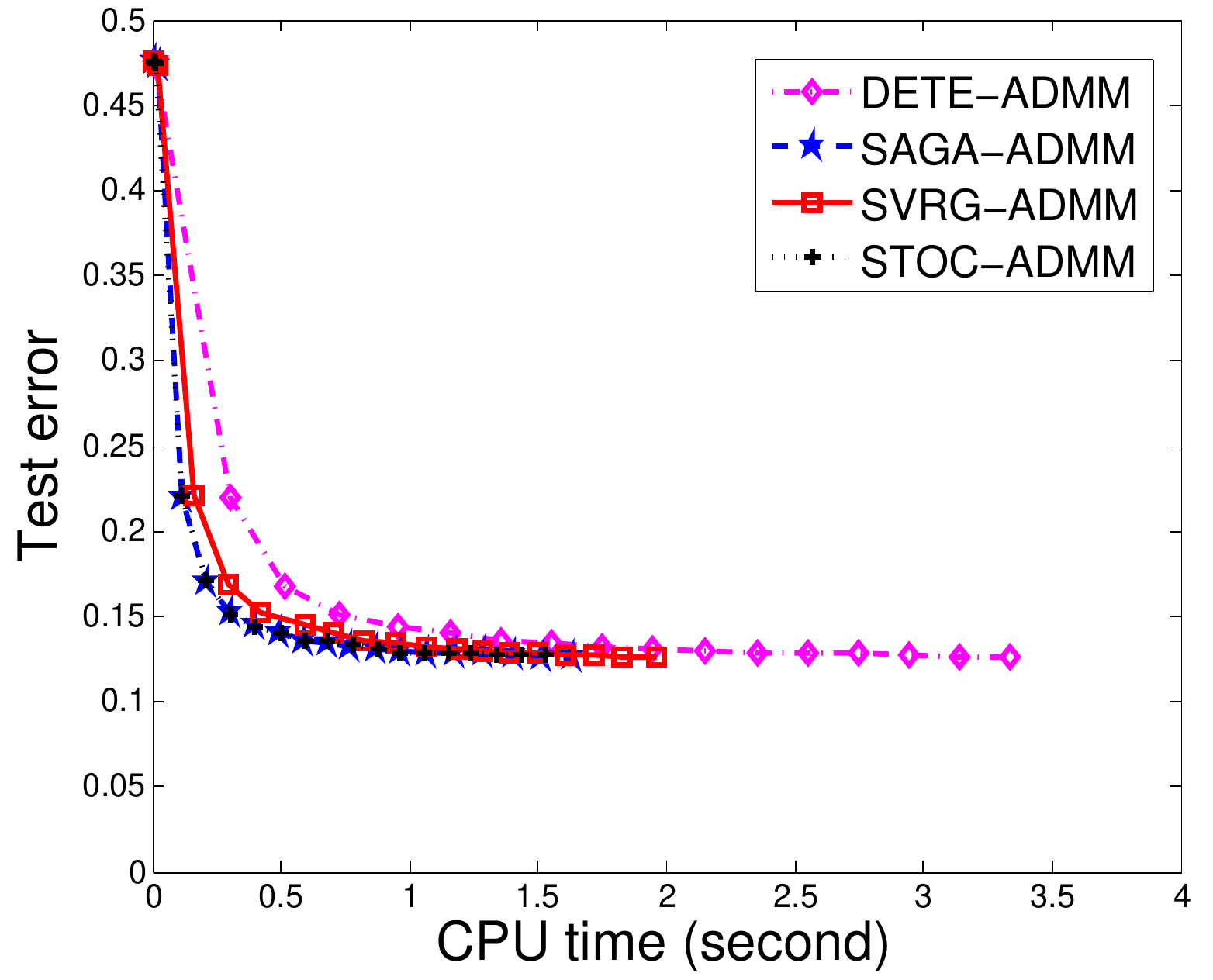}}
\subfigure[\emph{a9a}]{\includegraphics[width=0.19\textwidth]{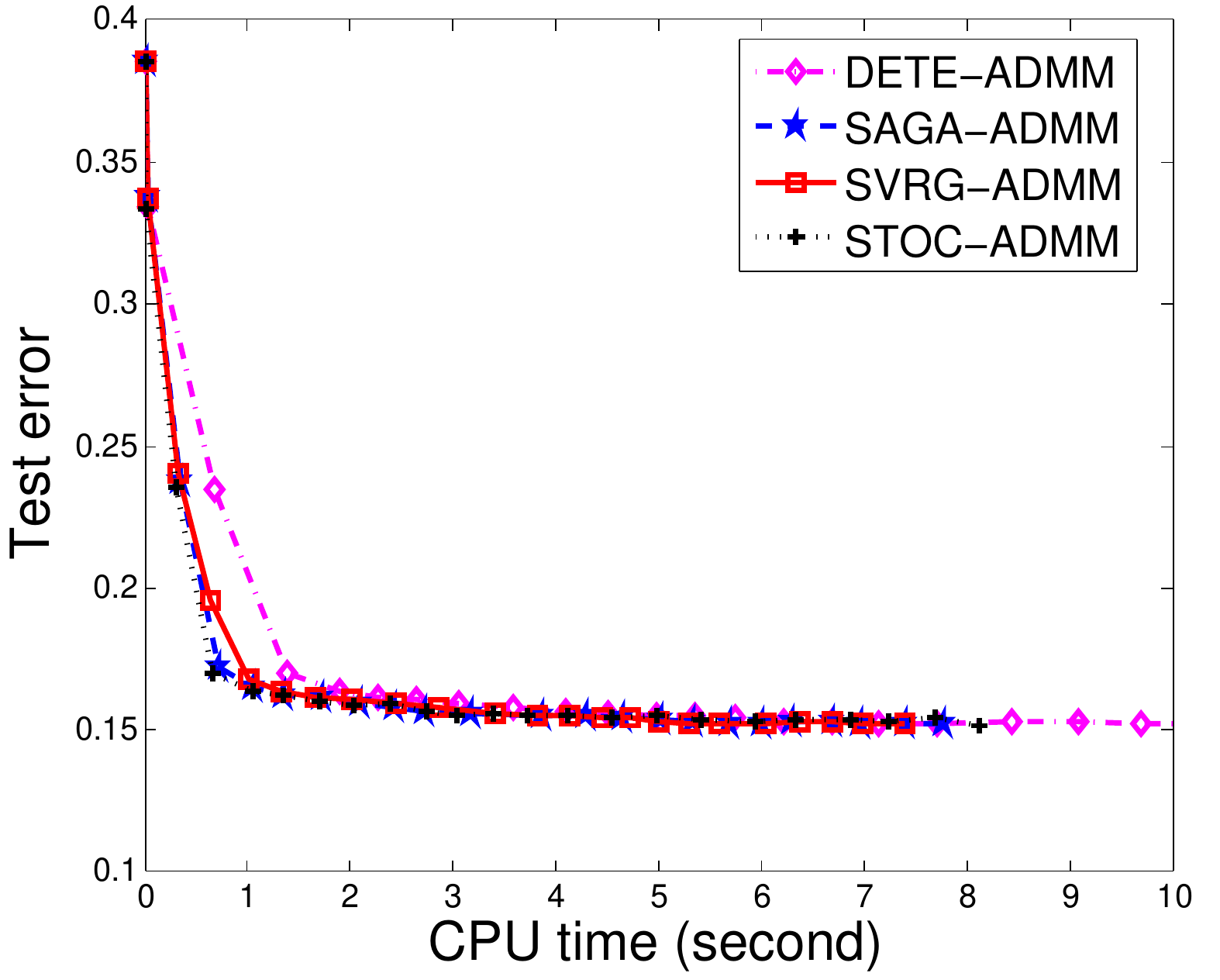}}
\subfigure[\emph{w8a}]{\includegraphics[width=0.19\textwidth]{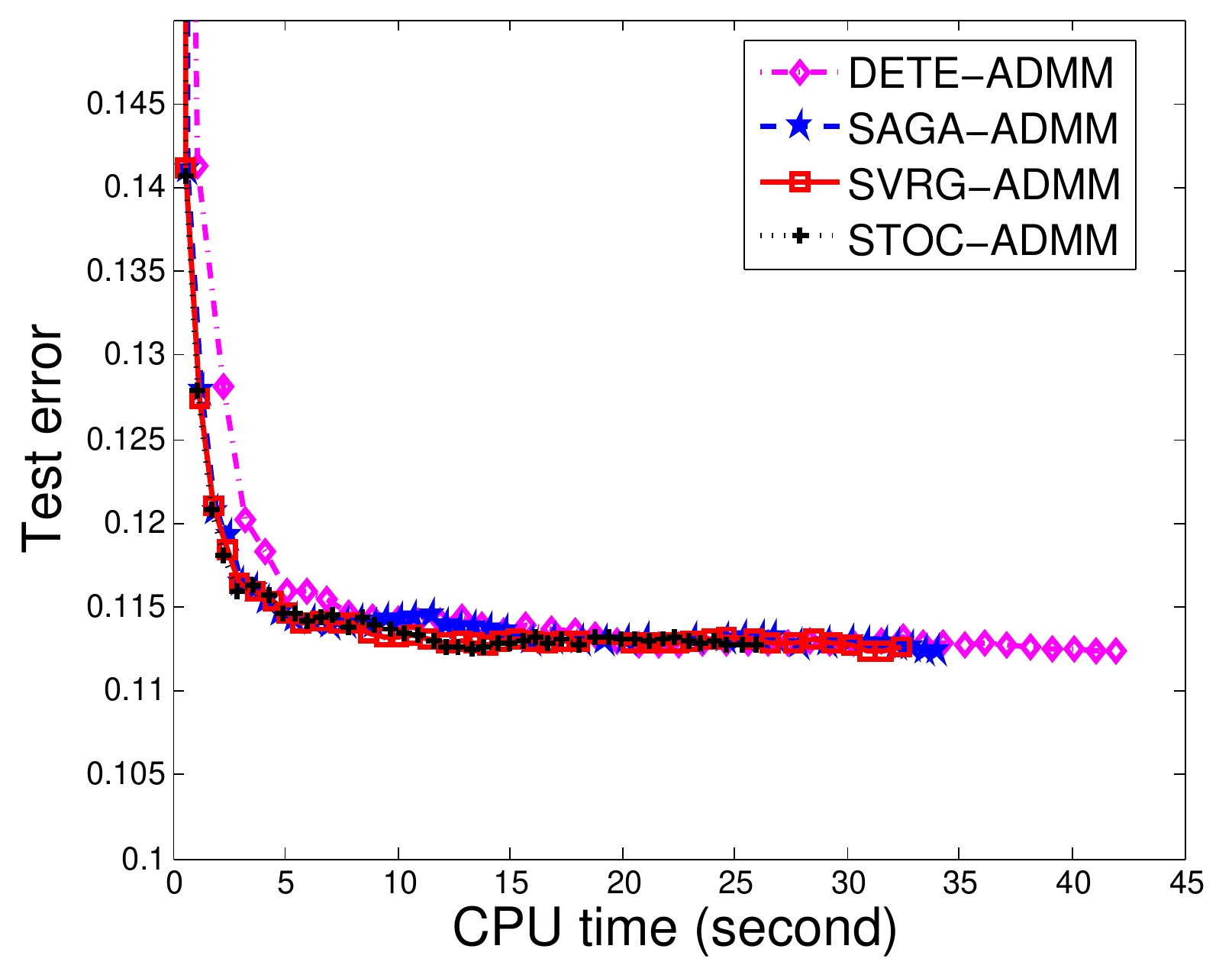}}
\subfigure[\emph{ijcnn1}]{\includegraphics[width=0.19\textwidth]{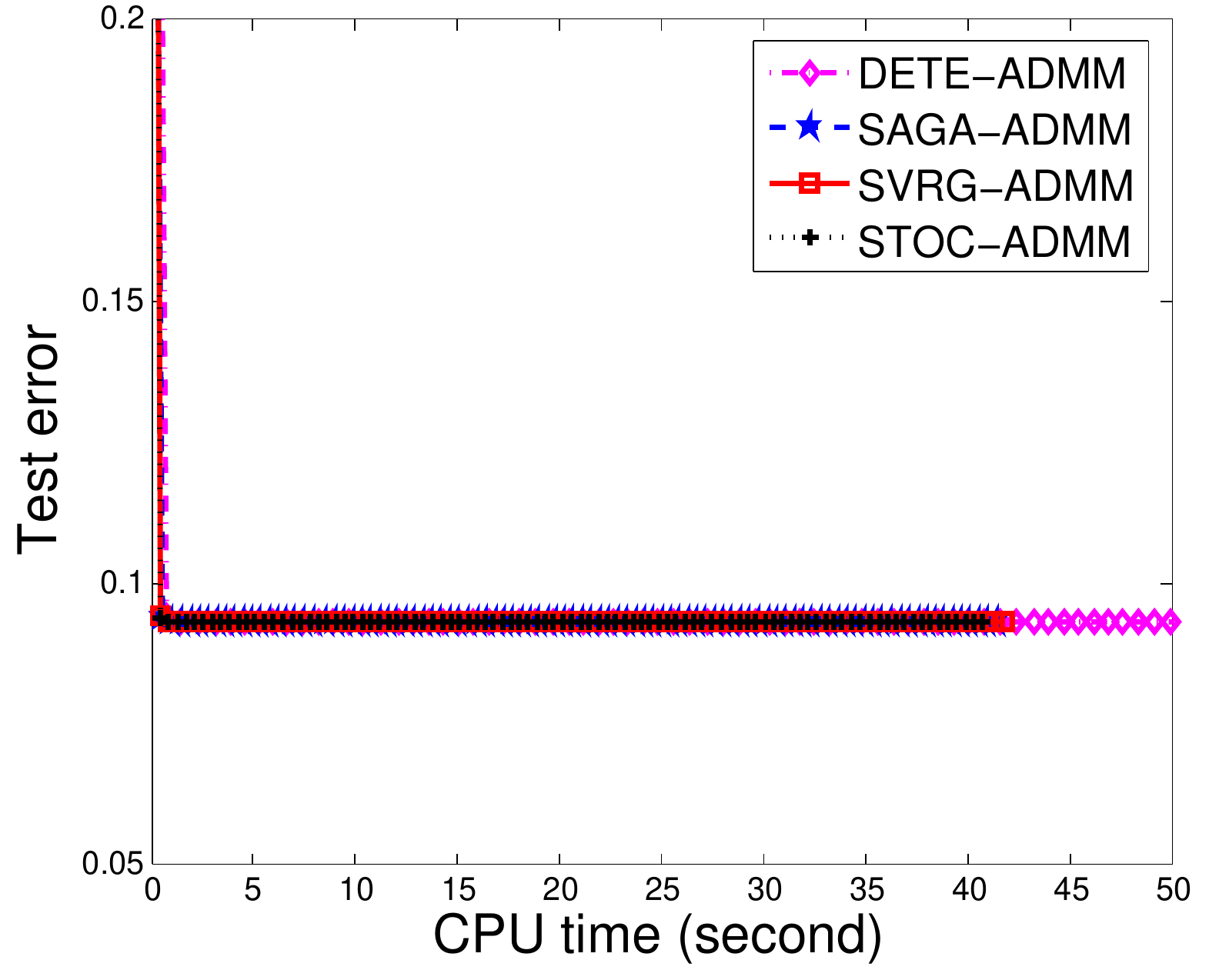}}
\subfigure[\emph{covtype.binary}]{\includegraphics[width=0.19\textwidth]{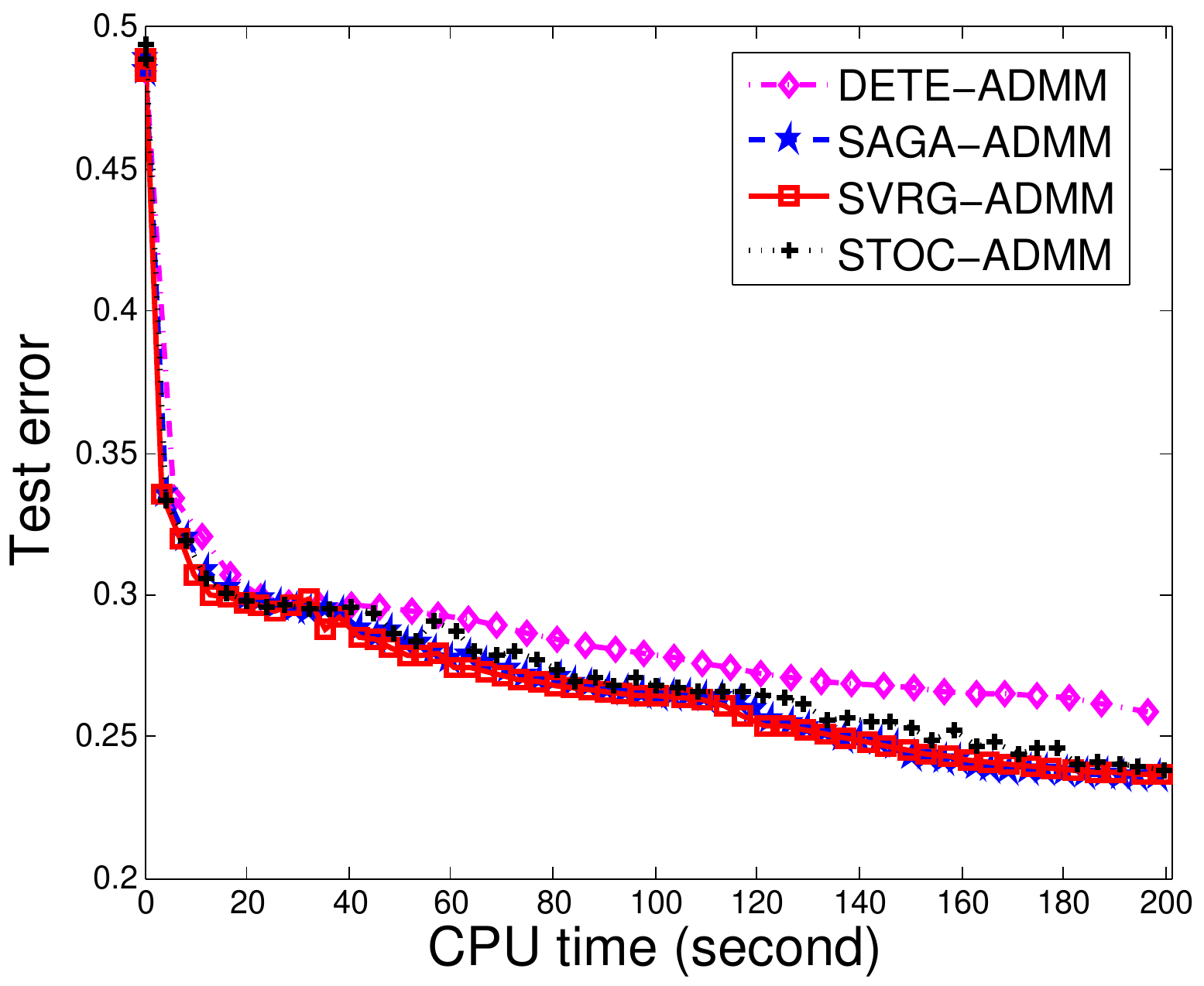}}
\caption{Test loss \emph{versus} CPU time of the \emph{nonconvex} \textbf{graph-guided binary classification model} on some real datasets.}
\label{fig:6}
\vspace{-1em}
\end{figure*}

Figs. \ref{fig:5} and \ref{fig:6} show that both the objective values and test loss of these
stochastic ADMMs faster decrease than those of the deterministic ADMM,
as CPU time consumed increases.
In particular, though the nonconvex \emph{STOC-ADMM} uses a fixed step size $\eta$,
it shows good performance in the nonconvex optimization with graph-guided fused lasso regularization,
and is comparable with both the nonconvex \emph{SVRG-ADMM} and \emph{SAGA-ADMM}.
Due to that \emph{ijcnn1} is a severely imbalanced data set (i.e., includes large negative samples),
the first iteration solution of all algorithms shows good performance on the testing error.

\begin{table}
  \centering
  \caption{Real data for graph-guided fused lasso } \label{tab:3}
  \begin{tabular}{c|c|c|c|c}
  \hline
  datasets & $\#training$ & $\#test$ & $\#features$ & $\#classes$ \\ \hline
  \emph{20news} & 8,121 & 8121 &  100 & 2\\
  \emph{a9a}   & 16,281 & 16,280 & 123 & 2\\
  \emph{w8a} & 32,350 & 32,350 & 300 & 2\\
  \emph{ijcnn1} & 63,351 & 63,351 & 22 & 2\\
  \emph{covtype.binary} & 290,506 & 290,506 & 54 & 2\\
  \hline
  \end{tabular}
\end{table}

\subsubsection{ Multi-task Learning }
Here we perform the multi-task learning with both sparse and low-rank penalty functions.
Specifically, given a set of training samples $(a_i,b_i)_{i=1}^n$,
where $a_i\in R^d$ and $b_i \in \{1,2,\cdots,m\}$. Let $\bar{b}_{i,c}=1$ if $b_i=c\in \{1,2,\cdots,m\}$,
and $\bar{b}_{i,c}=0$ otherwise.
Then we solve the following nonconvex problem
\begin{align} \label{eq:81}
 \min_{X\in R^{m\times d}} \frac{1}{n}\sum_{i=1}^n f_i(X) + \nu_1\sum_{i,j=1}^{m,d}\kappa(|X_{i,j}|) + \nu_2\|X\|_*,
\end{align}
where $f_i(X) = \log(\sum_{c=1}^m\exp(X_{c,.}^Ta_i)) - \sum_{c=1}^m\bar{b}_{i,c}X_{c,.}^Ta_i$ is a multinomial logistic loss
function, $\kappa(\alpha)=\beta\log(1+\frac{\alpha}{\theta})$ is the
nonconvex log-sum penalty function \cite{candes2008enhancing}, and $\|X\|_*$ denotes the nuclear norm of matrix $X$.
Here $\nu_1$ and $\nu_2$ are nonegative regularization parameters.
Following \cite{yao2016efficient}, we can transform the problem \eqref{eq:81} into the following problem
\begin{align} \label{eq:82}
 \min_{X\in R^{m\times d}} \frac{1}{n}\sum_{i=1}^n \bar{f}_i(X) + \breve{g}(X)
\end{align}
where $\bar{f}_i(X) = f(X) + \nu_1\big( \sum_{i,j=1}^{m,d}\kappa(|X_{i,j}|)-\kappa_0\|X\|_1 \big)$,
$\breve{g}(X) = \nu_1\kappa_0\|X\|_1+\nu_2\|X\|_*$, and $\kappa_0=\kappa'(0)$. By Proposition 2.3 in \cite{yao2016efficient},
$\bar{f}_i(X)$ is nonconvex and smooth, and $\breve{g}(X)$ is nonsmooth and convex.
To solve the problem \eqref{eq:82} by using ADMMs,
we introduce an auxiliary variable $Y$ with
the constraint $X=Y$, and given $A=[I;I]$, then $\breve{g}(AX)= \nu_1\kappa_0\|X\|_1+\nu_2\|Y\|_*$.

\begin{table}
  \centering
  \caption{Real data for multitask learning } \label{tab:4}
  \begin{tabular}{c|c|c|c|c}
  \hline
  datasets & $\#training$ & $\#test$ & $\#features$ & $\#classes$ \\ \hline
  \emph{letter} & 7,500 & 7,500 & 16 & 26 \\
  \emph{sensorless}   & 29,255 & 29,254 & 48 & 11 \\
  \emph{mnist} & 30,000 & 30,000 &  780 & 10 \\
  \emph{covtype} & 290,506 & 290,506 & 54 & 7 \\
  \emph{mnist8m} & 4,050,000 & 4,050,000 & 780 & 10 \\
  \hline
  \end{tabular}
\end{table}

\begin{figure*}[htbp]
\centering
\subfigure[\emph{letter}]{\includegraphics[width=0.19\textwidth]{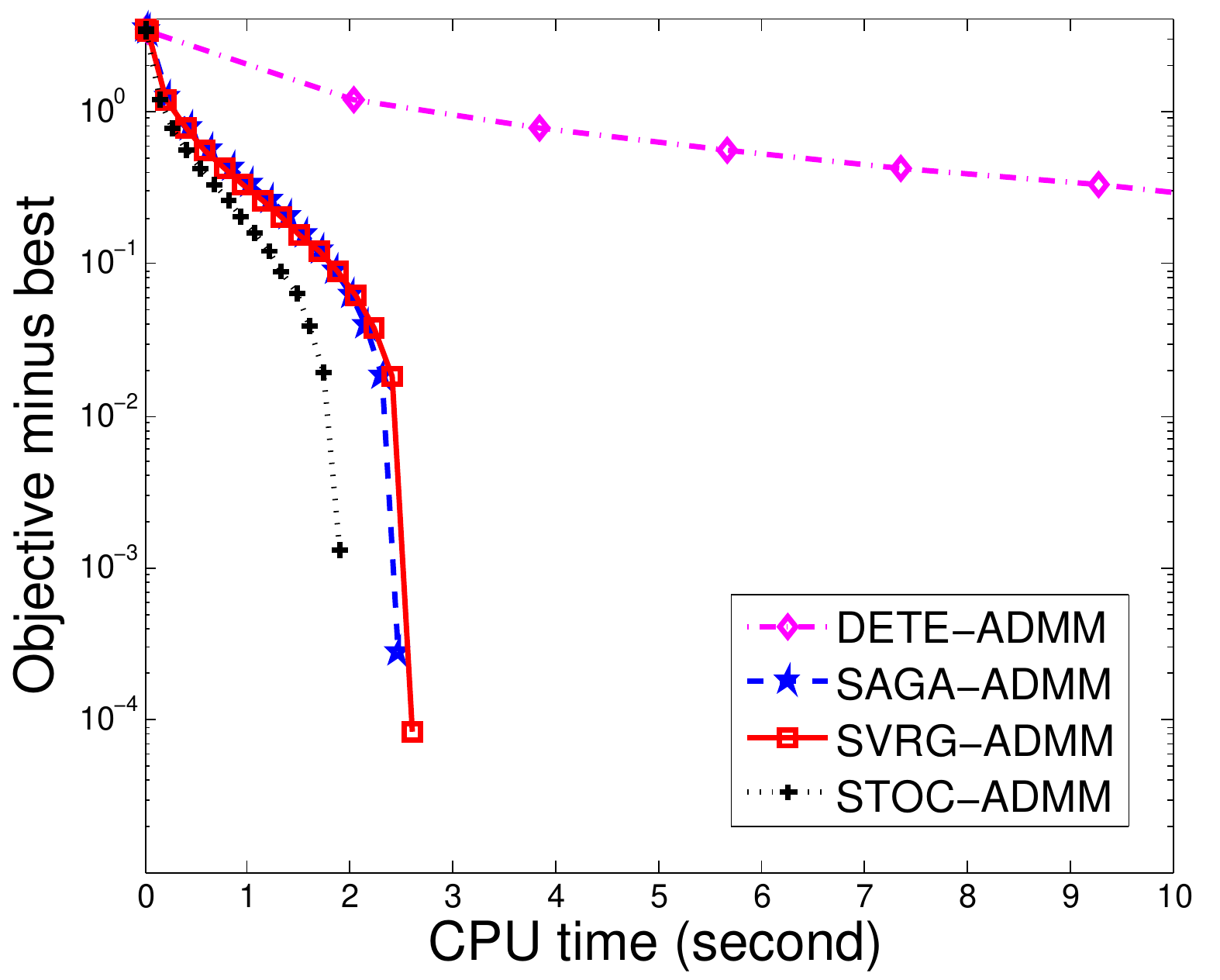}}
\subfigure[\emph{sensorless}]{\includegraphics[width=0.19\textwidth]{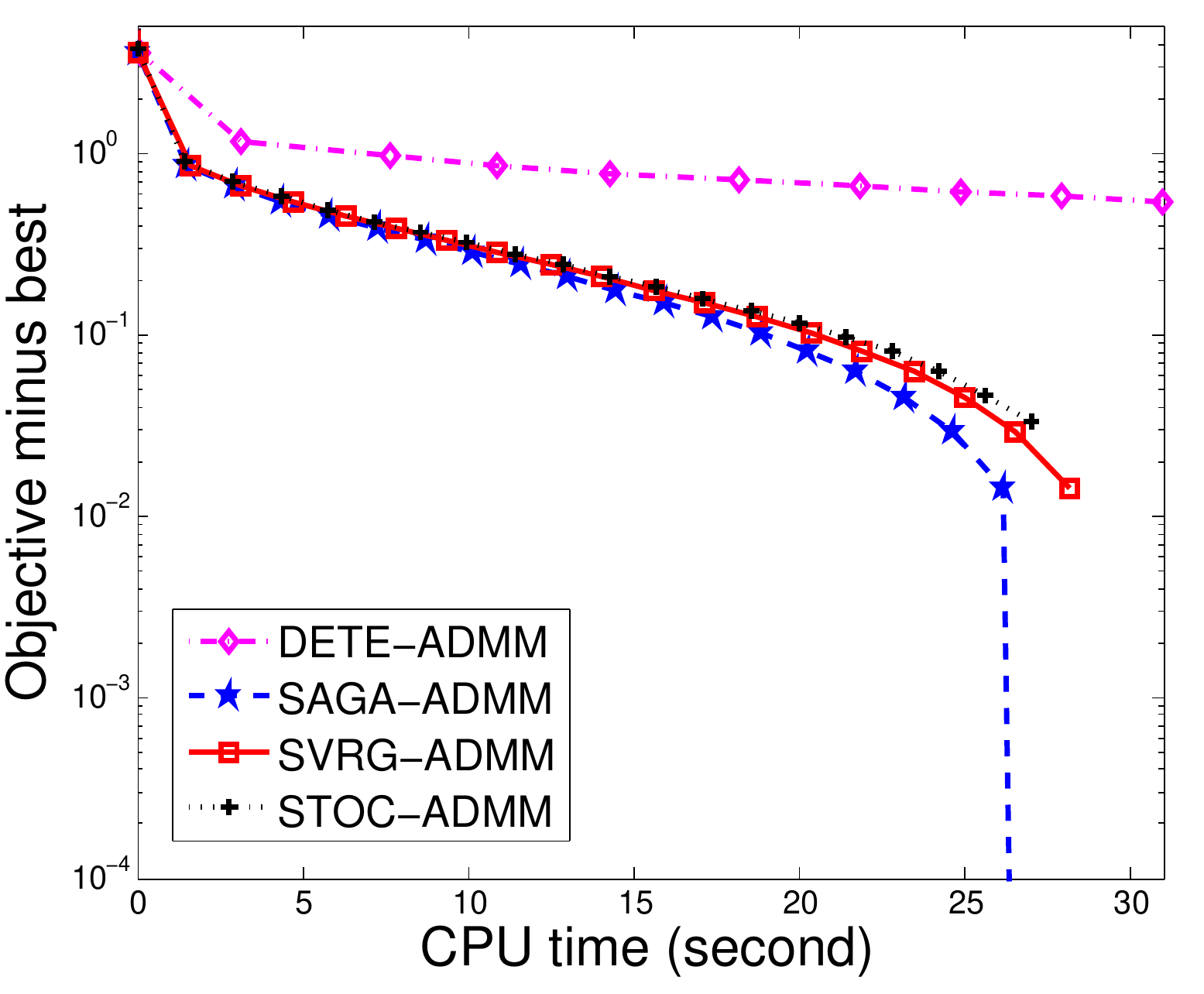}}
\subfigure[\emph{mnist}]{\includegraphics[width=0.19\textwidth]{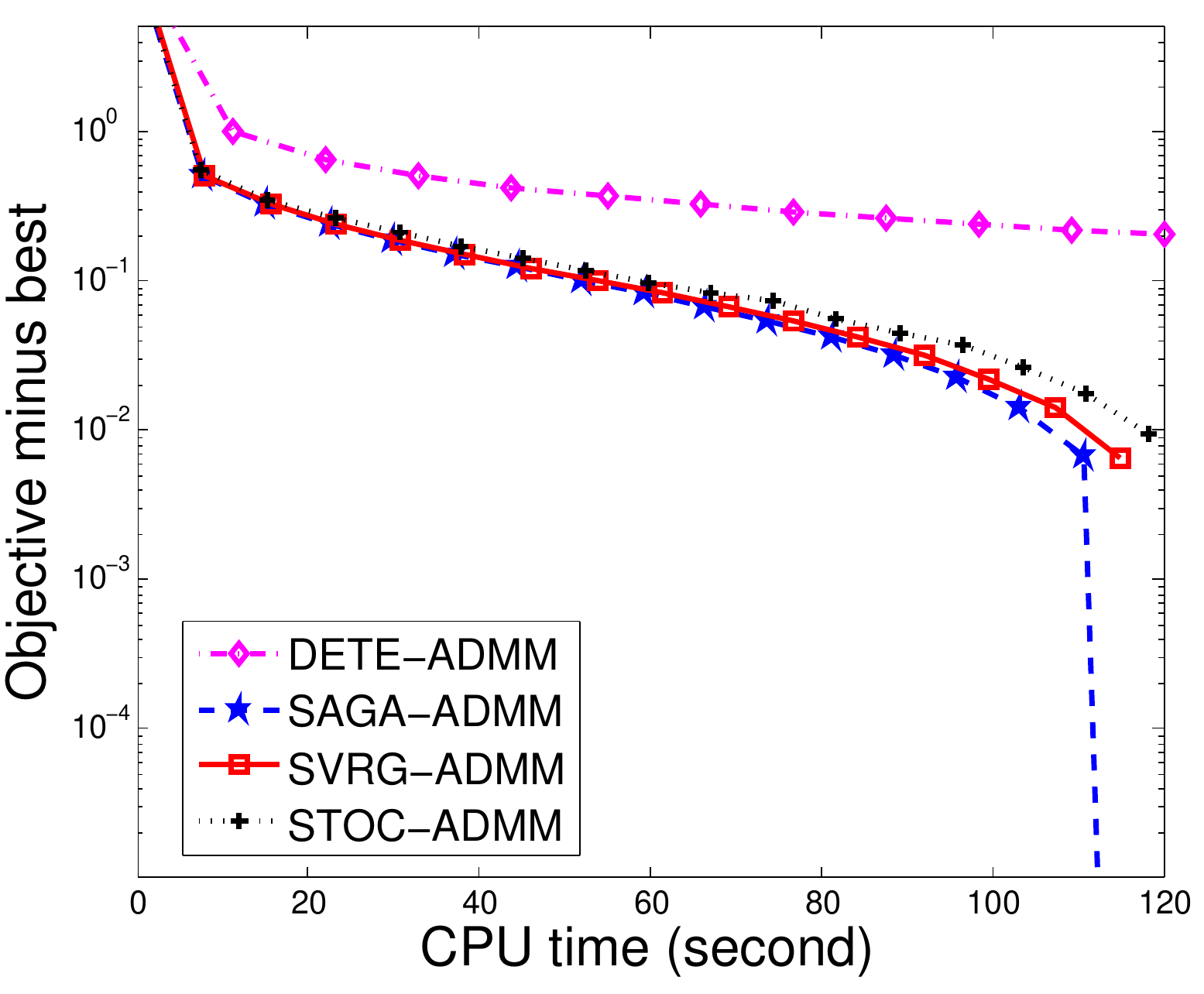}}
\subfigure[\emph{covtype}]{\includegraphics[width=0.19\textwidth]{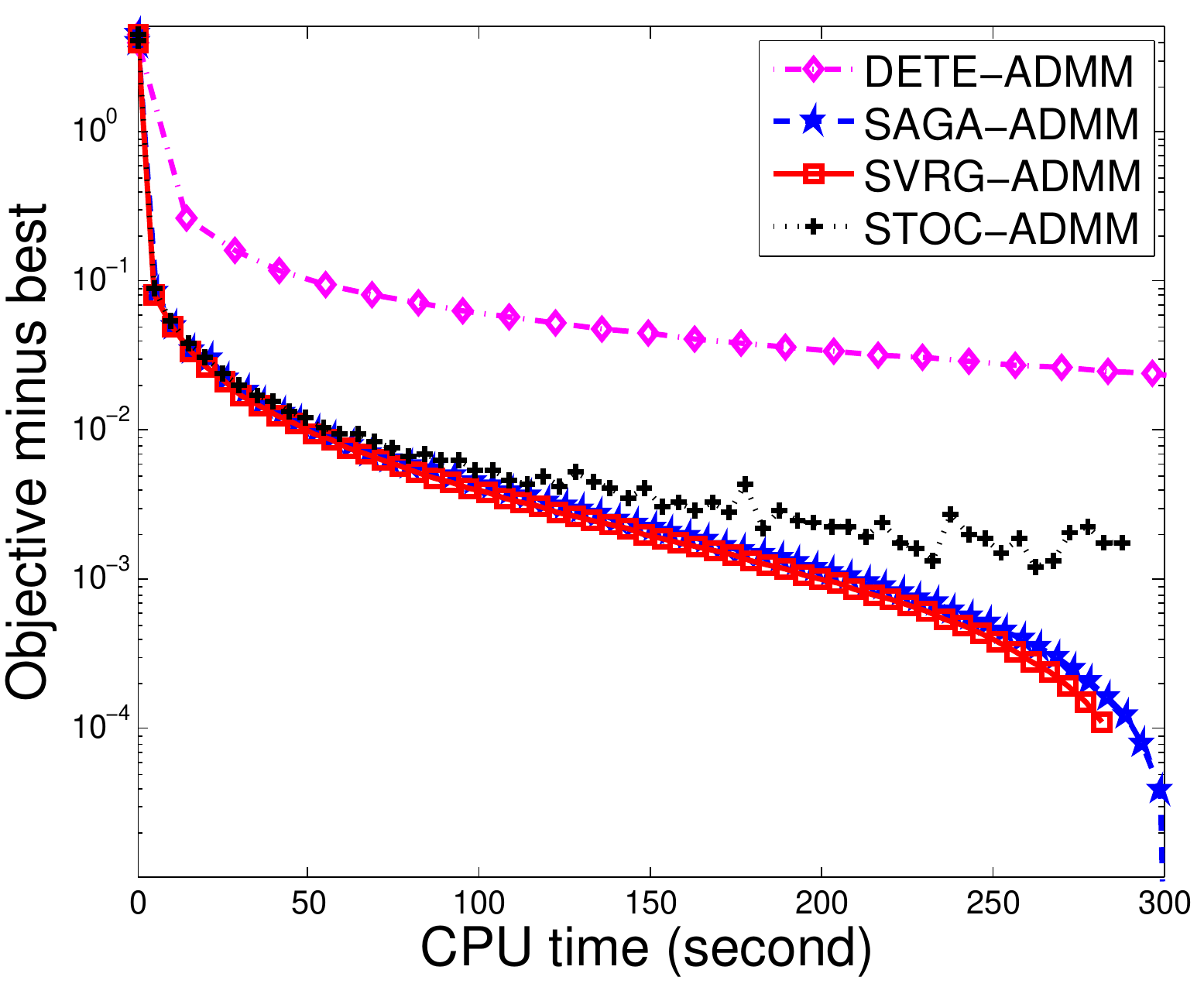}}
\subfigure[\emph{mnist8m}]{\includegraphics[width=0.19\textwidth]{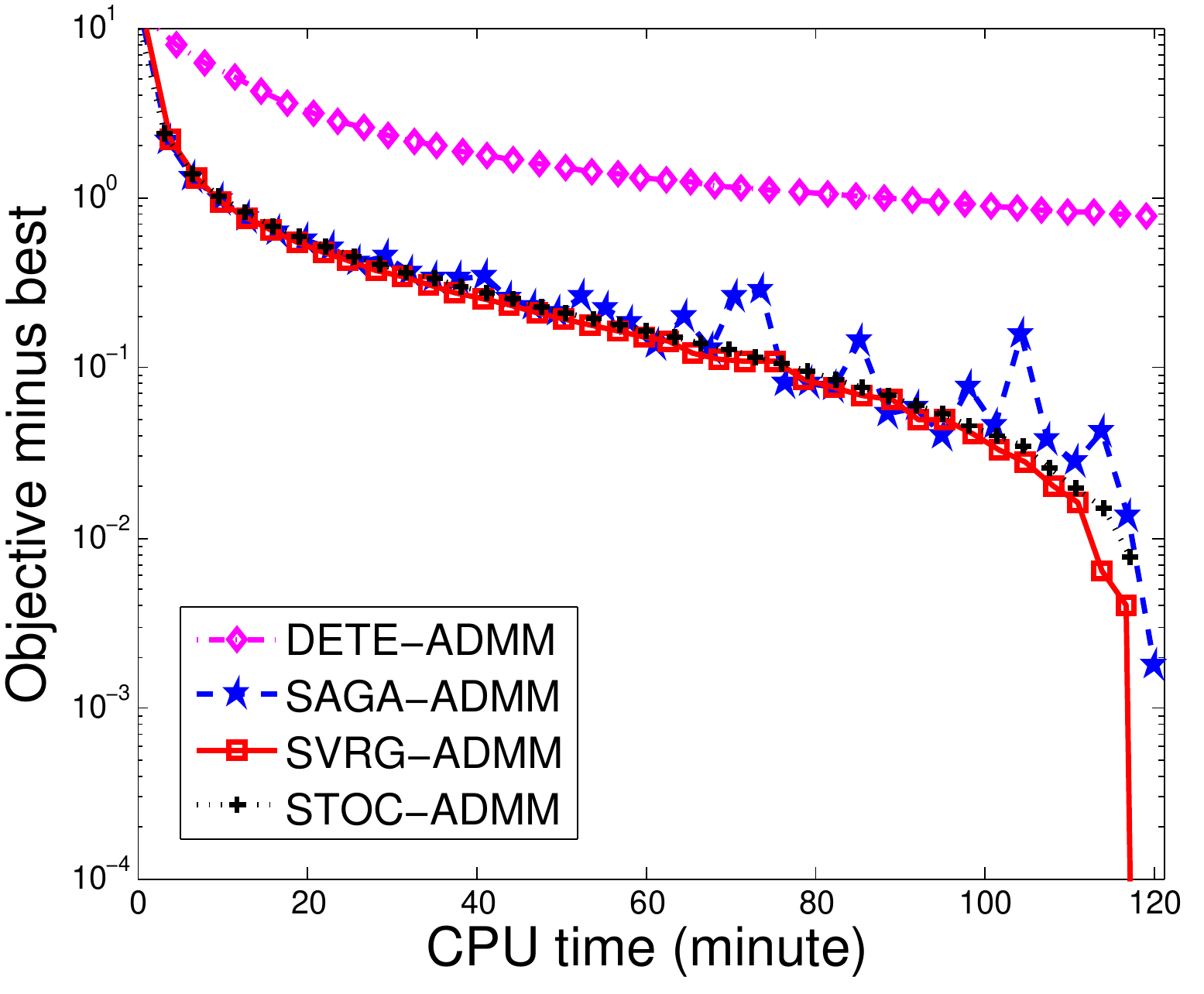}}
\caption{Objective value \emph{versus} CPU time of the \emph{nonconvex} \textbf{multi-task learning} on some real datasets.}
\label{fig:7}
\vspace{-1em}
\end{figure*}

\begin{figure*}[htbp]
\centering
\subfigure[\emph{letter}]{\includegraphics[width=0.19\textwidth]{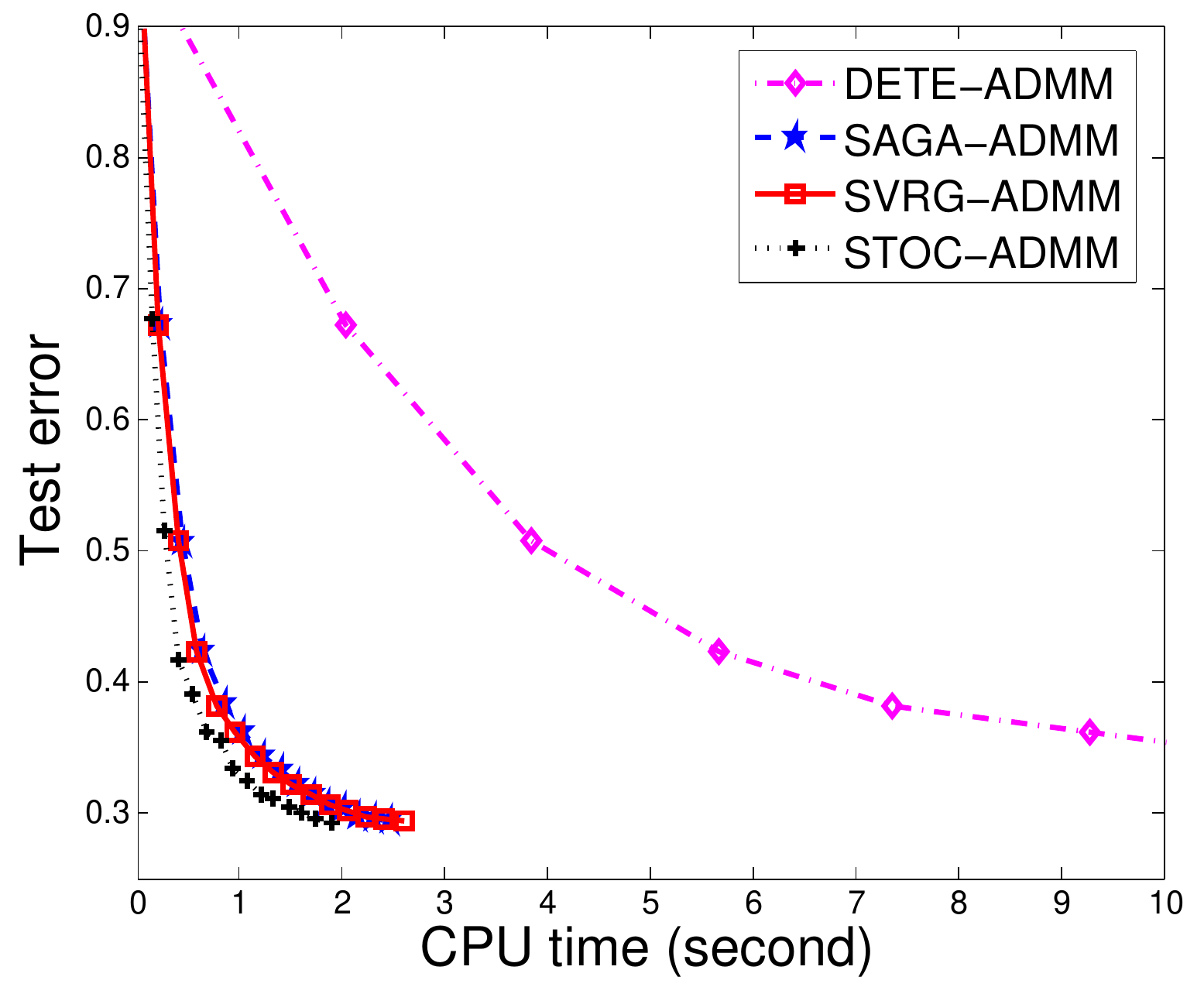}}
\subfigure[\emph{sensorless}]{\includegraphics[width=0.19\textwidth]{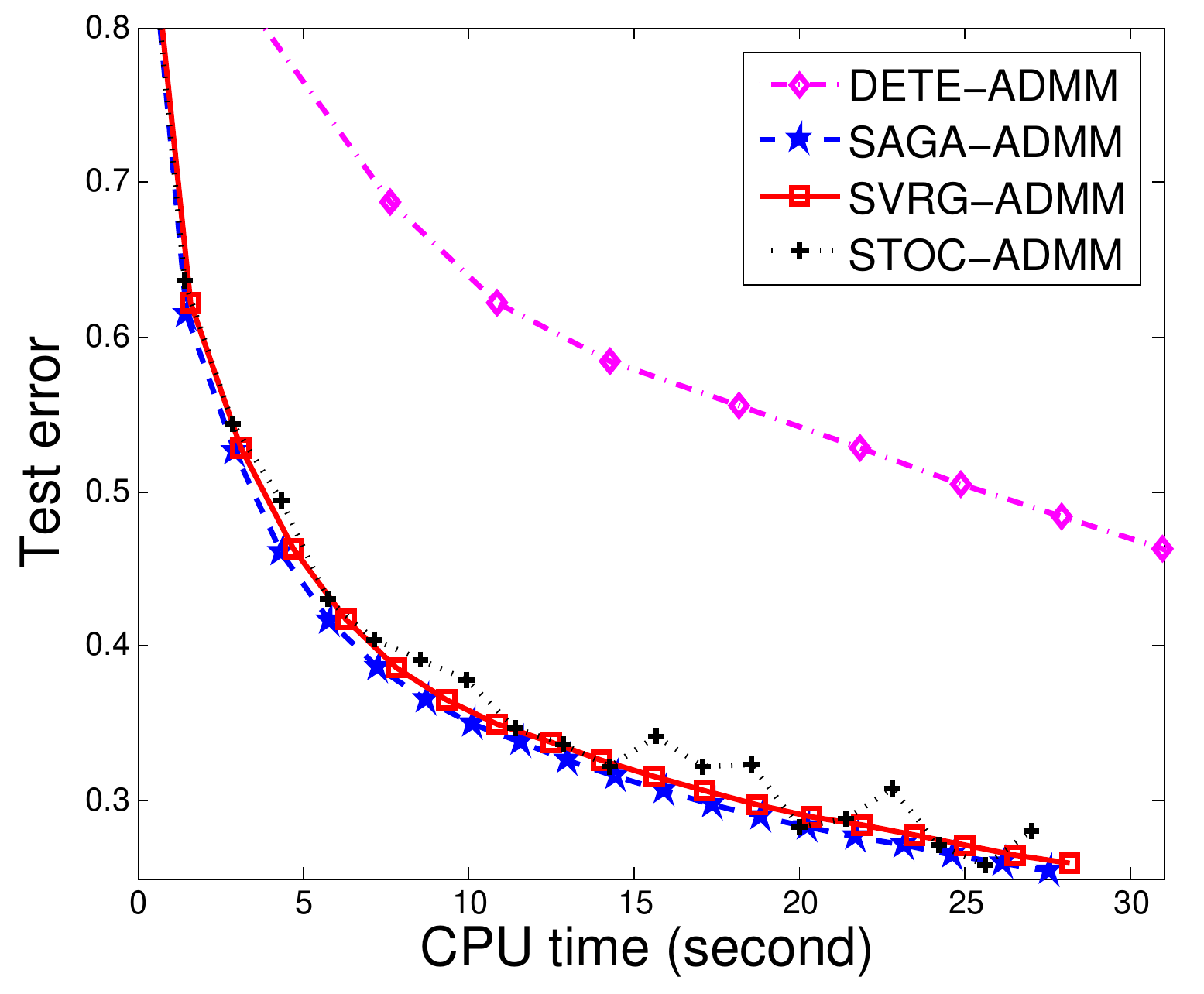}}
\subfigure[\emph{mnist}]{\includegraphics[width=0.19\textwidth]{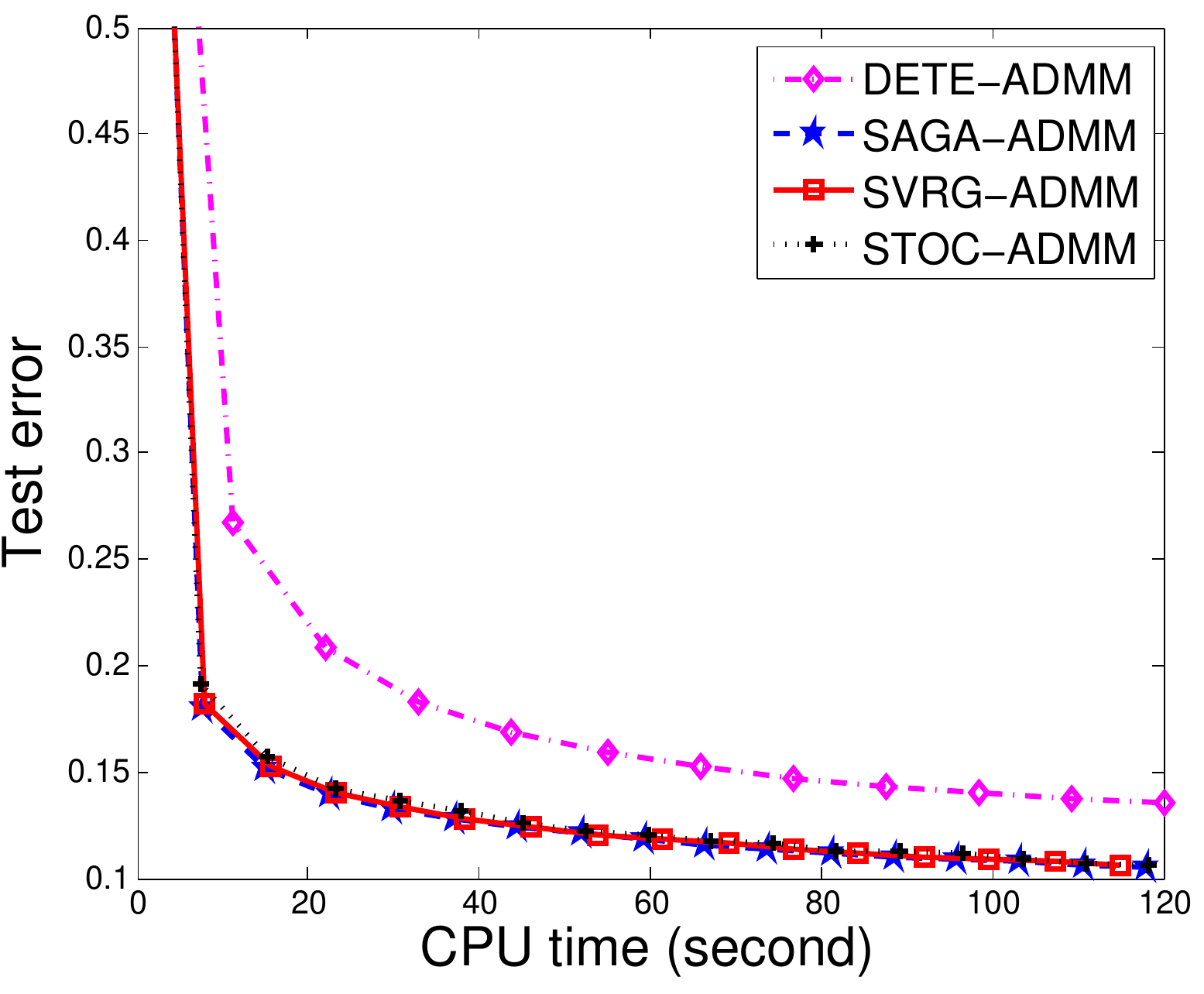}}
\subfigure[\emph{covtype}]{\includegraphics[width=0.19\textwidth]{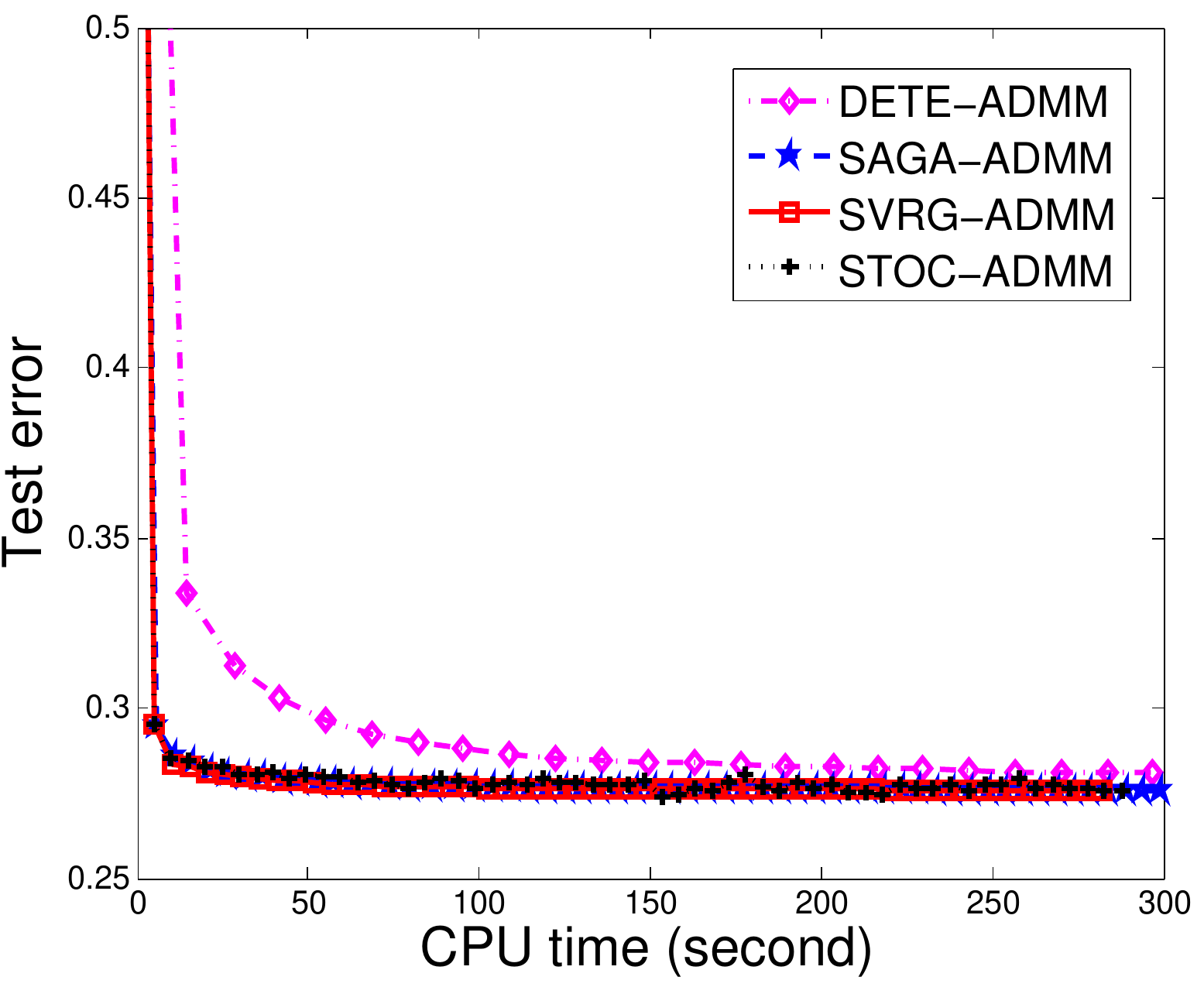}}
\subfigure[\emph{mnist8m}]{\includegraphics[width=0.19\textwidth]{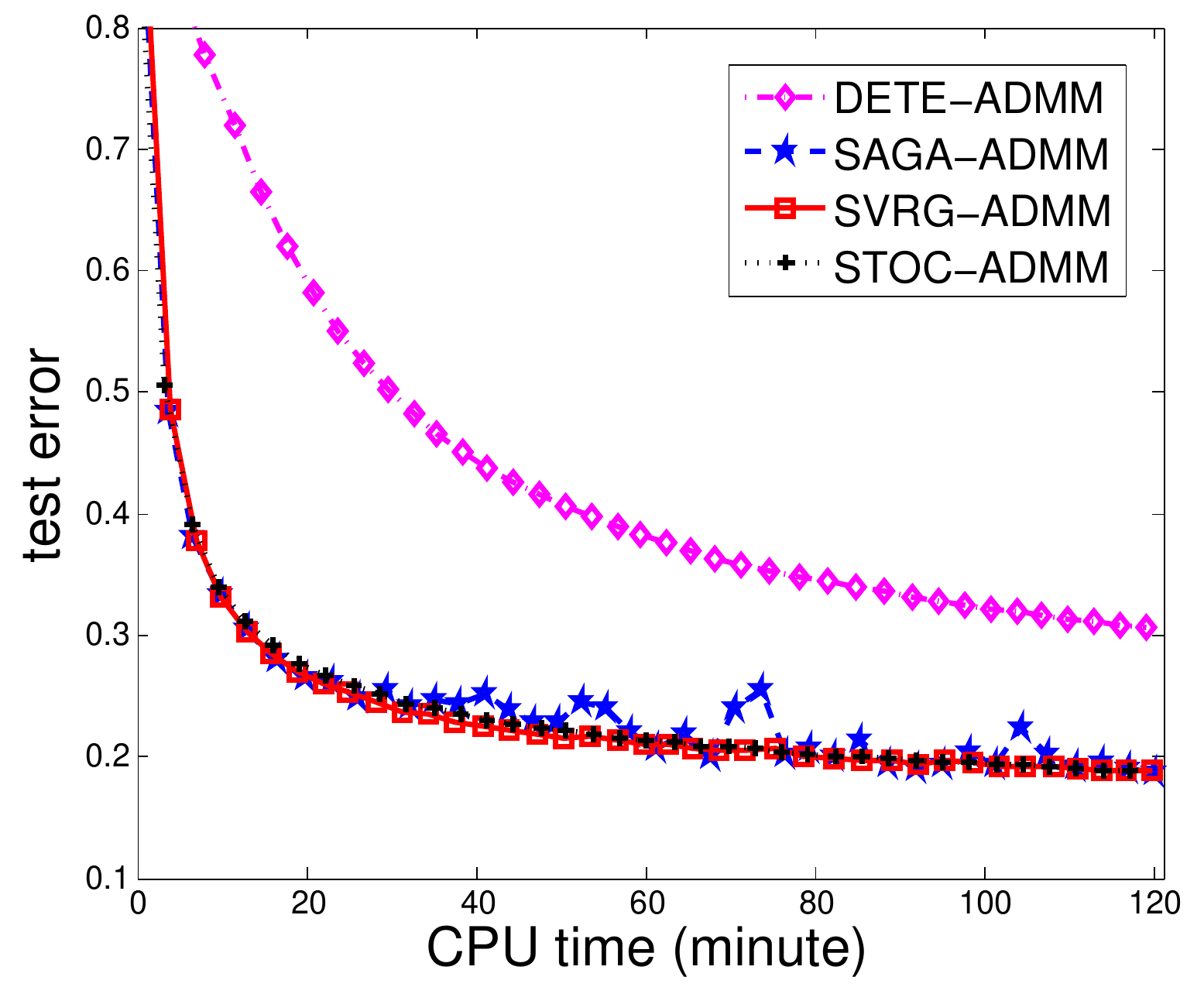}}
\caption{Test loss \emph{versus} CPU time of the \emph{nonconvex} \textbf{multi-task learning} on some real datasets.}
\label{fig:8}
\vspace{-1em}
\end{figure*}

In the experiment, we use some publicly available datasets\footnote{
\emph{letter}, \emph{sensorless}, \emph{covtype}, \emph{mnist} and \emph{mnist8m}
are from the LIBSVM website (www.csie.ntu.edu.tw/~cjlin/libsvmtools/datasets/).},
which are summarized in Table \ref{tab:4}. We use the mini-batch size of $M=100$ on \emph{letter},
$M=300$ on \emph{sensorless} and \emph{mnist}, $M=500$ on \emph{covtype},
and $M=1000$ on \emph{mnist8m}.
In the algorithms, we use the same initial solution $x_0$ from the standard normal distribution
and choose a fixed step size $\eta = 0.8$.
In the problem \eqref{eq:81}, we fix the regularization parameters $\nu_1 = 10^{-5}$ and $\nu_2 = 10^{-4}$.

Figs. \ref{fig:7} and \ref{fig:8} show that both objective values and test loss of the
stochastic ADMMs faster decrease than those of the deterministic ADMM,
as CPU time consumed increases.
In particular, though the nonconvex \emph{STOC-ADMM} uses a fixed step size $\eta$,
it shows good performance in the nonconvex multi-task learning with spare and low-rank regularization functions,
and is comparable with both the nonconvex \emph{SVRG-ADMM} and \emph{SAGA-ADMM}.
Due to large training samples, the stochastic gradient of
\emph{SAGA-ADMM} includes many old gradients, and slowly updates.

\subsection{ Varying $\rho$ }

In the subsection, we demonstrate the specific parameter selection for step size $\eta$ of stochastic gradient and
penalty parameter $\rho$ of augmented Lagrangian function.
Specifically, we give a fixed $\eta$, then find an optimal $\rho$.
In the experiment, we use the above simulated data
imposed the overlapping group lasso regularization function, and
set $n=40,000$, $d=400$. In the problem \eqref{eq:80},
we fix the regularization parameter $\nu=10^{-5}$.
In the algorithms, we fix the step size $\eta = 1$.

\begin{figure}[htbp]
\centering
\subfigure[Objective value]{\includegraphics[width=0.24\textwidth]{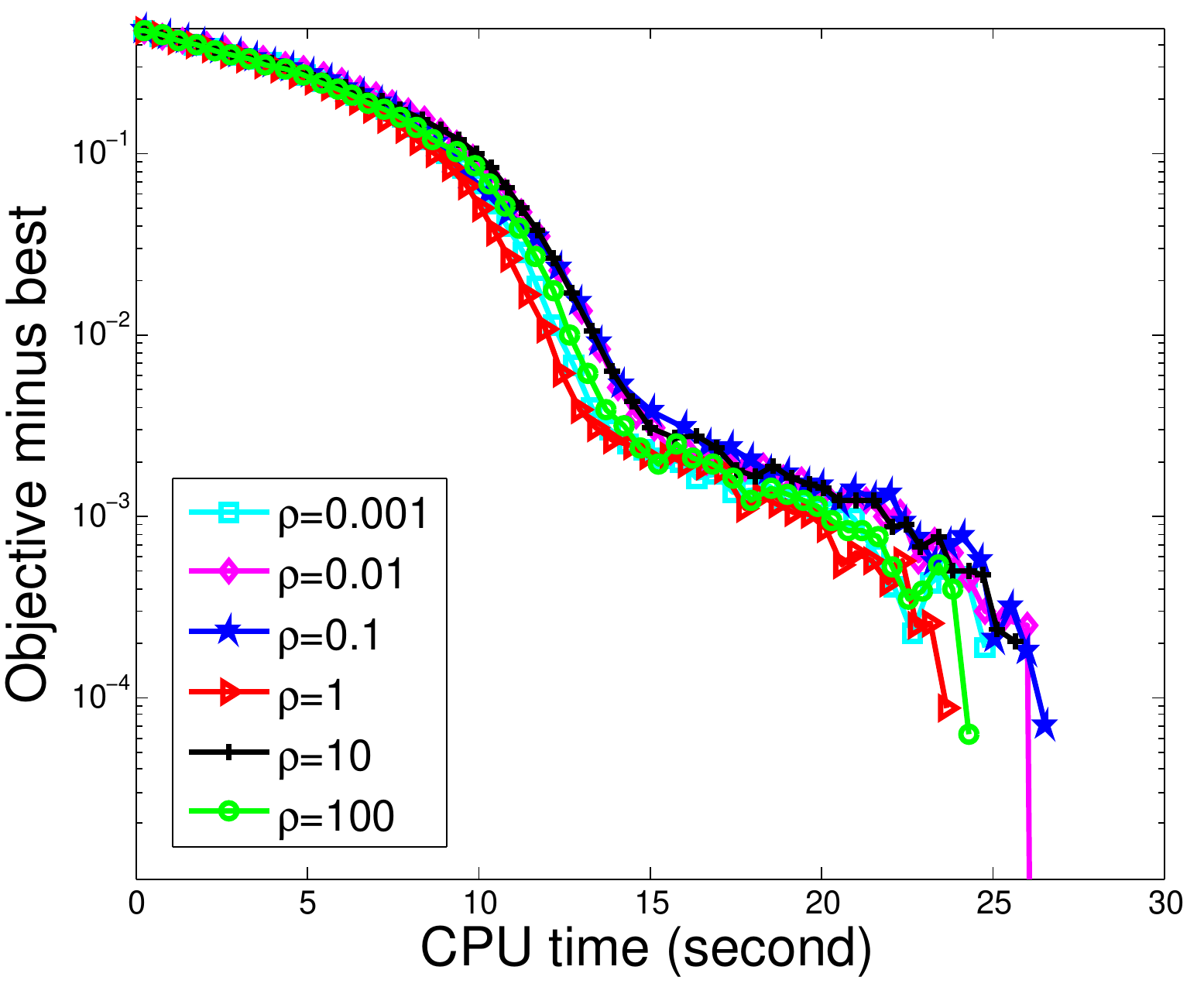}}
\subfigure[Test error]{\includegraphics[width=0.24\textwidth]{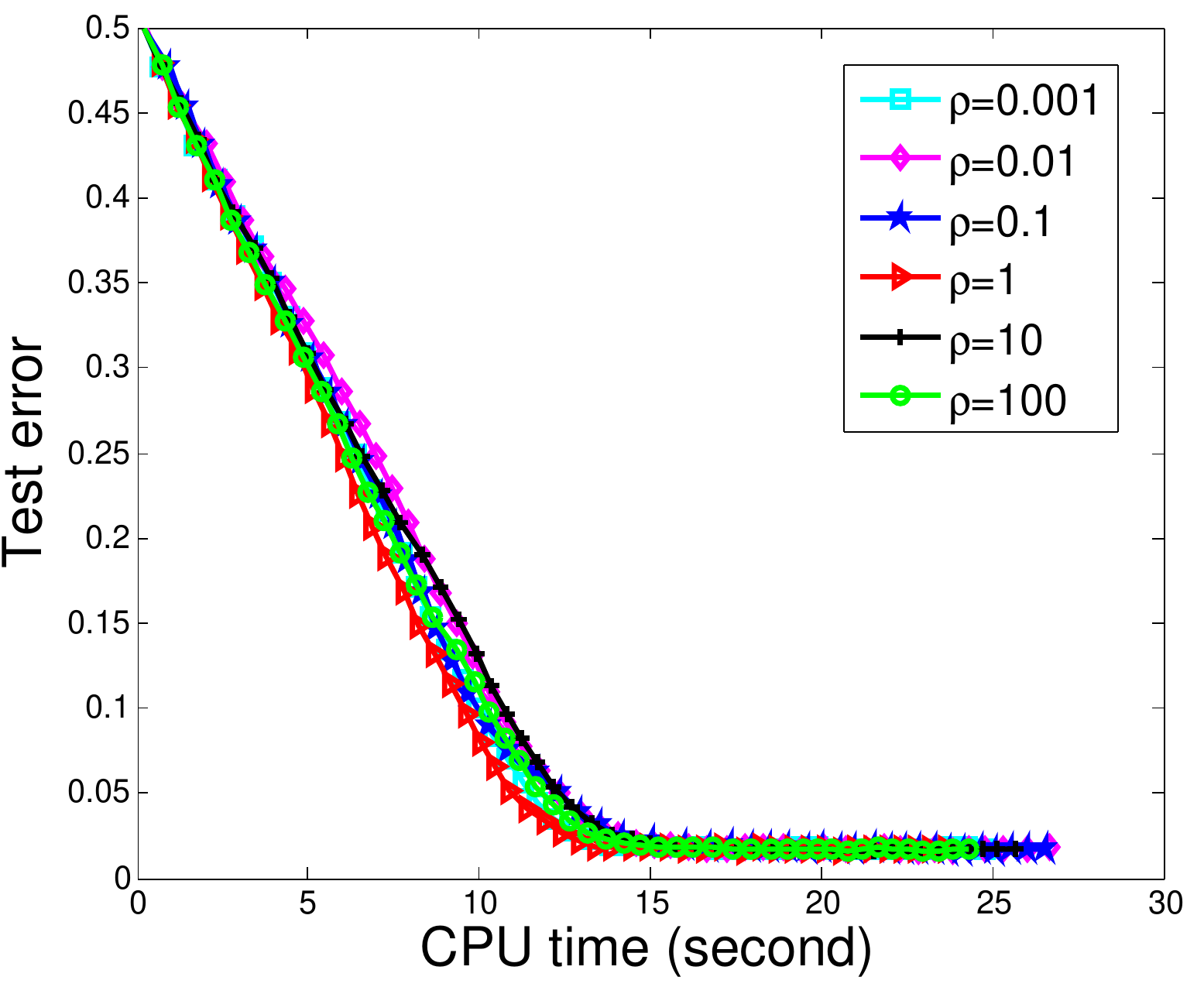}}
\caption{Performance of nonconvex STOC-ADMM at different $\rho$.}
\label{fig:9}
\vspace{-1em}
\end{figure}

\begin{figure}[htbp]
\centering
\subfigure[Objective value]{\includegraphics[width=0.24\textwidth]{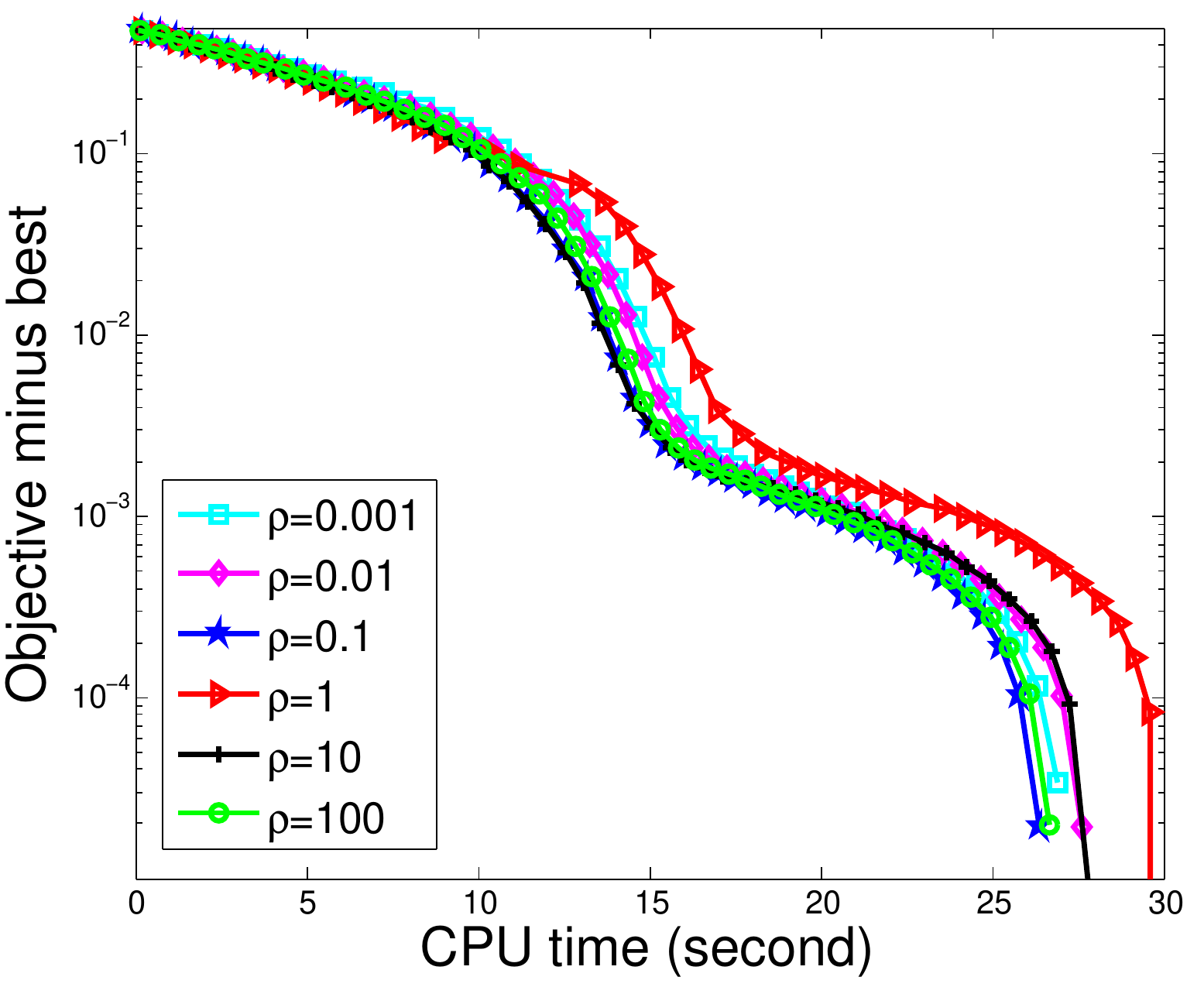}}
\subfigure[Test error]{\includegraphics[width=0.24\textwidth]{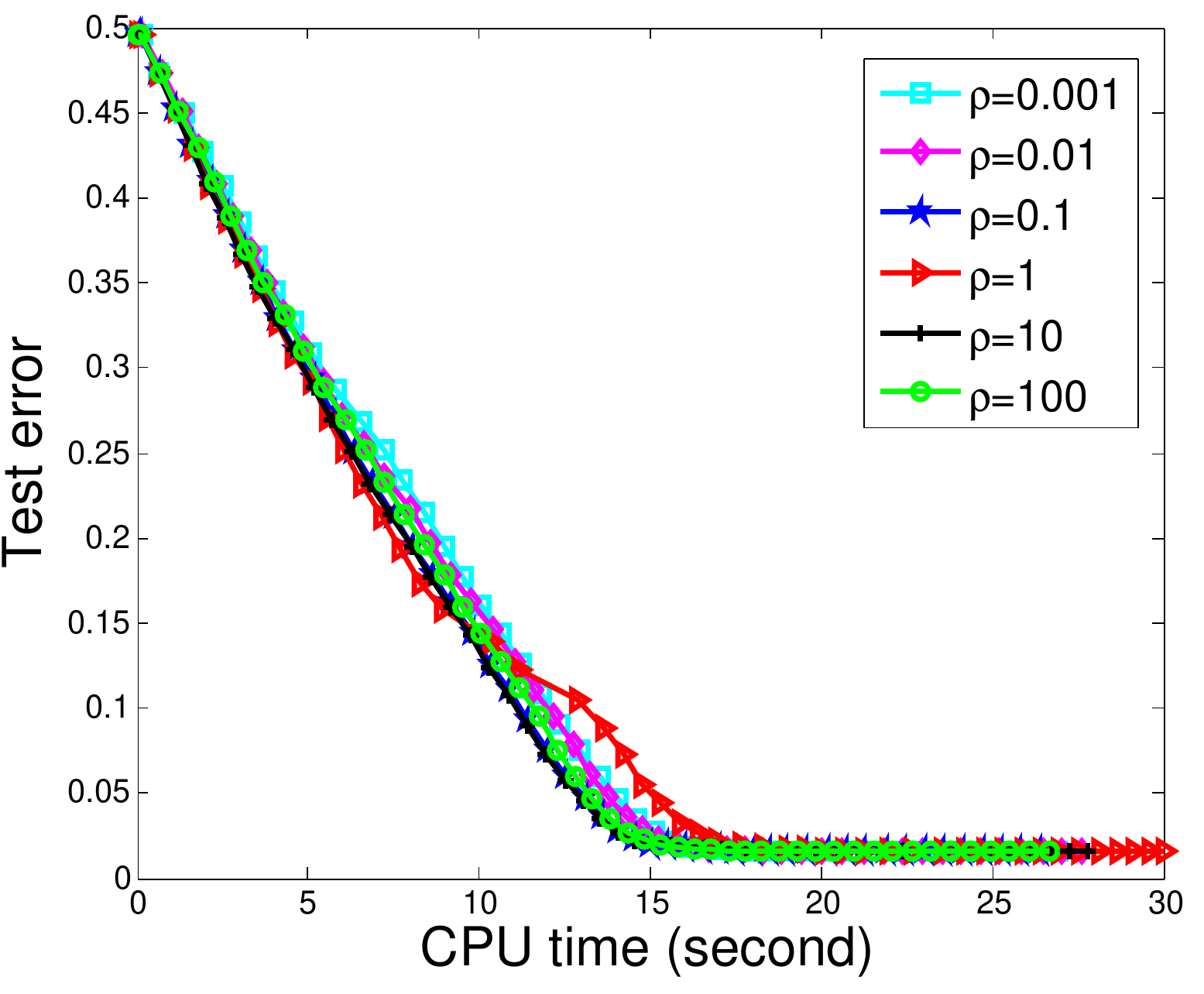}}
\caption{Performance of nonconvex SVRG-ADMM at different $\rho$.}
\label{fig:10}
\vspace{-1em}
\end{figure}

\begin{figure}[htbp]
\centering
\subfigure[Objective value]{\includegraphics[width=0.24\textwidth]{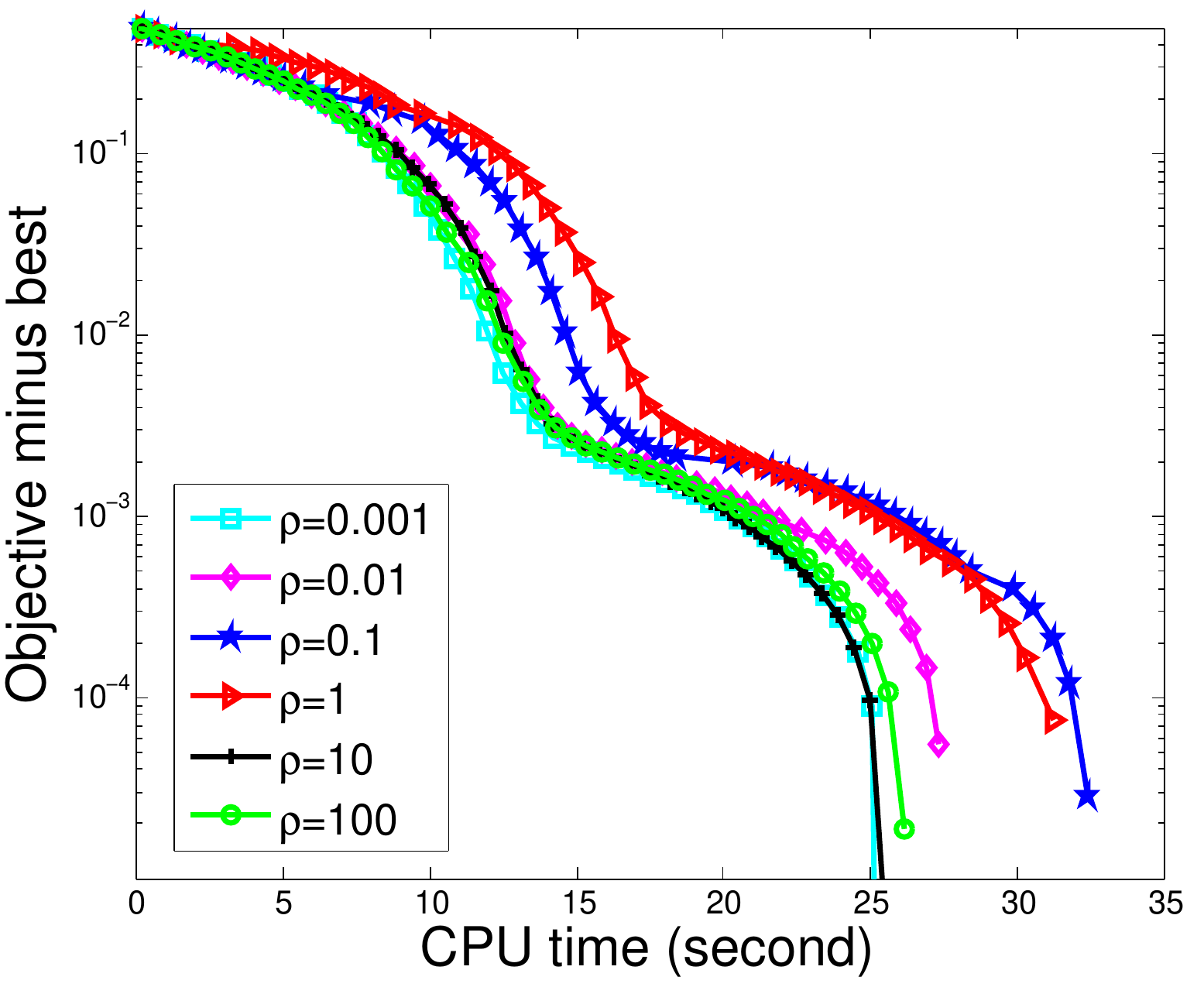}}
\subfigure[Test error]{\includegraphics[width=0.24\textwidth]{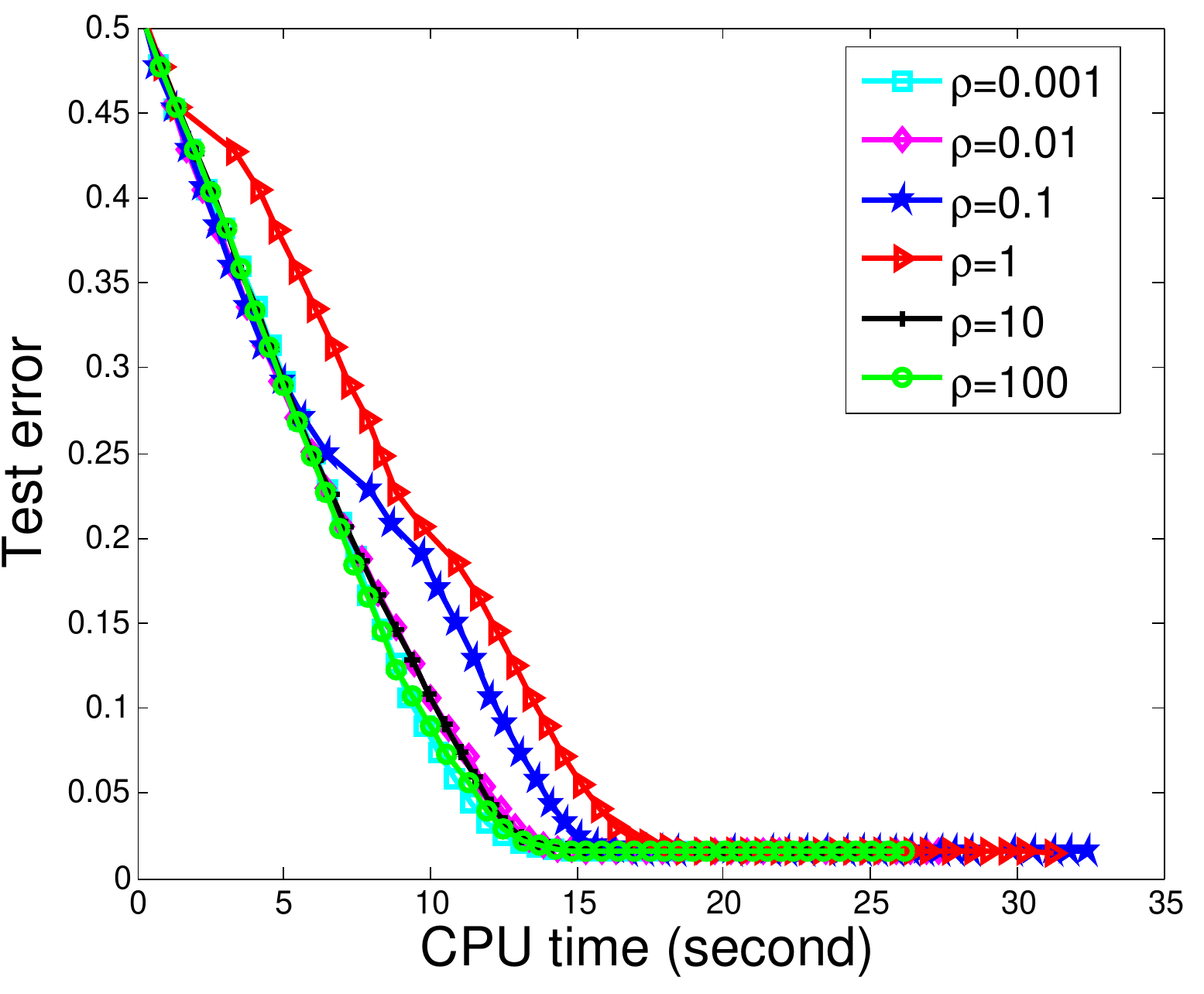}}
\caption{Performance of nonconvex SAGA-ADMM at different $\rho$.}
\label{fig:11}
\vspace{-1em}
\end{figure}

Figs. \ref{fig:9}, \ref{fig:10} and \ref{fig:11} show the objective value and
test error \emph{versus} CPU time with different $\rho$. From these results,
we can find that given an appropriate step size $\eta$, the proposed mini-batch stochastic algorithms have good performances
in a wide-range of parameter $\rho$.
In particular, when the parameter $\rho$ satisfies the above conditions \eqref{eq:le3}, \eqref{eq:le8} and \eqref{eq:le13},
these mini-batch algorithms show good performances.

\section{Conclusion}
In the paper, we have studied the mini-batch stochastic ADMMs
for the nonconvex nonsmooth optimization.
We have theoretically proved that, give mini-batch size $M=O(1/\epsilon)$,
the mini-batch stochastic ADMM without VR (STOC-ADMM) has the convergence rate of $O(1/T)$
to obtain an $\epsilon$-stationary point.
In theoretical analysis, the mini-batch size $M$ may be very large when $\epsilon$ is small.
However, the above extensive experimental results show that
STOC-ADMM still has good performances given a moderate $M$, and is comparable with
both SVRG-ADMM and SAGA-ADMM.
In particular, as long as the step size $\eta$ and the penalization
parameter $\rho$ satisfy the above condition \eqref{eq:le3} instead of $\eta=O(\frac{1}{t})$
used in the convex stochastic ADMM \cite{ouyang2013stochastic},
STOC-ADMM is convergent, and reaches a convergence rate of $O(1/T)$.

Moreover, we have extended the mini-batch stochastic gradient method to
both the non-convex SVRG-ADMM and SAGA-ADMM proposed in our initial manuscript \cite{huang2016stochastic},
and also proved that these mini-batch stochastic ADMMs reach the convergence rate of $O(1/T)$.
Though both SVRG-ADMM and SAGA-ADMM reach the convergence rate of $O(1/T)$ without the
condition on $M$, SVRG-ADMM requires frequently compute gradients over the full data, and
SAGA-ADMM requires memory of the
same size for storing the old gradients.
In the future work, we will develop a more efficient stochastic ADMM algorithm for automatically adapting to the system resources,
and yield the best performance in practice.
In addition, we will propose some accelerated
stochastic ADMMs for nonconvex optimization by using the
momentum techniques.

% if have a single appendix:
%\appendix[Proof of the Zonklar Equations]
% or
%\appendix  % for no appendix heading
% do not use \section anymore after \appendix, only \section*
% is possibly needed

% use appendices with more than one appendix
% then use \section to start each appendix
% you must declare a \section before using any
% \subsection or using \label (\appendices by itself
% starts a section numbered zero.)
%

%\appendices
%\section{Proof of the First Zonklar Equation}
%Appendix one text goes here.
%
%% you can choose not to have a title for an appendix
%% if you want by leaving the argument blank
%\section{}
%Appendix two text goes here.

% use section* for acknowledgment
%\ifCLASSOPTIONcompsoc
%  % The Computer Society usually uses the plural form
%  \section*{Acknowledgments}
%\else
%  % regular IEEE prefers the singular form
%  \section*{Acknowledgment}
%\fi
%
%This paper is supported in part by the key founding of National Natural Science Foundation of China (NSFC) under Grant No.61732006,
%and in part by the NSFC under Grant No.61672281.

% Can use something like this to put references on a page
% by themselves when using endfloat and the captionsoff option.
\ifCLASSOPTIONcaptionsoff
  \newpage
\fi

% trigger a \newpage just before the given reference
% number - used to balance the columns on the last page
% adjust value as needed - may need to be readjusted if
% the document is modified later
%\IEEEtriggeratref{8}
% The "triggered" command can be changed if desired:
%\IEEEtriggercmd{\enlargethispage{-5in}}

% references section

% can use a bibliography generated by BibTeX as a .bbl file
% BibTeX documentation can be easily obtained at:
% http://mirror.ctan.org/biblio/bibtex/contrib/doc/
% The IEEEtran BibTeX style support page is at:
% http://www.michaelshell.org/tex/ieeetran/bibtex/
%\bibliographystyle{IEEEtran}
% argument is your BibTeX string definitions and bibliography database(s)
%\bibliography{IEEEabrv,../bib/paper}
%
% <OR> manually copy in the resultant .bbl file
% set second argument of \begin to the number of references
% (used to reserve space for the reference number labels box)
\bibliographystyle{IEEEtran}
\bibliography{IEEEabrv,reference}

% biography section
%
% If you have an EPS/PDF photo (graphicx package needed) extra braces are
% needed around the contents of the optional argument to biography to prevent
% the LaTeX parser from getting confused when it sees the complicated
% \includegraphics command within an optional argument. (You could create
% your own custom macro containing the \includegraphics command to make things
% simpler here.)
%\begin{IEEEbiography}[{\includegraphics[width=1in,height=1.25in,clip,keepaspectratio]{mshell}}]{Michael Shell}
% or if you just want to reserve a space for a photo:

% You can push biographies down or up by placing
% a \vfill before or after them. The appropriate
% use of \vfill depends on what kind of text is
% on the last page and whether or not the columns
% are being equalized.

%\vfill

% Can be used to pull up biographies so that the bottom of the last one
% is flush with the other column.
%\enlargethispage{-5in}

\begin{onecolumn}

%\appendices

\begin{appendices}

%% Convergence Analysis of Nonconvex Mini-batch stochastic ADMM

\section{Convergence Analysis of Nonconvex Mini-batch STOC-ADMM}

\subsection{Proof of Lemma 1}
\label{app:lem-2}

\begin{proof}
 By the optimal condition of step 6 in Algorithm \ref{alg:1}, we have
 \begin{align}
  0 & = G(x_t,\xi_{\mathcal{I}_t})-A^T\lambda_t+\rho A^T(Ax_{t+1}+By_{t+1}-c) - \frac{H}{\eta}(x_t - x_{t+1}) \nonumber \\
  & = G(x_t,\xi_{\mathcal{I}_t})- A^T\lambda_{t+1} - \frac{H}{\eta}(x_t - x_{t+1}), \nonumber
 \end{align}
 where the second equality is due to step 7 in Algorithm \ref{alg:1}.
 Thus, we have
 \begin{align}
  A^T \lambda_{t+1} = G(x_t,\xi_{\mathcal{I}_t}) - \frac{H}{\eta}(x_t - x_{t+1}). 
 \end{align}
 It follows that
 \begin{align} \label{eq:14}
  \lambda_{t+1} = (A^T)^+ \big( G(x_t,\xi_{\mathcal{I}_t}) + \frac{H}{\eta}(x_{t+1}-x_t) \big),
 \end{align}
where $(A^T)^+$ is the pseudoinverse of $A^T$. By Assumption 5, without loss of generality, we use the full column matrix $A$. Then we have $(A^T)^+=A(A^TA)^{-1}$.
Using the equality \eqref{eq:14}, we have
 \begin{align}
 \|\lambda_{t+1}-\lambda_{t}\|^2 & = \|(A^T)^+ \big( G(x_t,\xi_{\mathcal{I}_t}) + \frac{H}{\eta}(x_{t+1}-x_t)  - G(x_{t-1},\xi_{\mathcal{I}_{t-1}}) - \frac{H}{\eta} (x_t -x_{t-1})\big)\|^2  \nonumber \\
 & \leq (\phi_{\min}^A)^{-1} \|G(x_t,\xi_{\mathcal{I}_t}) - G(x_{t-1},\xi_{\mathcal{I}_{t-1}}) + \frac{H}{\eta} (x_{t+1}-x_t) - \frac{H }{\eta}(x_t-x_{t-1})\|^2 \nonumber \\
 & = (\phi_{\min}^A)^{-1} \|G(x_t,\xi_{\mathcal{I}_t})-\nabla f(x_{t}) + \nabla f(x_{t}) - \nabla f(x_{t-1}) + \nabla f(x_{t-1})
    - G(x_{t-1},\xi_{\mathcal{I}_{t-1}}) \nonumber \\
 & \quad + \frac{H}{\eta} (x_{t+1}-x_t) - \frac{H }{\eta}(x_t-x_{t-1})\|^2 \nonumber \\
 & \mathop{\leq}^{(i)} \frac{5}{\phi_{\min}^A} \|G(x_t,\xi_{\mathcal{I}_t})-\nabla f(x_{t})\|^2
 + \frac{5}{\phi_{\min}^A} \|G(x_{t-1},\xi_{\mathcal{I}_{t-1}})-\nabla f(x_{t-1})\|^2
 + \frac{5(\phi^H_{\max})^2}{\phi_{\min}^A \eta^2} \|x_t-x_{t+1}\|^2 \nonumber \\
 & \quad + \frac{5((\phi^H_{\max})^2+\eta^2L^2)}{\phi_{\min}^A\eta^2} \|x_{t-1}-x_{t}\|^2,    \label{eq:15}
 \end{align}
 where the inequality $(i)$ holds by Assumption 2.

 Taking expectation conditioned on information $\mathcal{I}_t$ to \eqref{eq:15}, we have
 \begin{align}
\mathbb{E}\|\lambda_{t+1}-\lambda_{t}\|^2 &\leq \frac{5}{\phi_{\min}^A} \mathbb{E} \|G(x_t,\xi_{\mathcal{I}_t})-\nabla f(x_{t})\|^2
  + \frac{5}{\phi_{\min}^A} \mathbb{E} \|G(x_{t-1},\xi_{\mathcal{I}_{t-1}})-\nabla f(x_{t-1})\|^2 \nonumber \\
 & \quad + \frac{ 5\eta^2(\phi^H_{\max})^2 }{\phi_{\min}^A} \|x_t-x_{t+1}\|^2 + \frac{ 5(\eta^2(\phi^H_{\max})^2+L^2) }{\phi_{\min}^A} \|x_{t-1}-x_{t}\|^2 \nonumber \\
 & \mathop{\leq}^{(i)} \frac{10\sigma^2}{M \phi_{\min}^A}
 + \frac{5(\phi^H_{\max})^2}{\phi_{\min}^A\eta^2}\mathbb{E}\|x_{t+1}-x_t\|^2  + \frac{5((\phi^H_{\max})^2 + \eta^2L^2)}{\phi_{\min}^A\eta^2}\|x_{t}-x_{t-1}\|^2, \nonumber \\
 & = \zeta \|x_{t}-x_{t-1}\|^2 + \zeta_1 \mathbb{E}\|x_{t+1}-x_t\|^2 + \frac{10\sigma^2}{M \phi_{\min}^A} \nonumber
 \end{align}
 where the inequality $(i)$ holds by Assumption 1.
\end{proof}

\subsection{Proof of Lemma 2}
\label{app:lem-3}

\begin{proof}
By the step 5 of Algorithm \ref{alg:1}, we have
\begin{align}
 \mathcal {L}_\rho (x_t, y_{t+1},\lambda_t) \leq \mathcal {L}_\rho (x_t, y_{t},\lambda_t). \label{eq:18}
\end{align}
By the optimal condition of step 6 in Algorithm \ref{alg:1}, we have
\begin{align}
0 & =(x_t-x_{t+1})^T\big[G(x_t,\xi_{\mathcal{I}_t}) -A^T\lambda_t- \frac{H}{\eta}(x_t-x_{t+1}) + \rho A^T(Ax_{t+1}+By_{t+1}-c) \big] \nonumber \\
 & = (x_t-x_{t+1})^T\big[G(x_t,\xi_{\mathcal{I}_t}) - \nabla f(x_t) + \nabla f(x_t) -A^T\lambda_t - \frac{H}{\eta}(x_t-x_{t+1})
 +\rho A^T(Ax_{t+1}+By_{t+1}-c)\big] \nonumber \\
 & \mathop{\leq}^{(i)} f(x_t) - f(x_{t+1}) + \frac{L}{2}\|x_{t+1}-x_t\|^2 + (x_t -x_{t+1})^T\big(G(x_t,\xi_{\mathcal{I}_t}) - \nabla f(x_t)\big)  -\frac{1}{\eta}\|x_{t+1}-x_t\|^2_H  \nonumber \\
 &\quad - \lambda_t^T(A x_t- Ax_{t+1}) + \rho (Ax_t -Ax_{t+1})^T(Ax_{t+1}+By_{t+1}-c) \nonumber \\
 & \mathop{=}^{(ii)} f(x_t) - f(x_{t+1}) + \frac{L}{2}\|x_{t+1}-x_t\|^2 + (x_t -x_{t+1})^T\big(G(x_t,\xi_{\mathcal{I}_t}) - \nabla f(x_t)\big) -\frac{1}{\eta}\|x_{t+1}-x_t\|^2_H \nonumber \\
 & \quad - \lambda_t^T(A x_t+By_{t+1}-c) + \lambda_t^T(Ax_{t+1}+By_{t+1}-c)
  + \frac{\rho}{2}\|Ax_{t}+By_{t+1}-c\|^2 \nonumber \\
 & \quad - \frac{\rho}{2}\|Ax_{t+1}+By_{t+1}-c\|^2 -\frac{\rho}{2}\|Ax_t-Ax_{t+1}\|^2 \nonumber \\
 & = \mathcal {L}_\rho (x_t, y_{t+1},\lambda_t)- \mathcal {L}_\rho (x_{t+1}, y_{t+1},\lambda_t)
 + (x_t -x_{t+1})^T\big(G(x_t,\xi_{\mathcal{I}_t}) - \nabla f(x_t)\big) + \frac{L}{2}\|x_{t+1}-x_t\|^2\nonumber \\
 & \quad -\frac{1}{\eta}\|x_{t+1}-x_t\|^2_H -\frac{\rho}{2}\|Ax_t-Ax_{t+1}\|^2 \nonumber \\
 & \mathop{\leq}^{(iii)} \mathcal {L}_\rho (x_t, y_{t+1},\lambda_t)- \mathcal {L}_\rho (x_{t+1}, y_{t+1},\lambda_t)
 + \frac{1}{2}\|G(x_t,\xi_{\mathcal{I}_t}) - \nabla f(x_t)\|^2 \nonumber \\
 & \quad - (\frac{\phi_{\min}^H}{\eta} + \frac{\rho\phi_{\min}^A}{2}-\frac{L+1}{2})\|x_t-x_{t+1}\|^2, \label{eq:19}
\end{align}
where the inequality (i) holds by \eqref{eq:13}; the equality (ii) holds by using the equality
$(a-b)^T(b-c) = \frac{1}{2}(\|a-c\|^2-\|a-b\|^2-\|b-c\|^2)$ on the term $\rho (Ax_t -Ax_{t+1})^T(Ax_{t+1}+By_{t+1}-c)$;
the inequality (iii) holds by the Cauchy inequality.
Taking expectation conditioned on information $\mathcal{I}_t$ to \eqref{eq:19}, we have
\begin{align}
\mathbb{E} [\mathcal {L}_\rho (x_{t+1}, y_{t+1},\lambda_t)] \leq \mathcal {L}_\rho (x_t, y_{t+1},\lambda_t) - (\frac{\phi_{\min}^H}{\eta} + \frac{\rho\phi_{\min}^A}{2}-\frac{L+1}{2})\mathbb{E}\|x_{t+1}-x_t\|^2 + \frac{\sigma^2}{2M}. \label{eq:20}
\end{align}
By the step 7 of Algorithm \ref{alg:1}, we have
\begin{align}
\mathbb{E} [\mathcal {L}_\rho (x_{t+1}, y_{t+1},\lambda_{t+1})-\mathcal {L}_\rho (x_{t+1}, y_{t+1},\lambda_t)]
&= \frac{1}{\rho}\mathbb{E} \|\lambda_t-\lambda_{t+1}\|^2 \nonumber \\
& \leq \frac{5(L^2\eta^2 + (\phi_{\max}^H)^2)}{\phi_{\min}^A\eta^2\rho}\|x_{t}-x_{t-1}\|^2
   + \frac{ 5(\phi_{\max}^H)^2}{\phi_{\min}^A\eta^2\rho} \mathbb{E}\|x_{t+1}-x_t\|^2 \nonumber \\
& \quad + \frac{10\sigma^2}{M\phi_{\min}^A\rho}, \label{eq:21}
\end{align}
where the inequality $(i)$ holds by the Lemma \ref{lem:2}.
Combining \eqref{eq:18}, \eqref{eq:20} and \eqref{eq:21}, we have
\begin{align}
 \mathbb{E} [\mathcal {L}_\rho (x_{t+1}, y_{t+1},\lambda_{t+1})] \leq  & \mathcal {L}_\rho (x_t, y_{t},\lambda_t) + \frac{5(L^2\eta^2 + (\phi_{\max}^H)^2)}{\phi_{\min}^A\eta^2\rho}\|x_{t}-x_{t-1}\|^2 \nonumber \\
 & - (\frac{\phi_{\min}^H}{\eta} + \frac{\rho\phi_{\min}^A}{2}-\frac{L+1}{2}-\frac{ 5(\phi_{\max}^H)^2}{\phi_{\min}^A\eta^2\rho})\mathbb{E}\|x_t-x_{t+1}\|^2 + \frac{(\phi_{\min}^A\rho+20)\sigma^2}{2\phi_{\min}^A\rho M}. \nonumber
\end{align}

Next, we define a useful sequence $\big\{\Psi_t\big\}_{t=1}^T$ as follows:
\begin{align} 
\Psi_t =  \mathbb{E}\big[ \mathcal{L}_{\rho}(x_t,y_t,\lambda_t) + \frac{5(L^2\eta^2 + (\phi_{\max}^H)^2)}{\rho \phi_{\min}^A\eta^2}\|x_t-x_{t-1}\|^2 \big].
\end{align}
Then we have
\begin{align} \label{eq:A45}
\Psi_{t+1} - \Psi_t &\leq- \underbrace{\big(\frac{\phi_{\min}^H}{\eta} + \frac{\rho\phi_{\min}^A}{2}-\frac{L+1}{2}
 -\frac{5(L^2\eta^2 + 2(\phi_{\max}^H)^2)}{\phi_{\min}^A\eta^2\rho}\big)}_{\gamma} \mathbb{E}\|x_t-x_{t+1}\|^2 + \frac{(\phi_{\min}^A\rho+20)\sigma^2}{2\phi_{\min}^A\rho M}, \nonumber \\
& = -\gamma \mathbb{E}\|x_t-x_{t+1}\|^2 + \frac{(\phi_{\min}^A\rho+20)\sigma^2}{2\phi_{\min}^A\rho M}. 
\end{align}
Finally, using \eqref{eq:le3} and the properties of \emph{quadratic equation in one unknown}, we have $\gamma>0$.

Since $A$ is a full column rank matrix,
we have $(A^T)^+ = A(A^T A)^{-1}$.
It follows that $\sigma_{\max}\big(((A^T)^+)^T(A^T)^+\big) = \sigma_{\max}((A^TA)^{-1}) = \frac{1}{\phi_{\min}^A}$.
using \eqref{eq:14}, we have
\begin{align}
\mathcal{L}_{\rho} (x_{t+1},y_{t+1},\lambda_{t+1})
& = f(x_{t+1}) + g(y_{t+1}) - \lambda_{t+1}^T(Ax_{t+1} + By_{t+1} -c) + \frac{\rho}{2}\|Ax_{t+1} + By_{t+1} -c\|^2 \nonumber \\
& = f(x_{t+1}) + g(y_{t+1}) - \langle(A^T)^+ \big( G(x_t,\xi_{\mathcal{I}_t}) + \frac{H}{\eta}(x_{t+1}-x_t) \big), Ax_{t+1} + By_{t+1} -c\rangle
 + \frac{\rho}{2}\|Ax_{t+1} + By_{t+1} -c\|^2 \nonumber \\
& = f(x_{t+1}) + g(y_{t+1}) - \langle(A^T)^+ \big( G(x_t,\xi_{\mathcal{I}_t}) - \nabla f(x_{t}) + \nabla f(x_{t})+ \frac{H}{\eta}(x_{t+1}-x_t) \big), Ax_{t+1} + By_{t+1} -c\rangle  \nonumber \\
& \quad +  \frac{\rho}{2}\|Ax_{t+1} + By_{t+1} -c\|^2 \nonumber \\
& \geq f(x_{t+1}) + g(y_{t+1}) - \frac{2}{\phi^A_{\min}\rho}\|G(x_t,\xi_{\mathcal{I}_t}) - \nabla f(x_{t})\|^2 - \frac{2}{\phi^A_{\min}\rho}\|\nabla f(x_{t})\|^2
- \frac{2(\phi_{\max}^H)^2}{\phi^A_{\min}\eta^2\rho}\|x_{t+1}-x_t\|^2 \nonumber \\
& \quad + \frac{\rho}{8}\|Ax_{t+1} + By_{t+1} -c\|^2 \nonumber\\
& \geq f(x_{t+1}) + g(y_{t+1}) -\frac{2\sigma^2}{M\phi^A_{\min}\rho} - \frac{2\delta^2}{\phi^A_{\min}\rho} \!-\! \frac{2(\phi_{\max}^H)^2}{\phi^A_{\min}\eta^2\rho}\|x_{t+1}-x_t\|^2 \nonumber \\
& \geq f^* + g^* -\frac{2\sigma^2}{M\phi^A_{\min}\rho} - \frac{2\delta^2}{\phi^A_{\min}\rho} \!-\! \frac{2(\phi_{\max}^H)^2}{\phi^A_{\min}\eta^2\rho}\|x_{t+1}-x_t\|^2,
\end{align}
where the first inequality is obtained by applying $ \langle a, b\rangle \leq \frac{1}{2\beta}\|a\|^2 + \frac{\beta}{2}\|b\|^2$ to the terms
$\langle(A^T)^+(G(x_t,\xi_{\mathcal{I}_t}) - \nabla f(x_{t})), Ax_{t+1} + By_{t+1} -c\rangle$, $\langle(A^T)^+G(x_t,\xi_{\mathcal{I}_t}), Ax_{t+1} + By_{t+1} -c\rangle $ and
$\langle(A^T)^+\frac{H}{\eta}(x_{t+1}-x_t), Ax_{t+1} + By_{t+1} -c\rangle$ with $\beta = \frac{\rho}{4}$, respectively,
and the second inequality follows by the inequality \label{eq:11} and Assumption 3, and the third inequality holds by Assumption 4.
Using the definition of $\Psi_t$, we have
\begin{align}
 \Psi_{t+1}\geq f^* + g^* -\frac{2\sigma^2}{M\phi^A_{\min}\rho} - \frac{2\delta^2}{\phi^A_{\min}\rho}, \ \forall \ t=0,1,2,\cdots.
\end{align}
It follows that the function $\Psi_t$ is bounded from below. Let $\Psi^*$ denotes a lower bound of sequence $\{\Psi_t\}_{t=1}^T$.

Telescoping inequality \eqref{eq:A45} over $t$ from $0$ to $T$,
we have
\begin{align} \label{eq:A48}
\frac{1}{T} \sum_{t=0}^{T-1} \mathbb{E}\|x_t-x_{t+1}\|^2 \leq \frac{\Psi_0-\Psi^*}{\gamma T} + \frac{(\phi_{\min}^A\rho+20)\sigma^2}{2\gamma\phi_{\min}^A\rho M}.
\end{align}

\end{proof}

\subsection{Proof of Theorem 1}
\label{app:th-1}

\begin{proof}
First, we define a useful variable $\theta_t = \mathbb{E}\big[\|x_t-x_{t+1}\|^2+\|x_{t-1}-x_t\|^2\big]$.
By \eqref{eq:A48}, then we have
\begin{align}
t^*= \mathop{\arg\min}_{2 \leq t \leq T+1} \mathbb{E}[\theta_{t}] \leq \frac{2}{\gamma T}
(\Psi_1  - \Psi^*) + \frac{(\phi_{\min}^A\rho+20)\sigma^2}{2\gamma\phi_{\min}^A\rho M}. \nonumber
\end{align}

By \eqref{eq:14}, we have
\begin{align}
 \mathbb{E}\|A^T\lambda_{t+1}-\nabla f(x_{t+1})\|^2  
 & = \mathbb{E}\|G(x_{t},\xi_{\mathcal{I}_t}) - \nabla f(x_{t+1}) - \frac{H}{\eta} (x_t-x_{t+1})\|^2 \nonumber \\
 & = \mathbb{E}\|G(x_{t},\xi_{\mathcal{I}_t})-\nabla f(x_{t}) +\nabla f(x_{t})- \nabla f(x_{t+1})
  - \frac{H}{\eta}  (x_t-x_{t+1})\|^2  \nonumber \\
 & \leq  3\big( L^2+\frac{(\phi^H_{\max})^2}{\eta^2} \big)\|x_t-x_{t+1}\|^2 + \frac{3\sigma^2}{M} \nonumber \\
 & \leq 3\big( L^2+\frac{(\phi^H_{\max})^2}{\eta^2} \big)\theta_{t} + \frac{3\sigma^2}{M}.   \label{eq:26}
\end{align}
By the step 7 of Algorithm \ref{alg:1}, we have
\begin{align}
\mathbb{E}\|Ax_{t+1}+By_{t+1}-c\|^2 &= \frac{1}{\rho^2} \mathbb{E} \|\lambda_{t+1}-\lambda_t\|^2  \nonumber \\
 & \leq \frac{5(L^2\eta^2 + (\phi_{\max}^H)^2)}{\phi_{\min}^A\eta^2\rho^2}\|x_{t}-x_{t-1}\|^2
   + \frac{ 5(\phi_{\max}^H)^2}{\phi_{\min}^A\eta^2\rho^2} \mathbb{E}\|x_{t+1}-x_t\|^2 + \frac{10\sigma^2}{M\phi_{\min}^A\rho^2}, \nonumber \\
 & \leq \frac{5(L^2\eta^2+(\phi^H_{\max})^2)}{\phi^A_{\min}\rho^2\eta^2} \theta_{t} + \frac{10\sigma^2}{\phi_{\min}^A\rho^2 M} \nonumber \\
 & = \frac{\zeta}{\rho^2}\theta_{t} + \frac{10\sigma^2}{\phi_{\min}^A\rho^2 M} .  \label{eq:27}
\end{align}
By the step 5 of Algorithm \ref{alg:1},
there exists a subgradient $\mu \in \partial g(y_{t+1})$ such that
\begin{align}
 \mathbb{E}\big[\mbox{dist} (B^T\lambda_{t+1}, \partial g(y_{t+1}))^2\big]  & \leq \|\mu-B^T\lambda_{t+1}\|^2 \nonumber \\
 & = \|B^T\lambda_t-\rho B^T(Ax_t+By_{t+1}-c)-B^T\lambda_{t+1}\|^2 \nonumber \\
 & = \|\rho B^TA(x_{t+1}-x_{t})\|^2 \nonumber \\
 & \leq \rho^2\|B\|^2\|A\|^2\|x_{t+1}-x_t\|^2 \nonumber \\
 & \leq \rho^2\|B\|^2\|A\|^2 \theta_{t}.     \label{eq:28}
\end{align}
Finally, using the above bounds \eqref{eq:26}, \eqref{eq:27} and \eqref{eq:28}, and the definition \ref{def:1},
an $\epsilon$-stationary point of the problem \eqref{eq:1} holds in expectation.

\end{proof}

%% Convergence Analysis of Nonconvex Mini-batch SVRG-ADMM

\section{Convergence Analysis of Nonconvex Mini-batch SVRG-ADMM}

\subsection{Proof Lemma 4}
\label{app:lem-6}

\begin{proof}
 Since $\hat{\nabla}f(x^{s+1}_t)= \frac{1}{M} \sum_{i_t\in \mathcal{I}_t} \big( \nabla f_{i_t}(x^{s+1}_t)-\nabla f_{i_t}(\tilde{x}^{s}) \big)
 +\nabla f(\tilde{x}^{s})$, we have
 \begin{align}
  &\mathbb{E}\|\hat{\nabla}f(x^{s+1}_t)-\nabla f(x^{s+1}_t)\|^2 \nonumber \\
  & = \mathbb{E}\|\frac{1}{M} \sum_{i_t\in \mathcal{I}_t} \big(\nabla f_{i_t}(x^{s+1}_t)-\nabla f_{i_t}(\tilde{x}^{s})\big)+\nabla f(\tilde{x}^{s})-\nabla f(x_t^{s+1})\|^2 \nonumber \\
  & \mathop{=}^{(i)}  \mathbb{E}\|\frac{1}{M} \sum_{i_t\in \mathcal{I}_t} \big(\nabla f_{i_t}(x^{s+1}_t)-\nabla f_{i_t}(\tilde{x}^{s})\big)\|^2 -
    \|\nabla f(x_t^{s+1})-\nabla f(\tilde{x}^{s})\|^2 \nonumber \\
  & \leq \frac{1}{M^2} \sum_{i_t\in \mathcal{I}_t} \mathbb{E}\|\nabla f_{i_t}(x^{s+1}_t)-\nabla f_{i_t}(\tilde{x}^{s})\|^2 \nonumber \\
  & = \frac{1}{M^2} \sum_{i_t\in \mathcal{I}_t} \frac{1}{n}\sum_{i=1}^n \|\nabla f_{i_t}(x^{s+1}_t)-\nabla f_{i_t}(\tilde{x}^{s})\|^2 \nonumber \\
  & \mathop{\leq}^{(ii)} \frac{L^2}{M} \|x_t^{s+1}-\tilde{x}^s\|^2. \nonumber
 \end{align}
 where the equality (i) holds by the equality $\mathbb{E}(\xi-\mathbb{E}\xi)^2 = \mathbb{E}\xi^2-(\mathbb{E}\xi)^2$ for random variable $\xi$;
 the inequality (ii) holds by \eqref{eq:12}.
\end{proof}

\subsection{Proof of Lemma 5}
\label{app:lem-7}

\begin{proof}
 For simplicity, let $x^{s+1}_t = x_t$, $y^{s+1}_t = y_t$, $\lambda^{s+1}_t = \lambda_t$, and $\tilde{x}=\tilde{x}^s$.
 By the optimal condition of step 10 in Algorithm \ref{alg:2}, we have
 \begin{align}
  0 & = \hat{\nabla}f(x_t) - A^T\lambda_t + \rho A^T(Ax_{t+1}+By_{t+1}-c) - \frac{H}{\eta} (x_t - x_{t+1}) \nonumber \\
    & = \hat{\nabla}f(x_t) - A^T\lambda_{t+1}- \frac{H}{\eta} (x_t - x_{t+1}), \nonumber
 \end{align}
 where the second equality is due to step 11 in Algorithm \ref{alg:2}.
Then we have
 \begin{align}
  A^T\lambda_{t+1} = \hat{\nabla}f(x_t)- \frac{H}{\eta} (x_t - x_{t+1}). \nonumber
 \end{align}
 It follows that
 \begin{align} \label{eq:30}
  \lambda_{t+1} = (A^T)^+ \big( \hat{\nabla} f(x_t) + \frac{H}{\eta}(x_{t+1}-x_t) \big),
 \end{align}
where $(A^T)^+$ is the pseudoinverse of $A^T$. By Assumption 5, without loss of generality, we use the full column matrix $A$. So we have $(A^T)^+=A(A^TA)^{-1}$.
By \eqref{eq:30}, we have
 \begin{align}
\|\lambda_{t+1}-\lambda_{t}\|^2 & = \|(A^T)^+ \big( \hat{\nabla} f(x_t) + \frac{H}{\eta}(x_{t+1}-x_t) - \hat{\nabla} f(x_{t-1}) - \frac{H}{\eta}(x_{t}-x_{t-1}) \big)\|^2   \nonumber \\
 & \leq (\phi^A_{\min})^{-1}\|\hat{\nabla}f(x_{t})-\hat{\nabla}f(x_{t-1}) + \frac{H}{\eta}(x_{t+1}-x_t) - \frac{H}{\eta}(x_{t}-x_{t-1})\|^2 \nonumber \\
 & = (\phi^A_{\min})^{-1}\|\hat{\nabla}f(x_{t})-\nabla f(x_{t}) + \nabla f(x_{t}) - \nabla f(x_{t-1}) + \nabla f(x_{t-1}) - \hat{\nabla}f(x_{t-1}) \nonumber \\
 &  \quad + \frac{H}{\eta}(x_{t+1}-x_t) - \frac{H}{\eta}(x_{t}-x_{t-1}) \|^2 \nonumber \\
 & \mathop{\leq}^{(i)} \frac{5}{\phi^A_{\min}} \|\hat{\nabla}f(x_{t})-\nabla f(x_{t})\|^2 + \frac{5}{\phi^A_{\min}} \|\hat{\nabla}f(x_{t-1})-\nabla f(x_{t-1})\|^2
 + \frac{5(\phi^H_{\max})^2 }{\phi^A_{\min}\eta^2} \|x_t-x_{t+1}\|^2 \nonumber \\
 & \quad + \frac{5(L^2\eta^2+(\phi^H_{\max})^2) }{\phi^A_{\min}\eta^2} \|x_{t-1}-x_{t}\|^2, \label{eq:31}
 \end{align}
 where the inequality (i) holds by Assumption 2.

 Taking expectation conditioned on information $\mathcal{I}_t$ to \eqref{eq:31}, we have
 \begin{align}
\mathbb{E}\|\lambda_{t+1}-\lambda_{t}\|^2 &\leq \frac{5}{\phi^A_{\min}} \mathbb{E} \|\hat{\nabla}f(x_{t})-\nabla f(x_{t})\|^2
 + \frac{5}{\phi^A_{\min}} \mathbb{E} \|\hat{\nabla}f(x_{t-1})-\nabla f(x_{t-1})\|^2+ \frac{5\phi_{\max}^2}{\phi^A_{\min}\eta^2} \|x_t-x_{t+1}\|^2 \nonumber \\
 & \quad + \frac{5(L^2\eta^2+\phi_{\max}^2)}{\phi^A_{\min}\eta^2} \|x_{t-1}-x_{t}\|^2 \nonumber \\
 & \mathop{\leq}^{(i)} \frac{5L^2}{\phi^A_{\min}M} \mathbb{E}\|x_{t}-\tilde{x}\|^2 + \frac{5L^2}{\phi^A_{\min}M}\|x_{t-1}-\tilde{x}\|^2
 + \frac{5\phi_{\max}^2}{\phi^A_{\min}\eta^2} \mathbb{E}\|x_t-x_{t+1}\|^2 \nonumber \\
 & \quad + \frac{5(L^2\eta^2+\phi_{\max}^2)}{\phi^A_{\min}\eta^2}\|x_{t-1}-x_{t}\|^2, \nonumber
 \end{align}
 where the inequality $(i)$ holds by Lemma \ref{lem:6}.
\end{proof}

\subsection{Proof of Lemma 6}
\label{app:lem-8}

\begin{proof}
This proof includes two parts: First, we will prove the sequence
$\{(\Phi^{s}_{t})_{t=1}^m\}_{s=1}^S$
monotonically decreases over $t\in \{1,2,\cdots,m\}$ in each epoch $s\in \{1,2,\cdots,S\}$.
Second, we will prove $\Phi^{s}_{m}
\geq \Phi^{s+1}_{1}$ for $s\in \{1,2,\cdots,S\}$.

For simplicity, we omit the label of each epoch in the first part, i.e.,
let $x^{s+1}_t = x_t$, $y^{s+1}_t = y_t$, $\lambda^{s+1}_t = \lambda_t$ and $\tilde{x}^{s}=\tilde{x}$.
By the step 8 of Algorithm \ref{alg:2}, we have
\begin{align}
 \mathcal {L}_\rho (x_t, y_{t+1},\lambda_t) \leq \mathcal {L}_\rho (x_t, y_{t},\lambda_t). \label{eq:35}
\end{align}

By the optimal condition of step 10 in Algorithm \ref{alg:2}, we have
\begin{align}
0 &=(x_t-x_{t+1})^T\big[\hat{\nabla}f(x_t)-A^T\lambda_t+\rho(Ax_{t+1}+By_{t+1}-c) - \frac{H}{\eta}(x_t-x_{t+1})\big] \nonumber \\
 & = (x_t-x_{t+1})^T\big[\hat{\nabla}f(x_t) - \nabla f(x_t) + \nabla f(x_t) -A^T\lambda_t +\rho A^T(Ax_{t+1}+By_{t+1}-c)- \frac{H}{\eta}(x_t-x_{t+1})\big] \nonumber \\
 & \mathop{\leq}^{(i)} f(x_t) - f(x_{t+1})+ \frac{L}{2} \|x_{t+1}-x_t\|^2
  + (x_t -x_{t+1})^T\big(\hat{\nabla}f(x_t) - \nabla f(x_t)\big)
  - \frac{1}{\eta}\|x_{t+1}-x_t\|^2_H  \nonumber \\
 & \quad - \lambda_t^T(A x_t-Ax_{t+1}) + \rho (Ax_t -Ax_{t+1})^T(Ax_{t+1}+By_{t+1}-c) \nonumber \\
 & \mathop{=}^{(ii)} f(x_t) - f(x_{t+1}) +\frac{L}{2} \|x_{t+1}-x_t\|^2 + (x_t -x_{t+1})^T\big(\hat{\nabla}f(x_t) - \nabla f(x_t)\big) -\frac{1}{\eta}\|x_{t+1}-x_t\|^2_H \nonumber \\
 & \quad - \lambda_t^T(A x_t+By_{t+1}-c) + \lambda_t^T(Ax_{t+1}+By_{t+1}-c) + \frac{\rho}{2}\|Ax_{t}+By_{t+1}-c\|^2 \nonumber \\
 & \quad - \frac{\rho}{2}\|Ax_{t+1}+By_{t+1}-c\|^2 -\frac{\rho}{2}\|Ax_t -Ax_{t+1}\|^2 \nonumber \\
 & = \mathcal {L}_\rho (x_t, y_{t+1},\lambda_t)- \mathcal {L}_\rho (x_{t+1}, y_{t+1},\lambda_t)
  + (x_t -x_{t+1})^T\big(\hat{\nabla}f(x_t) - \nabla f(x_t)\big) \nonumber \\
 & \quad +\frac{L}{2} \|x_{t+1}-x_t\|^2  - \frac{1}{\eta}\|x_{t+1}-x_t\|^2_H -\frac{\rho}{2}\|Ax_t -Ax_{t+1}\|^2 \nonumber \\
 & \mathop{\leq}^{(iii)} \mathcal {L}_\rho (x_t, y_{t+1},\lambda_t)- \mathcal {L}_\rho (x_{t+1}, y_{t+1},\lambda_t)
  + \frac{1}{2}\|\hat{\nabla}f(x_t) - \nabla f(x_t)\|^2 \nonumber \\
 & \quad - (\frac{\phi^H_{\min}}{\eta} + \frac{\phi^A_{\min}\rho}{2}- \frac{L+1}{2})\|x_t-x_{t+1}\|^2,   \label{eq:36}
\end{align}
where the inequality (i) holds by \eqref{eq:13}; the equality (ii) holds by applying the equality
$(a-b)^T(b-c) = \frac{1}{2}(\|a-c\|^2-\|a-b\|^2-\|b-c\|^2)$ on the term $\rho (Ax_t -Ax_{t+1})^T(Ax_{t+1}+By_{t+1}-c)$;
the inequality (iii) holds by the Cauchy inequality.
Taking expectation conditioned on information $\mathcal{I}_t$ to \eqref{eq:36}, we have
\begin{align}
\mathbb{E} [\mathcal {L}_\rho (x_{t+1}, y_{t+1},\lambda_t)] & \leq \mathcal {L}_\rho (x_t, y_{t+1},\lambda_t)
 - (\frac{\phi^H_{\min}}{\eta} + \frac{\phi^A_{\min}\rho}{2}- \frac{L+1}{2})\|x_t-x_{t+1}\|^2 + \frac{L^2}{2M}\|x_{t}-\tilde{x}\|^2. \label{eq:37}
\end{align}

By the step 11 of Algorithm \ref{alg:2}, we have
\begin{align}
\mathbb{E} [\mathcal {L}_\rho (x_{t+1}, y_{t+1},\lambda_{t+1})-\mathcal {L}_\rho (x_{t+1}, y_{t+1},\lambda_t)]
&=\frac{1}{\rho}\mathbb{E} \|\lambda_{t+1}-\lambda_t\|^2 \nonumber \\
\mathop{\leq}^{(i)} & \frac{5L^2}{\rho\phi^A_{\min}M} \|x_{t}-\tilde{x}\|^2
   + \frac{5L^2}{\rho\phi^A_{\min}M}\|x_{t-1}-\tilde{x}\|^2
   + \frac{5(\phi^H_{\max})^2}{\rho\phi^A_{\min}\eta^2}\|x_{t+1}-x_t\|^2 \nonumber \\
  & + \frac{5(L^2\eta^2+(\phi^H_{\max})^2)}{\rho\phi^A_{\min}\eta^2}\|x_{t}-x_{t-1}\|^2, \label{eq:38}
\end{align}
where the inequality $(i)$ holds by Lemma \ref{lem:7}.

Combining \eqref{eq:35}, \eqref{eq:37} and \eqref{eq:38}, we have
\begin{align}
\mathbb{E} [\mathcal {L}_\rho (x_{t+1}, y_{t+1},\lambda_{t+1})] \leq & \mathcal {L}_\rho (x_t, y_{t},\lambda_t) + \frac{(10+\phi^A_{\min}\rho)L^2}{2\rho\phi^A_{\min}M} \|x_{t}-\tilde{x}\|^2
   + \frac{5L^2}{\rho\phi^A_{\min}M}\|x_{t-1}-\tilde{x}\|^2 \nonumber \\
  &+ \frac{5(L^2\eta^2+(\phi^H_{\max})^2)}{\phi^A_{\min}\eta^2\rho}\|x_{t}-x_{t-1}\|^2 \nonumber \\
  &- (\frac{\phi^H_{\min}}{\eta} + \frac{\phi^A_{\min}\rho}{2}- \frac{L+1}{2}-\frac{5(\phi^H_{\max})^2}{\phi^A_{\min}\eta^2\rho})\|x_{t+1}-x_t\|^2. \label{eq:39}
\end{align}
Next, considering $\mathbb{E}\|x_{t+1}-\tilde{x}\|^2$, we have
\begin{align}
\mathbb{E}\|x_{t+1}-\tilde{x}\|^2 & = \mathbb{E}\|x_{t+1}-x_t+x_t-\tilde{x}\|^2 \nonumber \\
& = \mathbb{E}[\|x_{t+1}-x_t\|^2+ 2(x_{t+1}-x_t)^T(x_t-\tilde{x} )+\|x_t-\tilde{x}\|^2] \nonumber \\
& \mathop{\leq}^{(i)} \mathbb{E}[\|x_{t+1}-x_t\|^2+ 2(\frac{1}{2\beta}\|x_{t+1}-x_t\|^2
+ \frac{\beta }{2}\|x_t-\tilde{x}\|^2 )+\|x_t-\tilde{x}\|^2] \nonumber \\
& = (1+\frac{1}{\beta})\|x_{t+1}-x_t\|^2 + (1+\beta)\|x_t-\tilde{x}\|^2,  \label{eq:40}
\end{align}
where the inequality $(i)$ is due to the Cauchy-Schwarz inequality, and $\beta>0$.
Combining \eqref{eq:39} and \eqref{eq:40}, then, we have
\begin{align}
 &\mathbb{E} \big[ \mathcal {L}_\rho (x_{t+1}, y_{t+1},\lambda_{t+1})+\frac{5\big(L^2\eta^2+(\phi^H_{\max})^2\big)}{\phi^A_{\min}\eta^2 \rho}\|x_{t+1}-x_t\|^2
  + h^s_{t+1}(\|x_{t+1}-\tilde{x}\|^2 + \|x_{t}-\tilde{x}\|^2) \big] \nonumber \\
 & \leq \mathcal {L}_\rho (x_{t}, y_{t},\lambda_{t}) + \frac{5\big(L^2\eta^2+(\phi^H_{\max})^2\big)}{\phi^A_{\min}\eta^2\rho} \|x_t - x_{t-1}\|^2
  + \big((2+\beta)h^s_{t+1} + \frac{(10+\phi^A_{\min}\rho)L^2}{2\phi^A_{\min}M\rho} \big) (\|x_t-\tilde{x}\|^2+\|x_{t-1}-\tilde{x}\|^2)  \nonumber \\
 & \quad -\big[\frac{\phi^H_{\min}}{\eta} + \frac{\phi^A_{\min}\rho}{2}- \frac{L+1}{2} -\frac{5\big(L^2\eta^2+2(\phi^H_{\max})^2\big)}{\phi^A_{\min}\eta^2 \rho}-(1+\frac{1}{\beta})h_{t+1} \big]\mathbb{E}\|x_{t+1}-x_t\|^2 \nonumber \\
 & \quad - \big( (2+\beta)h^s_{t+1} + \frac{L^2}{2M}\big)\|x_{t-1}-\tilde{x}\|^2, \nonumber
\end{align}
where $h^s_{t+1}>0$.
By the definition of the sequence $\{(\Phi^{s}_{t})_{t=1}^m\}_{s=1}^S$ \ref{eq:88}, we have
\begin{align}
  \Phi^s_{t+1} \leq \Phi^s_{t}
 -\Gamma^s_t \mathbb{E}\|x^s_{t+1}-x^s_t\|^2 - \big( (2+\beta)h^s_{t+1} + \frac{L^2}{2M}\big)\|x^s_{t-1}-\tilde{x}^{s-1}\|^2. \label{eq:41}
\end{align}
Then using \eqref{eq:le8} and the properties of \emph{quadratic equation in one unknown},
we have $\Gamma^s_t >0, \ \forall t\in \{1,2,\cdots, m\}$.
Thus, we prove the first part.

Next, we will prove the second part.
Since $\lambda^{s+1}_0 = \lambda^s_m$ and $x^{s+1}_0=x^s_m=\tilde{x}^s$, we have
\begin{align}
 \mathbb{E}\|\lambda^{s+1}_0-\lambda^{s+1}_1\|^2 & = \mathbb{E}\|\lambda^{s}_m-\lambda^{s+1}_1\|^2 \nonumber \\
 & \leq \frac{1}{\phi^A_{\min}} \mathbb{E}\|A^T\lambda^{s}_m-A^T\lambda^{s+1}_1\|^2 \nonumber \\
 & \mathop{=}^{(i)} \frac{1}{\phi^A_{\min}} \mathbb{E}\|\hat{\nabla}f(x^s_{m-1})-\hat{\nabla}f(x^{s+1}_{0}) - \frac{H}{\eta}(x^s_{m-1}-x^s_m)
  -\frac{H}{\eta}(x^{s+1}_{0}-x^{s+1}_1)\|^2  \nonumber \\
 & \mathop{=}^{(ii)} \frac{1}{\phi^A_{\min}} \mathbb{E}\|\hat{\nabla}f(x^s_{m-1})
   - \nabla f(x^s_{m-1}) + \nabla f(x^s_{m-1})-\nabla f(x^{s}_{m})\nonumber \\
 & \quad -\frac{H}{\eta} (x^s_{m-1}-x^s_m) - \frac{H}{\eta}(x^{s+1}_{0}-x^{s+1}_1)\|^2  \nonumber \\
 & \leq \frac{5L^2}{\phi^A_{\min}M} \|x^s_{m-1}-\tilde{x}^{s-1}\|^2 + \frac{5(L^2\eta^2+(\phi^H_{\max})^2)}{\phi^A_{\min}\eta^2} \|x^s_{m-1}-x^s_{m}\|^2 \nonumber \\
 & \quad + \frac{5(\phi^H_{\max})^2}{\phi^A_{\min}\eta^2} \|x^{s+1}_{0}-x^{s+1}_1\|^2, \label{eq:42}
\end{align}
where the equality $(i)$ holds by \eqref{eq:30}, and the equality $(ii)$ holds by the following result:
\begin{align}
  \hat{\nabla}f(x^{s+1}_{0}) & = \frac{1}{M}\sum_{i_t\in \mathcal{I}_t} \big(\nabla f_{i_t}(x^{s+1}_{0}) - \nabla f_{i_t}(\tilde{x}^{s})\big)
   + \nabla f(\tilde{x}^{s}) \nonumber \\
  & = \frac{1}{M}\sum_{i_t\in \mathcal{I}_t} \big( \nabla f_{i_t}(x^{s}_{m}) - \nabla f_{i_t}(x^{s}_m) \big) + \nabla f(x^{s}_m) \nonumber \\
  & = \nabla f(x^{s}_m). \nonumber
\end{align}
By \eqref{eq:35}, we have
\begin{align}
 \mathcal {L}_\rho (x^{s+1}_0, y^{s+1}_{1},\lambda^{s+1}_0)&\leq \mathcal {L}_\rho (x^{s+1}_0, y^{s+1}_{0},\lambda^{s+1}_0)
  = \mathcal {L}_\rho (x^{s}_m, y^{s}_{m},\lambda^{s}_m).  \label{eq:43}
\end{align}
By \eqref{eq:37}, we have
\begin{align}
\mathbb{E} [\mathcal {L}_\rho (x^{s+1}_{1}, y^{s+1}_{1},\lambda^{s+1}_0)] &\leq \mathcal {L}_\rho (x^{s+1}_0, y^{s+1}_{1},\lambda^{s+1}_0) \nonumber \\
 & - \big(\frac{\phi^H_{\min}}{\eta} + \frac{\phi^A_{\min}\rho}{2} - \frac{L+1}{2}\big)\mathbb{E}\| x^{s+1}_{1}-x^{s+1}_0\|^2.  \label{eq:44}
\end{align}
By \eqref{eq:42}, we have
\begin{align}
  \mathbb{E} [\mathcal {L}_\rho (x^{s+1}_{1}, y^{s+1}_{1},\lambda^{s+1}_{1})-\mathcal {L}_\rho (x^{s+1}_{1}, y^{s+1}_{1},\lambda^{s+1}_0)]
 & = \frac{1}{\rho} \mathbb{E} \|\lambda^{s+1}_0 - \lambda^{s+1}_1\|  \nonumber \\
 & \leq \frac{5L^2}{\phi^A_{\min}\rho M} \|x^s_{m-1}-\tilde{x}^{s-1}\|^2 + \frac{5(L^2\eta^2+(\phi^H_{\max})^2)}{\phi^A_{\min}\eta^2\rho} \|x^s_{m-1}-x^s_{m}\|^2 \nonumber \\
 & \quad + \frac{5(\phi^H_{\max})^2}{\phi^A_{\min}\eta^2\rho} \|x^{s+1}_{0}-x^{s+1}_1\|^2.  \label{eq:45}
\end{align}
Combining \eqref{eq:43}, \eqref{eq:44} and \eqref{eq:45}, we have
\begin{align}
\mathbb{E} [\mathcal {L}_\rho (x^{s+1}_{1}, y^{s+1}_{1},\lambda^{s+1}_{1})] &\leq \frac{5L^2}{\phi^A_{\min}\rho M} \|x^s_{m-1}-\tilde{x}^{s-1}\|^2 + \frac{5(L^2\eta^2+(\phi^H_{\max})^2)}{\phi^A_{\min}\eta^2\rho} \|x^s_{m-1}-x^s_{m}\|^2 \nonumber \\
 & \quad - \big(\frac{\phi^H_{\min}}{\eta} + \frac{\phi^A_{\min}\rho}{2} - \frac{L+1}{2}-\frac{5(\phi^H_{\max})^2}{\phi^A_{\min}\eta^2\rho}\big)\mathbb{E}\| x^{s+1}_{1}-x^{s+1}_0\|^2. \nonumber
\end{align}
Using $h_1^{s+1}=\frac{10L^2}{\phi^A_{\min}\rho M}$, then, we have
\begin{align}
 &\mathbb{E} \big[\mathcal {L}_\rho (x^{s+1}_{1}, y^{s+1}_{1},\lambda^{s+1}_{1})
 + h^{s+1}_1\big( \|x^{s+1}_1-\tilde{x}^s\|^2 + \|x^{s+1}_0-\tilde{x}^s\|^2 \big)
 + \frac{5(L^2\eta^2+(\phi^H_{\max})^2)}{\phi^A_{\min}\eta^2\rho} \|x^{s+1}_{1}-x^{s+1}_{0}\|^2 \big] \nonumber \\
 & \leq \mathcal {L}_\rho (x^{s}_m, y^{s}_{m},\lambda^{s}_m) + \frac{10L^2}{\phi^A_{\min}\rho M} \big( \|x^s_{m}-\tilde{x}^{s-1}\|^2 + \|x^s_{m-1}-\tilde{x}^{s-1}\|^2 \big)
  + \frac{5(L^2\eta^2+(\phi^H_{\max})^2)}{\phi^A_{\min}\eta^2\rho} \|x^s_{m-1}-x^s_{m}\|^2 \nonumber \\
 & \quad - \big(\frac{\phi^H_{\min}}{\eta} + \frac{\phi^A_{\min}\rho}{2} - h_1^{s+1}-\frac{L+1}{2}-\frac{5(L^2\eta^2+2(\phi^H_{\max})^2)}{\phi^A_{\min}\eta^2\rho}\big)\mathbb{E}\|x^{s+1}_0 - x^{s+1}_{1}\|^2 - \frac{5L^2}{\phi^A_{\min}\rho M} \|x^s_{m-1}-\tilde{x}^{s-1}\|^2  \nonumber \\
 & \quad -\frac{10L^2}{\phi^A_{\min}\rho M} \|x^s_{m}-\tilde{x}^{s-1}\|^2. \nonumber
\end{align}
By the definition of the sequence $\big\{(\Phi^{s}_{t})_{t=1}^m\big\}_{s=1}^S$, we have
\begin{align}
 \Phi^{s+1}_{1} \leq \Phi^{s}_{m}
  - \Gamma^{s}_m\mathbb{E}\|x^{s+1}_0 - x^{s+1}_{1}\|^2 -\frac{5L^2}{\phi^A_{\min}\rho M} \|x^s_{m-1}-\tilde{x}^{s-1}\|^2. \label{eq:46}
\end{align}
Then using \eqref{eq:le8} and the properties of \emph{quadratic equation in one unknown},
we have $\Gamma^{s}_m>0,\ \forall s \geq 1$.
Finally, we prove that the sequence $\{(\Phi^{s}_{t})_{t=1}^m\}_{s=1}^S$
monotonically decreases.
\end{proof}

\subsection{Proof of Theorem 2}
\label{app:th-2}

\begin{proof}
By \eqref{eq:41} and \eqref{eq:46}, we have, for $s\in \{1,2,\cdots, S\}$ and $t\in\{1,2,\cdots,m\}$,
\begin{align}
 \Phi_{t+1}^{s}  \leq  \Phi_{t}^{s}
 -\Gamma_t^s\mathbb{E}\|x_{t+1}^{s}-x_t^{s}\|^2
- \big((2+\beta)h^{s}_{t+1}+\frac{L^2}{2M}\big)\|x^{s}_{t-1}-\tilde{x}^{s-1}\|^2,  \label{eq:50}
\end{align}
and
\begin{align}
\Phi^{s+1}_{1}  \leq  \Phi^{s}_{m}
  - \Gamma^{s}_m\mathbb{E}\|x^{s+1}_0 - x^{s+1}_{1}\|^2
  -\frac{5L^2}{\phi^A_{\min}\rho M} \|x^s_{m}-\tilde{x}^{s-1}\|^2.   \label{eq:51}
\end{align}
Summing \eqref{eq:50} and \eqref{eq:51} over $s\in \{1,2,\cdots, S\}$ and $t\in\{1,2,\cdots,m\}$, we have
\begin{align}
 \Phi_{m}^{S} - \Phi_{1}^{1}
 \leq -\gamma\sum_{s=1}^{S}\sum_{t=1}^{m}\mathbb{E}\|x_{t}^{s}-x_{t-1}^{s}\|^2
  - \omega \sum_{s=1}^{S}\sum_{t=1}^{m} \|x_{t-1}^{s}-\tilde{x}^{s-1}\|^2   \label{eq:52}
\end{align}
where $\gamma = \min_{(s,t)} \Gamma^s_t$,
and $\omega = \min_{(s,t)}\{(2+\beta)h^{s}_{t+1}+\frac{L^2}{2M},\frac{5L^2}{\phi^A_{\min}\rho M} \}$.
By Lemma \ref{lem:9}, there exists a constant $\Phi^*$
such that $\Phi^{s}_{t} \geq \Phi^*$.
By \eqref{eq:48} and \eqref{eq:52}, then, we have
\begin{align}
 (s^*,t^*) = \mathop{\arg\min}_{1 \leq s \leq S,\ 1 \leq t \leq m} \hat{\theta}^{s}_t \leq \frac{2}{\tau T} (\Phi^{1}_{1}- \Phi^*),
\end{align}
where $\tau=\min(\gamma,\omega)$, and $T=mS$.

By \eqref{eq:30}, we have
\begin{align}
 & \mathbb{E}\|A^T\lambda^{s}_{t+1}-\nabla f(x^{s}_{t+1})\|^2  \nonumber \\
 & = \mathbb{E}\|\hat{\nabla}f(x^{s}_{t})-\nabla f(x^{s}_{t+1})-\frac{H}{\eta }(x^{s}_{t}-x^{s}_{t+1})\|^2 \nonumber \\
 & = \mathbb{E}\|\hat{\nabla}f(x^{s}_{t})-\nabla f(x^{s}_{t}) +\nabla f(x^{s}_{t})- \nabla f(x^{s}_{t+1})
 - \frac{H}{\eta}(x^{s}_{t}-x^{s}_{t+1})\|^2 \nonumber \\
 & \leq \frac{3L^2}{M}\|x^{s}_{t}-\tilde{x}^{s-1}\|^2  + 3\big(L^2+\frac{(\phi^H_{\max})^2}{\eta^2}\big)\|x^{s}_{t}-x^{s}_{t+1}\|^2 \nonumber \\
 & \leq  3\big(L^2+\frac{(\phi^H_{\max})^2}{\eta^2}\big)\hat{\theta}^{s}_t. \label{eq:54}
\end{align}
By Lemma \ref{lem:7}, we have
\begin{align}
 \mathbb{E}\|Ax^{s}_{t+1}+By^{s}_{t+1}-c\|^2 & = \frac{1}{\rho^2}\|\lambda^{s}_{t+1}-\lambda^{s}_{t}\|^2  \nonumber \\
 & \leq \frac{5L^2}{\phi^A_{\min}\rho^2M} \mathbb{E} \|x_{t}^{s}-\tilde{x}^{s-1}\|^2
 + \frac{5L^2}{\phi^A_{\min}\rho^2M}\|x_{t-1}^{s}-\tilde{x}^{s-1}\|^2  \nonumber \\
 & + \frac{5(\phi^H_{\max})^2}{\phi^A_{\min}\rho^2\eta^2} \mathbb{E}\|x_{t+1}^{s}-x_t^{s}\|^2
 + \frac{5(L^2\eta^2+(\phi^H_{\max})^2)}{\phi^A_{\min}\rho^2\eta^2}\|x_{t}^{s}-x_{t-1}^{s}\|^2 \nonumber \\
 & \leq \frac{5(L^2\eta^2+\phi^2_{\max})}{\phi^A_{\min}\rho^2\eta^2} \hat{\theta}^s_{t} = \frac{\zeta}{\rho^2}  \hat{\theta}^s_{t}. \label{eq:55}
\end{align}
By the step 8 of Algorithm \ref{alg:2},
there exists a sub-gradient $\mu \in \partial g(y_{t+1}^{s})$ such that
\begin{align}
 \mathbb{E}\big[\mbox{dist}( B^T\lambda_{t+1}^{s}, \partial g(y_{t+1}^{s}))^2\big]  &\leq \|\mu-B^T\lambda_{t+1}^{s}\|^2 \nonumber \\
 & = \|B^T\lambda^s_{t}-\rho B^T(Ax^s_{t}+By_{t+1}^{s}-c)-B^T\lambda_{t+1}^{s}\|^2 \nonumber \\
 & = \|\rho B^TA(x^{s}_{t+1}-x^{s}_{t})\|^2 \nonumber \\
 & \leq \rho^2\|B\|^2\|A\|^2\|x^{s}_{t+1}-x^{s}_{t}\|^2 \nonumber \\
 & \leq \rho^2\|B\|^2\|A\|^2 \hat{\theta}^s_{t}.   \label{eq:56}
\end{align}
Finally, using the above bounds \eqref{eq:54}, \eqref{eq:55} and \eqref{eq:56}, and the definition \ref{def:1},
an $\epsilon$-stationary point of the problem \eqref{eq:1} holds in expectation.

\end{proof}

%% convergence analysis of nonconvex mini-batch SAGA-ADMM

\section{Convergence Analysis of Nonconvex Mini-batch SAGA-ADMM}

\subsection{Proof of Lemma 8}
\label{app:lem-11}
\begin{proof}
 Since $\psi_{t}=\frac{1}{n}\sum_{j=1}^n\nabla f_j(z^t_j)$, we have
  \begin{align}
\mathbb{E}\|\Delta_t\|^2 &= \mathbb{E}\|\frac{1}{M} \sum_{i_t\in \mathcal{I}_t} \big(\nabla f_{i_t}(x_t)-\nabla f_{i_t}(z^{t}_{i_t}) \big)+\psi_t-\nabla f(x_t)\|^2 \nonumber \\
  & \mathop{=}^{(i)}  \mathbb{E}\|\frac{1}{M} \sum_{i_t\in \mathcal{I}_t} \big(\nabla f_{i_t}(x_t)-\nabla f_{i_t}(z^{t}_{i_t}) \big)\|^2 - \|\nabla f(x_t)-\psi_t\|^2 \nonumber \\
  & \leq \frac{1}{M^2} \sum_{i_t\in \mathcal{I}_t} \mathbb{E}\|\nabla f_{i_t}(x_t)-\nabla f_{i_t}(z^{t}_{i_t})\|^2 \nonumber \\
  & =  \frac{1}{M^2} \sum_{i_t\in \mathcal{I}_t}\frac{1}{n} \sum_{i=1}^n \|\nabla f_{i_t}(x_t)-\nabla f_{i_t}(z^{t}_{i})\|^2 \nonumber \\
  & \mathop{\leq}^{(ii)} \frac{L^2}{nM} \sum_{i=1}^n \|x_t-z^t_i\|^2. \nonumber
 \end{align}
 where the equality (i) holds by the equality $\mathbb{E}(\xi-\mathbb{E}\xi)^2 = \mathbb{E}\xi^2-(\mathbb{E}\xi)^2$ for random variable $\xi$, and
 $\mathbb{E}[\nabla f_{i_t}(z^{t}_{i_t})]=\frac{1}{n}\sum_{j=1}^n\nabla f_j(z^t_j)=\psi_{t}$;
 the inequality (ii) holds by \eqref{eq:12}.
\end{proof}

\subsection{Proof of Lemma 10}
\label{app:lem-13}

\begin{proof}
By the step 5 of Algorithm \ref{alg:3}, we have
\begin{align}
 \mathcal {L}_\rho (x_t, y_{t+1},\lambda_t) \leq \mathcal {L}_\rho (x_t, y_{t},\lambda_t).  \label{eq:61}
\end{align}
By the optimal condition of step 7 in Algorithm \ref{alg:3}, we have
\begin{align}
0 & =(x_t-x_{t+1})^T\big[\hat{\nabla}f(x_t)+\rho A^T(Ax_{t+1}+By_{t+1}-c) -A^T\lambda_t- \frac{H}{\eta}(x_t-x_{t+1})\big] \nonumber \\
 & = (x_t-x_{t+1})^T\big[\hat{\nabla}f(x_t) - \nabla f(x_t) + \nabla f(x_t) -A^T\lambda_t - \frac{H}{\eta}(x_t-x_{t+1})
 +\rho A^T(Ax_{t+1}+By_{t+1}-c)\big] \nonumber \\
 & \mathop{\leq}^{(i)} f(x_t) - f(x_{t+1}) + \frac{L}{2}\|x_{t+1}-x_t\|^2 + (x_t -x_{t+1})^T\big(\hat{\nabla}f(x_t)-\nabla f(x_t)\big) -\frac{1}{\eta}\|x_{t+1}-x_t\|^2_H \nonumber \\
 &\quad - \lambda_t^T(Ax_{t+1}-A x_t) + \rho (Ax_t -Ax_{t+1})^T(Ax_{t+1}+By_{t+1}-c) \nonumber \\
 & \mathop{=}^{(ii)} f(x_t) - f(x_{t+1}) + \frac{L}{2}\|x_{t+1}-x_t\|^2 + (x_t -x_{t+1})^T\big(\hat{\nabla}f(x_t)-\nabla f(x_t)\big) -\frac{1}{\eta}\|x_{t+1}-x_t\|^2_H \nonumber \\
 & \quad - \lambda_t^T(A x_t+By_{t+1}-c) + \lambda_t^T(Ax_{t+1}+By_{t+1}-c)
  + \frac{\rho}{2}\|Ax_{t}+By_{t+1}-c\|^2 \nonumber \\
 & \quad - \frac{\rho}{2}\|Ax_{t+1}+By_{t+1}-c\|^2 -\frac{\rho}{2}\|Ax_t-Ax_{t+1}\|^2 \nonumber \\
 & = \mathcal {L}_\rho (x_t, y_{t+1},\lambda_t)- \mathcal {L}_\rho (x_{t+1}, y_{t+1},\lambda_t)
 + (x_t -x_{t+1})^T\big(\hat{\nabla}f(x_t) - \nabla f(x_t)\big) \nonumber \\
 & \quad + \frac{L}{2}\|x_{t+1}-x_t\|^2 -\frac{1}{\eta}\|x_{t+1}-x_t\|^2_H -\frac{\rho}{2}\|Ax_t-Ax_{t+1}\|^2 \nonumber \\
 & \mathop{\leq}^{(iii)} \mathcal {L}_\rho (x_t, y_{t+1},\lambda_t)- \mathcal {L}_\rho (x_{t+1}, y_{t+1},\lambda_t)
 + \frac{1}{2}\|\hat{\nabla}f(x_t) - \nabla f(x_t)\|^2 \nonumber \\
 & \quad - (\frac{\phi^H_{\min}}{\eta} + \frac{\phi^A_{\min}\rho}{2}-\frac{L+1}{2} )\|x_t-x_{t+1}\|^2,   \label{eq:62}
\end{align}
where the inequality (i) holds by \eqref{eq:13}; the equality (ii) holds by applying the equality
$(a-b)^T(b-c) = \frac{1}{2}(\|a-c\|^2-\|a-b\|^2-\|b-c\|^2)$ on the term $\rho (Ax_t -Ax_{t+1})^T(Ax_{t+1}+By_{t+1}-c)$;
the inequality (iii) holds by the Cauchy inequality.
Taking expectation conditioned on information $\mathcal{I}_t$ to \eqref{eq:62},  we have
\begin{align}
\mathbb{E} [\mathcal {L}_\rho (x_{t+1}, y_{t+1},\lambda_t)] \leq & \mathcal {L}_\rho (x_t, y_{t+1},\lambda_t)
+ \frac{L^2}{2 nM} \sum_{i=1}^n \mathbb{E} \|x_{t}-z^{t}_i\|^2- (\frac{\phi^H_{\min}}{\eta} + \frac{\phi^A_{\min}\rho}{2}-\frac{L+1}{2} )\|x_t-x_{t+1}\|^2.  \label{eq:63}
\end{align}

By the step 8 of Algorithm \ref{alg:3}, and taking expectation conditioned on $\mathcal{I}_t$, we have
\begin{align}
\mathbb{E} [\mathcal {L}_\rho (x_{t+1}, y_{t+1},\lambda_{t+1})-\mathcal {L}_\rho (x_{t+1}, y_{t+1},\lambda_t)]
&= \frac{1}{\rho}\mathbb{E} \|\lambda_t-\lambda_{t+1}\|^2 \nonumber \\
& \mathop{\leq}^{(i)}\frac{5L^2}{\rho\phi^A_{\min}M n} \sum_{i=1}^n \mathbb{E} \|x_{t}-z^{t}_i\|^2
   + \frac{5L^2}{\rho\phi^A_{\min}M n} \sum_{i=1}^n \|x_{t-1}-z^{t-1}_i\|^2\nonumber \\
& + \frac{5(\phi^H_{\max})^2}{\phi^A_{\min}\eta^2\rho}\mathbb{E}\|x_{t+1}-x_t\|^2
+ \frac{5(L^2\eta^2+(\phi^H_{\max})^2)}{\phi^A_{\min} \eta^2\rho}\|x_{t}-x_{t-1}\|^2,  \label{eq:64}
\end{align}
where the inequality $(i)$ holds by Lemma \ref{lem:12}.
Combining \eqref{eq:61}, \eqref{eq:63} and \eqref{eq:64}, we have
\begin{align}
 \mathbb{E} [\mathcal {L}_\rho (x_{t+1}, y_{t+1},\lambda_{t+1})] \leq  & \mathcal {L}_\rho (x_t, y_{t},\lambda_t) + \frac{10L^2+\phi^A_{\min}\rho L^2}{2\rho\phi^A_{\min}M n} \sum_{i=1}^n \mathbb{E} \|x_{t}-z^{t}_i\|^2 + \frac{5L^2}{\rho\phi^A_{\min}M n} \sum_{i=1}^n \|x_{t-1}-z^{t-1}_i\|^2 \nonumber \\
& + \frac{5(L^2\eta^2+(\phi^H_{\max})^2)}{\phi^A_{\min} \eta^2\rho}\|x_{t}-x_{t-1}\|^2
 - \big(\frac{\phi^H_{\min}}{\eta} + \frac{\phi^A_{\min}\rho}{2}-\frac{L+1}{2} - \frac{5(\phi^H_{\max})^2}{\phi^A_{\min}\eta^2\rho} \big)\mathbb{E}\|x_{t+1}-x_t\|^2.  \label{eq:65}
\end{align}

Next, we give an upper bound of $\frac{1}{n}\sum_{i=1}^n \mathbb{E}\|x_{t+1}-z^{t+1}_i\|^2$.
Using the step 9 in Algorithm \ref{alg:3}, we have
\begin{align}\label{eq:66}
\frac{1}{n}\sum_{i=1}^n \mathbb{E}\|x_{t+1}-z^{t+1}_i\|^2 &=\frac{1}{n}\sum_{i=1}^n
\big[ \frac{M}{n}\mathbb{E}\|x_{t+1}-x_{t+1}\|^2 + \frac{n-M}{n}\mathbb{E}\|x_{t+1} - z^{t}_i\|^2 \big] \nonumber \\
& = \frac{n-M}{n} \frac{1}{n}\sum_{i=1}^n\mathbb{E}\|x_{t+1} - z^{t}_i\|^2.
\end{align}
The term $\mathbb{E}\|x_{t+1} - z^{t}_i\|^2$ in \eqref{eq:66} can be bounded below:
\begin{align}
\mathbb{E}\|x_{t+1} - z^{t}_i\|^2 & = \mathbb{E}\|x_{t+1}-x_{t}+x_{t}-z^{t}_i\|^2 \nonumber \\
& = \mathbb{E}[\|x_{t+1}-x_{t}\|^2+ 2(x_{t+1}-x_{t})^T(x_{t}-z^{t}_i ) +\|x_{t-1}-z^{t}_i\|^2] \nonumber \\
& \mathop{\leq}^{(i)} \mathbb{E}[\|x_{t+1}-x_{t}\|^2 + 2(\frac{1}{2\beta}\mathbb{E}\|x_{t+1}-x_{t}\|^2 + \frac{\beta }{2}\|x_{t}-z^{t}_i\|^2 )
 +\|x_{t}-z^{t}_i\|^2] \nonumber \\
& = (1+\frac{1}{\beta})\mathbb{E}\|x_{t+1}-x_{t}\|^2 + (1+\beta)\|x_{t}-z^{t}_i\|^2,
\end{align}
where $\beta>0$, and the inequality $(i)$ is due to Cauchy-Schwarz inequality.
Thus, we have
\begin{align}
 &\frac{1}{n}\sum_{i=1}^n \mathbb{E}\|x_{t+1}-z^{t+1}_i\|^2 \leq \frac{n-M}{n}(1+\frac{1}{\beta})\mathbb{E}\|x_{t+1}-x_{t}\|^2
 + \frac{n-M}{n} (1+\beta) \frac{1}{n}\sum_{i=1}^n \|x_{t}-z^{t}_i\|^2 \label{eq:68}
\end{align}
Combining \eqref{eq:65} and \eqref{eq:68}, we have
\begin{align}
 &\mathbb{E} \big[\mathcal {L}_\rho (x_{t+1}, y_{t+1},\lambda_{t+1}) + \frac{5(L^2\eta^2+(\phi^H_{\max})^2)}{\phi^A_{\min} \eta^2\rho}\|x_{t+1}-x_t\|^2
  + \frac{\alpha_{t+1}}{n}\sum_{i=1}^n (\|x_{t+1}-z^{t+1}_i\|^2 + \|x_{t}-z^{t}_i\|^2) \big ] \nonumber \\
 & \leq \mathcal {L}_\rho (x_{t}, y_{t},\lambda_{t}) + \frac{5(L^2\eta^2+(\phi^H_{\max})^2)}{\phi^A_{\min} \eta^2\rho} \|x_t - x_{t-1}\|^2 \nonumber\\
 & \quad + \big(\frac{10L^2+\phi^A_{\min}\rho L^2}{2\rho\phi^A_{\min}M} + (\frac{2n-M}{n}+\frac{n-M}{n}\beta)\alpha_{t+1}\big) \frac{1}{n}\sum_{i=1}^n \big( \|x_{t}-z^{t}_i\|^2
  + \|x_{t-1}-z^{t-1}_i\|^2 \big) \nonumber\\
 & \quad -\big(\frac{\phi^H_{\min}}{\eta} + \frac{\phi^A_{\min}\rho}{2}-\frac{L+1}{2} - \frac{5L^2\eta^2+10(\phi^H_{\max})^2}{\phi^A_{\min} \eta^2\rho}
   -\frac{n-M}{n}(1+\frac{1}{\beta})\alpha_{t+1} \big) \|x_{t+1}-x_t\|^2 \nonumber\\
 & \quad -\big( \frac{L^2}{2M} + (\frac{2n-M}{n}+\frac{n-M}{n}\beta)\alpha_{t+1}\big)\frac{1}{n}\sum_{i=1}^n \|x_{t-1}-z^{t-1}_i\|^2.
\end{align}
By the definition of the sequence $\{\Theta_t\}_{t=1}^T$ \eqref{eq:89}, we have
\begin{align}
\Theta_{t+1} \leq  \Theta_{t} -\Gamma_{t}\|x_{t+1}-x_t\|^2
   -\big(\frac{L^2}{2M} + (\frac{2n-M}{n}+\frac{n-M}{n}\beta)\alpha_{t+1}\big)\frac{1}{n}\sum_{i=1}^n \|x_{t-1}-z^{t-1}_i\|^2. \label{eq:71}
\end{align}
Using \eqref{eq:le13} and the properties of \emph{quadratic equation in one unknown},
then we have $\Gamma_t >0$.
Finally, we prove that the sequence $\{\Theta_t\}_{t=1}^T$  monotonically decreases.
\end{proof}

\subsection{Proof of Theorem 3}
\label{app:th-3}

\begin{proof}
By \eqref{eq:71}, we have, for $t\in\{1,2,\cdots,T\}$
\begin{align}
 \Theta_{t+1} \leq  \Theta_{t} -\Gamma_{t}\|x_{t+1}-x_t\|^2
   -\big(\frac{L^2}{2M} + (\frac{2n-M}{n}+\frac{n-M}{n}\beta)\alpha_{t+1}\big)\frac{1}{n}\sum_{i=1}^n \|x_{t-1}-z^{t-1}_i\|^2.  \label{eq:74}
\end{align}
Summing \eqref{eq:74} over $t=1,2,\cdots,T$, we have
\begin{align}
 \Theta_{T}  &\leq \Theta_{1} - \gamma\sum_{t=1}^{T} \mathbb{E}\|x_{t+1}-x_t\|^2
   - \omega \sum_{t=1}^T\frac{1}{n}\sum_{i=1}^n\|x_{t-1}-z^{t-1}_i\|^2, \label{eq:75}
\end{align}
where $\gamma = \min_t \Gamma_t $ and $\omega=\min_t \frac{L^2}{2M} + (\frac{2n-M}{n}+\frac{n-M}{n}\beta)\alpha_{t+1}$.
By Lemma \ref{lem:14}, there exists a constant $\Theta^*$ such that $ \Theta_{t} \geq \Theta^*$
holds for $\forall t \geq 1$. By \eqref{eq:72} and \eqref{eq:75}, then, we have
\begin{align}
t^* = \mathop{\arg\min}_{2\leq t \leq T+1} \tilde{\theta}_{t} \leq \frac{2}{\tau T} (\Theta_1 - \Theta^*),
\end{align}
where $\tau=\min(\gamma,\omega)$.

Next, by the optimal condition
of step 7 in Algorithm \ref{alg:3}, we have
\begin{align}
 & \mathbb{E}\|A^T\lambda_{t+1}-\nabla f(x_{t+1})\|^2  \nonumber \\
 & = \mathbb{E}\|\hat{\nabla}f(x_{t}) - \nabla f(x_{t+1}) - \frac{H }{\eta}(x_t-x_{t+1})\|^2 \nonumber \\
 & = \mathbb{E}\|\hat{\nabla}f(x_{t})-\nabla f(x_{t}) +\nabla f(x_{t})- \nabla f(x_{t+1})
  - \frac{H}{ \eta} (x_t-x_{t+1})\|^2  \nonumber \\
 & \leq\frac{3L^2}{nM}\sum_{i=1}^n\|x_{t}-z^t_i\|^2  + 3\big(L^2+\frac{(\phi^H_{\max})^2}{\eta^2}\big)\|x_t-x_{t+1}\|^2 \nonumber \\
 & \leq 3\big(L^2+\frac{(\phi^H_{\max})^2}{\eta^2}\big)\tilde{\theta}_{t}. \label{eq:77}
\end{align}
By Lemma \ref{lem:12}, we have
\begin{align}
\mathbb{E}\|Ax_{t+1}+By_{t+1}-c\|^2 &= \frac{1}{\rho^2}\|\lambda_{t+1}-\lambda_t\|^2  \nonumber \\
 & \leq \frac{5L^2}{\phi^A_{\min} nM\rho^2} \sum_{i=1}^n \mathbb{E} \|x_{t}-z^{t}_i\|^2
   + \frac{5L^2}{\phi^A_{\min} nM\rho^2} \sum_{i=1}^n \|x_{t-1}-z^{t-1}_i\|^2 \nonumber \\
  & + \frac{5(\phi^H_{\max})^2}{\phi^A_{\min}\eta^2\rho^2}\mathbb{E}\|x_{t+1}-x_t\|^2
  + \frac{5(L^2\eta^2+(\phi^H_{\max})^2)}{\phi^A_{\min} \eta^2\rho^2}\|x_{t}-x_{t-1}\|^2 \nonumber \\
  & \leq \frac{5(L^2\eta^2+(\phi^H_{\max})^2)}{\phi^A_{\min} \eta^2\rho^2} \tilde{\theta}_t = \frac{\zeta}{\rho^2}\tilde{\theta}_t.  \label{eq:78}
\end{align}
By the step 5 of Algorithm \ref{alg:3},
there exists a subgradient $\mu \in \partial g(y_{t+1})$ such that
\begin{align}
 \mathbb{E}\big[\mbox{dist} (B^T\lambda_{t+1}, \partial g(y_{t+1}))^2\big]  & \leq \|\mu-B^T\lambda_{t+1}\|^2 \nonumber \\
 & = \|B^T\lambda_t-\rho B^T(Ax_t+By_{t+1}-c)-B^T\lambda_{t+1}\|^2 \nonumber \\
 & = \|\rho B^TA(x_{t+1}-x_{t})\|^2 \nonumber \\
 & \leq \rho^2\|B\|^2\|A\|^2\|x_{t+1}-x_t\|^2 \nonumber \\
 & \leq \rho^2\|B\|^2\|A\|^2 \tilde{\theta}_{t}.   \label{eq:79}
\end{align}
Finally, using the above bounds \eqref{eq:77}, \eqref{eq:78} and \eqref{eq:79}, and the definition \ref{def:1},
an $\epsilon$-stationary point of the problem \eqref{eq:1} holds in expectation.
\end{proof}

\end{appendices}

\end{onecolumn}
% that's all folks
\end{document}